\pdfoutput=1
\documentclass[11pt]{article}
\usepackage{amsmath,amssymb,amscd}
\usepackage{amsfonts}
\usepackage{microtype}
\usepackage{soul}
\usepackage{tikz-cd}
\usepackage{pb-diagram,pb-xy}
\usepackage[pagewise]{lineno}

\usepackage{hyperref}
\hypersetup{
    colorlinks=true,
    linkcolor=blue,
    filecolor=magenta,
    urlcolor=cyan,
}

\def\N{\mathbb{N}}

\def\R{\mathbb{R}}

\def\C{C^{\infty}(M, \R)}

\newtheorem{definition}{Definition}[section]
\newtheorem{lemma}[definition]{Lemma}
\newtheorem{proposition}[definition]{Proposition}
\newtheorem{theorem}[definition]{Theorem}
\newtheorem{remark}[definition]{Remark}

\newtheorem{remarks}[definition]{Remarks}
\newtheorem{examples}[definition]{Examples}
\newtheorem{example}[definition]{Example}
\newtheorem{comment}[definition]{Comment}
\newtheorem{some remarks}[definition]{Some remarks}
\newtheorem{some remarks a}[definition]{Some remarks on the rank of $\alpha$}
\newtheorem{remark on C}[definition]{A remark on the existence of $C$}
\newtheorem{dual dorfman}[definition]{Dual Dorfman connection}
\newenvironment{proof}{\noindent{\bf Proof.}}{\hfill $\blacklozenge$}

\def\lcf{\lbrack\! \lbrack}
\def\rcf{\rbrack\! \rbrack}

\def\lan{\langle \negthickspace \langle}
\def\ran{\rangle \negthickspace \rangle}

\begin{document}

\title{Courant-Dorfman algebras of differential operators and Dorfman connections of Courant algebroids}

\vspace{5mm}

\author{Panagiotis Batakidis and Fani Petalidou
\\
\\
\emph{Department of Mathematics}
\\
\emph{Aristotle University of Thessaloniki} \\
\emph{54124 Thessaloniki, Greece}}

\date{}
\maketitle

\vskip 20 mm

\begin{abstract}
\noindent
We construct an algebra and a complex of multidifferential operators on tensor products of a Courant algebroid $E$ with values in the endomorphism bundle of a smooth vector bundle $B$, predual of $E$, extending the standard complex of the Courant-Dorfman algebra of $E$. Also, we study Dorfman connections of $E$ on $B$, and show that the Cartan calculus, curvatures of induced connections and basic differential geometric identities of them make sense in this algebra.
\end{abstract}

\vspace{5mm} \noindent {\emph{Keywords: }}Courant algebroid; Courant-Dorfman algebra; standard cohomology; differential operators, Dorfman connections.

\vspace{3mm} \noindent MSC (2020): 53C05, 53D17, 58J10, 58A15, 16S32, 17B35, 17B55, 70G45.

\microtypesetup{protrusion=false}
\tableofcontents
\microtypesetup{protrusion=true}

\section{Introduction}

\vspace{1mm}
\noindent
Since their introduction by Liu, Weinstein and Xu \cite{lwx1}, Courant algebroids have enjoyed much attention due to their strong relation to Dirac structures \cite{c}, which are the natural geometric framework for the study of dynamical systems with constraints \cite{dorf}, graded differential geometry \cite{ro, ro2002}, higher structures \cite{Ur, BR} and generalized complex geometry \cite{Gualtieri}, along to their applications in string theory \cite{d, JV1} and T-duality \cite{s, SV}, among others. In the initial geometric formulation, a Courant algebroid is a vector bundle $E$ equipped with a fiberwise inner product, a bracket and an anchor map to the
tangent bundle of the base manifold, satisfying certain compatibility conditions. In the formulation of graded geometry, it is realized as a $NQ$-symplectic manifold of degree $2$, allowing one to study the differential geometry of such objects in a more concise way. For example, it has been proved \cite{ro2002, sev1999} that the \emph{standard cohomology} of the \emph{standard complex} of a Courant algebroid is isomorphic to the de Rham cohomology of the corresponding symplectic graded manifold.

\vspace{1mm}
\noindent
On the other hand, it is well known that the extension of the classical notion of (affine) connection from the tangent bundle $TM$  of a smooth manifold $M$ to a Lie algebroid \cite{rui-con} is widely used in differential geometry. Following the development of generalized geometry, the attention turned to connections of the generalized tangent bundle $TM\oplus T^*M$ of $M$ and more generally to \textit{Courant algebroid connections}. The latter are linear connections of a Courant algebroid to some vector bundle, first introduced by Alekseev and Xu \cite{Al-Xu}. Cueca and Mehta \cite{CM} used the version of the standard cochain complex of the algebraic definition of a Courant algebroid defined by Keller and Waldmann in \cite{kel-Wald} to develop a theory of linear Courant algebroid connections in a way that mirrors the classical theory of connections. As the left Leibniz property and skew-symmetry of the Courant bracket fail, such connections produce tensorial problems for the curvature \cite{Al-Xu, ABD, Gualt-branes}. One strategy to amend this is to consider variants, corrections and simplifications of the notion of curvature in order to produce a tensor. In the work of Aschieri et al. \cite{ABD} who study the graded geometric point of view of connections related to exact Courant algebroids, the authors correct the natural definitions for curvature and torsion by introducing $K$-curvature and $K$-torsion to get sections of vector bundles but do not have nice properties when restricted to Dirac subbundles.

\vspace{1mm}
\noindent
We take up the approach of working with a more general object, namely \emph{Dorfman connections of Courant algebroids $E$ on predual vector bundles $B$}, modifying and adapting their original definitions for dull algebroids \cite{mjl}. The reason to look at these connections is that, as will be explained in a subsequent paper, they are more suitable for the context of Manin pairs, i.e. pairs of Courant algebroids and Dirac subbundles. Another reason is that we do not opt for correcting or redefining the curvature in order to make it a tensor. For example,  we take into account the role of the contraction operator $i_f$ of degree $-2$ in our generalization of Roytenberg's complex, even though this produces additional terms in all computations. Ignoring this contraction corrects the tensorial anomalies in the naive definition of curvature and simplifies a few aspects in the theory of Courant algebroid connections \cite{CM, grtz-stn}.  Another major difference in our work is that Dorfman connections are not linear, so their behavior is completely different than that of a Courant algebroid connection. As a result, one first needs to develop the theory of such nonlinear connections, i.e., show their existence and confirm that a Cartan calculus holds similarly to the linear case. The second task is to construct a cohomology theory where the basic objects associated to a connection, its covariant derivative and curvature, make sense. Understanding these objects as differential operators on $E$, we fulfill this task by generalizing Roytenberg's Courant-Dorfman algebra $\mathcal{C}(\mathcal{E};\mathcal{R})$ \cite{ro2009} to an algebra $\mathfrak{D}(\mathcal{E}; \mathcal{R})$ of multidifferential operators and then equip it with a cochain complex structure. This is a new cohomology theory associated to Courant algebroids.

\vspace{1mm}
\noindent
The paper is structured as follows. In Section \ref{section-Courant cohom}, we first recall the basic notions related to Courant algebroids and Roytenberg's Courant-Dorfman algebra  and its cohomology (subsections \ref{subsec-cohomol-royt}, \ref{cohomology CD}). In  subsection \ref{cohomology-CD-Diff} we construct the algebra $\mathfrak{D}(\mathcal{E};\mathcal{R})$ of multidifferential operators on tensor products of $E$ (Definition \ref{cd op}). The main results of this part are gathered in the following Theorem.

\begin{theorem}
For a Courant algebroid $E$, \\
(i) the differential, contractions and Cartan calculus in $\mathcal{C}(\mathcal{E};\mathcal{R})$, extend to analogous operators and calculus in $\mathfrak{D}(\mathcal{E};\mathcal{R})$ and make it a cochain complex;\\
(ii) the projection map from $(\mathfrak{D}(\mathcal{E}; \mathcal{R}), d)$ to the subcomplex $\big((\mathcal{D}_{m, m-1}^p(\mathcal{E}; \mathcal{R}))_{p\in \N, m\in \N^\ast}, \partial_L \big)$ of the Loday cochain complex $\big(\mathcal{D}(\mathcal{E}; \mathcal{R}), \partial_L\big)$ of $E$ is a chain map.
\end{theorem}
The proof is covered in Propositions \ref{proposition-symbols}, \ref{diff cdop} and \ref{prop bp to loday}. In the rest of the paper it is shown that this construction is the natural environment for the development of a theory of Dorfman connections. Section \ref{section-Dorfman conn} is dedicated to Dorfman connections of Courant algebroids $E$ on predual vector bundles $B$. We start the discussion by modifying the definition of a predual vector bundle given in \cite{mjl} with Definition \ref{pre-dual} followed by a detailed description of the maps involved. The original definition of a Dorfman connection given in \cite{mjl} is also modified in Definition \ref{dorfman connection}. We then discuss the related notions of dual connection, curvature of a Dorfman connection, and prove the following.

\begin{theorem}
Given a Courant algebroid $E$ and a predual vector bundle $B$, \\
(i) the space of Dorfman connections of $E$ on $B$ is nonempty;\\
(ii) the set of Dorfman connections carries an affine structure;\\
(iii) each Dorfman connection defines a linear connection of the predual bundle $B$ on $E$; \\
(iv) the curvature of a Dorfman connection satisfies the Bianchi identity.
\end{theorem}
The proof is spread in Proposition \ref{existence Dorfman connection} (existence), Proposition \ref{affine space} (affine structure), Proposition \ref{prop - D - connection} (induced linear connection), and Proposition \ref{Bianchi} (Bianchi identity). In subsection \ref{section examples}, we provide concrete examples of Dorfman connections.

\vspace{3mm}
\noindent
\textbf{Notation:} Let $M$ be a smooth $n$-dimensional manifold, $TM$ and $T^\ast M$ its tangent and cotangent bundle, respectively,
and $\C$ the space of smooth functions on $M$. For each $p\in \mathbb{N}$, we denote by $\Omega^p$ the space
of smooth sections of $\bigwedge^p T^\ast M$. By convention, for $p<0$ we set $\Omega^p = \{0\}$, $\Omega^0 = \C$ and $\Omega = \oplus_{p\in \mathbb{Z}}\Omega^p$.
For an arbitrary smooth vector bundle $E \to M$, the space of smooth sections of $E$ is written as $\Gamma(E)$ and, for any $x\in M$, $E_x$ denotes the fibre of $E$ over $x$. If $E$ and $E'$ are two smooth vector bundles over $M$, we will frequently use the same letter to denote a vector bundle map $F: E \to E'$ and the induced $\C$-linear map $F : \Gamma(E) \to \Gamma(E')$ on spaces of smooth sections.

\vspace{3mm}
\noindent
\textbf{Acknowledgements} We would like to sincerely thank Iakovos Androulidakis and Ping Xu for their interest and valuable comments on the paper and Camille Laurent-Gengoux for discussions on earlier stages of this work. The paper is dedicated to the memory of our colleagues and friends Pantelis A. Damianou, Michel Marias and K. H. Mackenzie whose nobleness and mathematical work will continue to inspire us.

\section{Courant algebroids, Dirac structures, and their cohomologies}\label{section-Courant cohom}

\subsection{Courant algebroids and Dirac structures}
In \cite{lwx1} Liu, Weinstein and Xu introduced the notion of a \emph{Courant algebroid} in order to generalize the notion of the Drinfel'd double of a Lie bialgebra to the
notion of the \emph{double $A\oplus A^{\ast}$ of a Lie bialgebroid} $(A, A^\ast)$ (the last notion defined by Mackenzie and Xu in \cite{mx}).
This structure consists of a smooth vector bundle $E\to M$ together with a skew-symmetric bracket $[\cdot, \cdot]$ on the space $\Gamma(E)$ whose "Jacobi anomaly"
has an explicit expression in terms of a bundle map $E\to TM$ and a field of nondegenerate symmetric bilinear forms on $E$.
It leads furthermore to a Courant algebroid structure on $E = A\oplus A^\ast$. In his thesis \cite{ro}, Roytenberg reformulated the notion of Courant algebroid introducing a non skew-symmetric bracket $\lcf \cdot,  \cdot \rcf$ on $\Gamma(E)$
satisfying certain axioms and proved the equivalence of the two definitions. The bracket $\lcf \cdot,  \cdot \rcf$ is derived by adding
a symmetric part to the initial skew-symmetric bracket which is, in an
appropriate sense, a coboundary. When $E = TM \oplus T^\ast M$, the skew-symmetric bracket on $\Gamma(E)$ is given,
for sections $X + \zeta, Y+\eta \in \Gamma(TM\oplus T^\ast M)$, by the formula
\begin{equation*}\label{Courant bracket}
[X+\zeta, Y+\eta] = [X,Y] + \mathcal{L}_X\eta - \mathcal{L}_Y\zeta - \frac{1}{2}d(i_X\eta - i_Y\zeta).
\end{equation*}
This coincides with the bracket introduced by Ted Courant in \cite{c}
while the new bracket is
\begin{equation}  \label{Courant - Dorfman bracket}
\lcf X+\zeta, Y+\eta \rcf = [X,Y] + \mathcal{L}_X\eta - i_Yd\zeta.
\end{equation}
This coincides with the expression of the bracket considered by Dorfman in the context of complexes over Lie algebras,
in order to chara\-cterize Dirac structures \cite{dor}. The non skew-symmetric bracket $\lcf \cdot,  \cdot \rcf$ on $\Gamma(E)$ is named \emph{Courant-Dorfman bracket}.
For more details about the history of Courant algebroids, one may consult the insightful paper of Kosmann-Schwarzbach \cite{yks}. After the remarks of Uchino \cite{uch}, a Courant algebroid is defined as follows.

\begin{definition}\label{def Courant algebroid}
A \emph{Courant algebroid} over a smooth manifold $M$ is a constant rank vector bundle $E$ over $M$ equipped with: (i) a fiberwise nondegenerate symmetric bilinear form $\langle \cdot, \cdot \rangle$
on the bundle, (ii) a $\R$-bilinear bracket $\lcf \cdot, \cdot \rcf$ on $\Gamma(E)$, and (iii) a smooth vector bundle map $\rho : E \to TM$,
called the \emph{anchor map}\footnote{The map induced by $\rho : E \to TM$ on the spaces of smooth sections is also denoted by $\rho:\;\Gamma(E)\to\Gamma(TM)$.}, with the following properties:
\begin{enumerate}
\item
The bracket $\lcf \cdot, \cdot \rcf$ satisfies the Jacobi identity in Leibniz form
\begin{equation} \label{Jacobi-Courant}
\lcf e_1, \lcf e_2, e_3 \rcf \rcf = \lcf \lcf e_1, e_2 \rcf, e_3 \rcf + \lcf e_2, \lcf e_1, e_3 \rcf \rcf, \quad \textrm{for} \,\, \textrm{any} \,\, e_1,e_2,e_3 \in \Gamma(E).
\end{equation}
\item
The structures $\langle \cdot, \cdot \rangle$ and $\lcf \cdot, \cdot \rcf$ on $E$ are compatible in the sense that, for $e_1,e_2,e_3 \in \Gamma(E)$,
\begin{equation*} \label{compatibility <,> and [,]}
\rho(e_1)\langle e_2, e_3\rangle = \langle \lcf e_1, e_2 \rcf, e_3\rangle + \langle e_2 , \lcf e_1, e_3 \rcf\rangle.
\end{equation*}
\item
For $e_1, e_2 \in \Gamma(E)$,
\begin{equation*}
\lcf e_1, e_2 \rcf + \lcf e_2, e_1 \rcf = d_E\langle e_1, e_2\rangle,
\end{equation*}
where $d_E: \C \to \Gamma(E)$ is the map defined, for $f\in \C$ and $e\in \Gamma(E)$, by
\begin{equation*}\label{def-dE}
\langle d_Ef, e\rangle  = \langle e, d_Ef \rangle = \rho(e)(f),
\end{equation*}
i.e. $d_E = {g^{\flat}}^{-1}\circ \rho^\ast \circ d$, where $g^{\flat} : E \to E^\ast$ is the vector bundle map defined by $\langle \cdot, \cdot \rangle$.
\end{enumerate}
\end{definition}
From the above axioms, we get \cite{uch}:
\begin{enumerate}
\item[4.]
The anchor map $\rho : \big(\Gamma(E), \lcf \cdot, \cdot \rcf\big) \to \big(\Gamma(TM), [\cdot,\cdot]\big)$ is a homomorphism, i.e.
\begin{equation*}\label{morphisme rho}
\rho(\lcf e_1, e_2 \rcf) = [\rho(e_1), \rho(e_2)].
\end{equation*}
\item[5.]
The right Leibniz identity is satisfied:
\begin{linenomath*} \begin{equation*}\label{right-Leibniz}
\lcf e_1, fe_2 \rcf = f\lcf e_1, e_2 \rcf + \rho(e_1)(f)e_2.
\end{equation*} \end{linenomath*}
\end{enumerate}
Furthermore, the left Leibniz identity takes the form
\begin{linenomath*} \begin{equation*}\label{left-Leibniz}
\lcf fe_1, e_2 \rcf = f\lcf e_1, e_2 \rcf - \rho(e_2)(f)e_1 + \langle e_1,e_2\rangle d_Ef.
\end{equation*} \end{linenomath*}
For $e\in \Gamma(E)$ and $\alpha \in \Omega^1$, with a simple calculation one also derives \cite{uch} the identities:
\begin{enumerate}
\item[6.]
$\rho \circ \rho^\ast =0$,
\item[7.]
$\lcf e, \rho^\ast(\alpha) \rcf = \rho^\ast (\mathcal{L}_{\rho(e)}\alpha)$,
\item[8.]
$\lcf \rho^\ast(\alpha), e\rcf = -\rho^\ast (i_{\rho(e)}d\alpha)$,
\end{enumerate}
where $\mathcal{L}$ and $i$ denote, respectively, the classical Lie derivative and the contraction of differential forms with vector fields.

\vspace{0.5mm}
\noindent
Finally, the initial skewsymmetric bracket $[\cdot, \cdot]$ and the non skew-symmetric bracket $\lcf \cdot, \cdot \rcf$ on $\Gamma(E)$ are related by the formula
\begin{linenomath*} \begin{equation}\label{skew sym - non sym}
\lcf e_1, e_2 \rcf = [e_1, e_2] + \frac{1}{2}d_E \langle e_1, e_2\rangle.
\end{equation} \end{linenomath*}

\begin{remark}\label{courant as loday}
\emph{The first condition in Definition \ref{def Courant algebroid} makes $(\Gamma(E), \lcf \cdot, \cdot \rcf)$ a left Loday (initially called Leibniz) algebra \cite{lod, lod-pir}. For later use, note that Courant and Lie algebroids are examples of Loday algebroids\footnote{A \emph{Loday algebroid} structure on a vector bundle $E\to M$ is defined by a Loday bracket on the $\C$-module $\Gamma(E)$ which is a bidifferential operator of total order $\leq 1$ and for which the adjoint operator $\lcf e, \cdot\rcf : \Gamma(E)\to\Gamma(E)$ is a derivative endomorphism \cite{grab-k-pon}.} \cite[Example 5.4 \& Theorem 4.8 for $\alpha = 0$, respectively]{grab-k-pon}.}
\end{remark}

\vspace{1mm}
\noindent
Certain subbundles of Courant algebroids have a privileged role in Differential Geometry and Hamiltonian Mechanics \cite{hb-brief}. These are the \emph{Dirac subbundles} defined as follows.

\begin{definition}\label{def-Dirac}
Let $(E, \lcf \cdot, \cdot\rcf, \langle \cdot, \cdot \rangle, \rho)$ be a Courant algebroid over $M$. A \emph{Dirac structure} is a subbundle $L\subset E$
that is maximally isotropic with respect to $\langle \cdot, \cdot \rangle$ and the space $\Gamma(L)$ of smooth sections is
closed under the  bracket $\lcf \cdot, \cdot\rcf$, i.e. $\lcf \Gamma(L), \Gamma(L)\rcf \subseteq \Gamma(L)$.
\end{definition}

\begin{remarks}
\vspace{-2mm}
\noindent
\begin{enumerate}
\item
\emph{The first condition in Definition \ref{def-Dirac} is equivalent to $L = L^\bot$, where $L^\bot$ denotes the subbundle of $E$ that is orthogonal
to $L$ with respect to $\langle \cdot, \cdot \rangle$. It is also equivalent to the conditions $\langle \cdot, \cdot \rangle\vert_{L\times L} = 0$ and $\mathrm{rank}L = \frac{1}{2}\mathrm{rank}E$.
Thus, $L$ is both isotropic and coisotropic and so, in analogy with the terminology in Symplectic Geometry, it is also said to be a \emph{Lagrangian subbundle} of $E$.
Note that $E$ admits Lagrangian subbundles $L$ if and only if the pairing $\langle \cdot, \cdot \rangle$ has split signature $(\frac{1}{2}\mathrm{rank}E, \frac{1}{2}\mathrm{rank}E)$.}

\noindent
\emph{We can also prove that $L$ admits maximal isotropic complements. A choice of such a maximal isotropic subbundle $L'$ complement of $L$ relatively to $E$, i.e. $E = L \oplus L'$,
determines an isomorphism between the dual bundle $L^\ast$ of $L$ and the subbundle $L'$.}
\item
\emph{The second condition implies that the restrictions of the bracket $\lcf \cdot, \cdot\rcf$ and the anchor $\rho$ to
$\Gamma(L)$ turn $\big(L, \lcf \cdot, \cdot\rcf \vert_{\Gamma(L) \times \Gamma(L)}, \rho \vert_{\Gamma(L)} \big)$ into a Lie algebroid over $M$.
The generalized distribution $\rho(L)\subset TM$ is integrable and defines a generalized foliation of $M$.}
\end{enumerate}
\end{remarks}

\vspace{2mm}
\noindent
Below, we give some classical examples of Courant algebroids and Dirac structures.

\begin{example}[Courant algebroid over a point]\label{Cour over point}
{\rm{A Courant algebroid over a point, i.e. $M=\{p\}$, is just a quadratic Lie algebra $\mathfrak{g}$, that is a Lie algebra endowed with a nondegenerate symmetric bilinear form $\langle \cdot, \cdot \rangle$ invariant under the adjoint representation, namely, $\langle \mathrm{ad}_uv,w\rangle + \langle v, \mathrm{ad}_uw\rangle =0$, for all $u,v,w \in \mathfrak{g}$. If $\dim\mathfrak{g}$ is even, a Dirac subspace of $\mathfrak{g}$ is a Lagrangian Lie subalgebra $\mathfrak{l}$ of $\mathfrak{g}$.}}
\end{example}

\begin{example}[Standard Courant algebroid]\label{Standard Courant algebroid}
{\rm{Consider the vector bundle $E = TM \oplus T^\ast M$ over $M$ equipped with: (i) the nondegenerate symmetric fiberwise bilinear form $\langle \cdot, \cdot \rangle$ given, at each point $x\in M$ and for all $X + \zeta, Y+\eta \in T_xM \oplus T_x^\ast M$, by}}
\begin{linenomath*} \begin{equation*}\label{standard inner product}
\langle X+\zeta, Y+\eta \rangle = \langle \eta, X\rangle + \langle \zeta, Y\rangle,
\end{equation*} \end{linenomath*}
{\rm{(ii) the vector bundle map $\rho : TM \oplus T^\ast M \to TM$ projecting on the first summand, (iii) the Courant--Dorfman bracket \eqref{Courant - Dorfman bracket} on the space $\Gamma(E)$ of smooth sections of $E$, and (iv) the map $d_E : \C \to \Gamma(TM \oplus T^\ast M)$ defined by $d_Ef = (0, df)$. The above data define a Courant algebroid structure on $TM\oplus T^\ast M$ which is called}} standard.

\noindent
\rm{Some classical examples of Dirac structures $L$ of $TM\oplus T^\ast M$ are following.}
\begin{enumerate}
\item
{\rm{The graph $L=\mathrm{graph}\omega^\flat = \{(X, \omega^\flat(X)) \, / \, X\in TM )\}$ of the vector bundle map $\omega^\flat : TM \to T^\ast M$ defined by a (pre)-symplectic $2$-form $\omega$ on $M$ \cite{c}.}}
\item
{\rm{The graph $L= \mathrm{graph}\Pi^\# = \{(\Pi^\#(\eta), \eta) \, / \, \eta\in T^\ast M )\}$ of the vector bundle map $\Pi^\# : T^\ast M \to TM$ defined by a Poisson bivector field $\Pi$ on $M$ \cite{c}.}}
\item
{\rm{The subbundle $L=F \oplus F^0$, where $F\subseteq TM$ is an involutive regular distribution on $M$ and $F^0\subseteq T^\ast M$
its annihilator in $T^\ast M$ \cite{hb-brief}. Clearly, $F = F^\bot$  defines a regular foliation of $M$.
The involutivity of $F$ is equivalent to the integrability of $L$. Therefore, regular foliations of a manifold can be viewed as particular cases of Dirac structures.}}
\end{enumerate}
\end{example}

\begin{example}[The double of a Lie bialgebroid]\label{double lie bialgebroids} {\rm{Let $\big((A, [\cdot, \cdot]_A, a), (A^\ast, [\cdot,\cdot]_{A^\ast}, a_\ast)\big)$ be a Lie bialgebroid over a smooth manifold $M$. This is a pair of Lie algebroids
in duality verifying, for all $X,Y \in \Gamma(A)$, the compatibility condition
\begin{linenomath*} \begin{equation*}
d_{\ast}[X,Y]_A = [d_{\ast}X, Y]_A + [X, d_{\ast}Y]_A,
\end{equation*} \end{linenomath*}
where $d_{\ast}$ denotes the differential operator defined on $\Gamma(\bigwedge A)$ by the Lie algebroid structure of $A^\ast$ \cite{mx,mck,yks1}. The vector bundle $E = A\oplus A^\ast$ has a Courant algebroid structure defined by (i) the natural nondegenerate bilinear form}
\begin{linenomath*} \begin{equation*}
\langle X + \zeta, Y + \eta   \rangle = \langle \zeta, Y\rangle + \langle \eta, X\rangle, \quad \quad X + \zeta, \, Y + \eta \in \Gamma(E),
\end{equation*} \end{linenomath*}
(ii) the anchor map $\rho = a + a_{\ast}$, (iii) the operator $d_E = d_{\ast} + d$, where $d: \C \to \Gamma(A^\ast)$ is the usual differential operator associated to Lie algebroid structure of $A$, and (iv) the bracket
\begin{linenomath*} \begin{equation*}\label{br-bialgebroids}
\lcf X+\zeta, Y+\eta \rcf = ([X,Y]_A + \mathcal{L}_{\ast_\zeta}Y - i_{\eta}d_\ast X ) + ([\zeta, \eta]_{A^\ast} + \mathcal{L}_X\eta - i_Yd\zeta).
\end{equation*} \end{linenomath*}
The subbundles $A$ and $A^\ast$ are Dirac subbundles of $E=A\oplus A^\ast$.

\vspace{1mm}
\noindent
In the case where $(A, A^\ast)$ is a Lie-quasi, a quasi-Lie, or a proto-bialgebroid, the vector bundle $E = A \oplus A^\ast$ has again a Courant algebroid structure \cite{ro2002LMP, yks2}.}
\end{example}

\begin{example}\label{manin pairs complex manifold} \emph{Let $M$ be a smooth mani\-fold of dimension $2n$ equipped with an almost complex structure $J$, i.e., a vector bundle isomorphism $J : TM \to TM$ such that $J^2 = - Id$. Its complexified tangent bundle $TM \otimes \mathbb{C}$ is decomposed as
\begin{equation}\label{complex dec}
TM \otimes \mathbb{C} = T^{1,0}M \oplus T^{0,1}M,
\end{equation}
where $T^{1,0}M$ and $T^{0,1}M$ are, respectively, the $+i$ and $-i$ - eigenbundles of $J$. By duality, we have that
\begin{equation*}\label{dual complex dec}
T^{\ast} M \otimes \mathbb{C} = T_{1,0}^{\ast} M \oplus T^{\ast}_{0,1}M,
\end{equation*}
where $T_{1,0}^{\ast} M =(T^{0,1}M)^0$ and $T^{\ast}_{0,1}M = (T^{1,0}M)^0$ are, respectively, the annihilators of $T^{1,0}M$ and $T^{0,1}M$ with respect the usual pairing between vector bundles in duality.
The space $\Omega^c$ of complex differential forms on $M$ carries a bigrading and we write  $\Omega^c = \bigoplus^n_{p,q = 0}\Omega^{p,q}$, where $\Omega^{p,q}=\Gamma(\bigwedge^p T^\ast_{1,0}M\otimes\bigwedge^q T^\ast_{0,1}M)$. We say then that $J$ is a complex structure (or that $J$ is integrable)
if and only if its Nijenhuis torsion $N_J$\,\footnote{$N_J : TM \otimes TM \to TM$ is the vector bundle map associated to $J$ and given, for any pair $(X,Y)$ of vector fields on $M$, by the formula
\begin{equation*}
N_J(X,Y) = [JX,JY]-J[X,JY]-J[JX,Y]+J^2[X,Y] = [JX,JY]-J[X,JY]-J[JX,Y]-[X,Y].
\end{equation*}} vanishes on $M$. It is well known \cite[Theorem 2.8]{Kob-Nom} that the integrability of $J$ is equivalent to the involutivity of $T^{1,0}M$, of $T^{0,1}M$, and the fact that the differential operator $d$ on $\Omega^c$ may be written as a sum $d = \partial + \bar{\partial}$ of differential operators, where $\partial : \Omega^{p,q} \to \Omega^{p+1,q}$, and $\bar{\partial} : \Omega^{p,q} \to \Omega^{p,q+1}$. It is then easy to see that
\begin{equation*}
d^2 = 0 \quad \Leftrightarrow \quad  \partial^2 = 0, \quad \bar{\partial}^2 = 0, \quad \partial \circ \bar{\partial} + \bar{\partial} \circ \partial =0.
\end{equation*}}

\noindent
\emph{Consider the complex vector bundle $E=(TM \oplus T^\ast M)\otimes \mathbb{C}$ over a complex manifold $(M, J)$ endowed with the usual symmetric nondegenerate $\mathbb{C}$-bilinear form $\langle \cdot, \cdot \rangle$ on the fibers of $E$ with values in $\mathbb{C}$. Furthermore, equip $\Gamma(E)$ with the $\mathbb{C}$-bilinear bracket $\lcf \cdot, \cdot \rcf$ given by the complex version of \eqref{Courant - Dorfman bracket}, the natural anchor map $\rho : E\to TM\otimes \mathbb{C}$ and the map $d_E : C^\infty(M, \mathbb{C}) \to \Gamma(E)$, $d_E ={g^{\flat}}^{-1}\circ  \rho^\ast \circ d$. With the above data, $E$ is a complex Courant algebroid \cite{grtz-stn}. Following the third case of Examples \ref{Standard Courant algebroid}, the subbundles $L = T^{0,1}M \oplus (T^{0,1}M)^0 = T^{0,1}M \oplus T^\ast_{1,0}M$ and $L' = T^{1,0}M \oplus (T^{1,0}M)^0 = T^{1,0}M \oplus T^\ast_{0,1}M$ are (complex) Dirac subbundles of $E$, which are in fact transversal to each other. The pair $(L,L')$ is then a complex Lie bialgebroid \cite{lwx1}.}
\end{example}

\subsection{Cohomologies of a Courant algebroid}\label{subsec-cohomol-royt}
In \cite{ro2009} Roytenberg defined and studied the notion of \emph{Courant-Dorfman algebra} which is an algebraic analogue of Courant algebroids.
The relation is analogous to that of Lie-Rinehart algebras to Lie algebroids and of Poisson algebras to Poisson manifolds \cite{hbs}.

\vspace{1mm}
\noindent
A \emph{Courant-Dorfman algebra} consists of a commutative algebra $\mathcal{R}$, an $\mathcal{R}$-module $\mathcal{E}$ equipped with a pseudo-metric $\langle \cdot, \cdot \rangle$,
an $\mathcal{E}$-valued derivation $\partial$ of $\mathcal{R}$ and a Courant-Dorfman bracket $\lcf \cdot, \cdot \rcf$ satisfying compatibility conditions generalizing
those defining a Courant algebroid. Given a Courant-Dorfman algebra $(\mathcal{E}; \mathcal{R}, \lcf \cdot, \cdot \rcf, \langle \cdot, \cdot \rangle, \partial)$ (or just $(\mathcal{E};\mathcal{R})$),
a certain graded commutative $\mathcal{R}$-algebra $\mathcal{C}(\mathcal{E};\mathcal{R})$ endowed with a differential $d$ is then defined
and the resulting cochain complex $(\mathcal{C}(\mathcal{E};\mathcal{R}), d)$ is called the \emph{standard complex} of $(\mathcal{E};\mathcal{R})$. It is an analogue for a Courant-Dorfman algebra
of the de Rham complex of a Lie-Rinehart algebra. We first recall the notions of universal enveloping and convolution algebra,
and then present the structure of the Courant-Dorfman algebra of a Courant algebroid $(E, \lcf \cdot, \cdot\rcf, \langle \cdot, \cdot \rangle, \rho)$ following \cite{ro2009}.

\subsubsection{Universal enveloping and convolution algebra}\label{roy1}
Let $\mathbb{K}$ be a commutative ring containing $\frac{1}{2}$, $V$ and $W$ two $\mathbb{K}$-modules and $(\cdot, \cdot): V\otimes V \to W$ a symmetric bilinear form. Consider the graded $\mathbb{K}$-module $L = V[1]\oplus W[2]$  endowed with the bracket $[\cdot, \cdot] : L \times L \to L$ given, for any $v, v_1,v_2 \in V[1]$ and $w, w_1,w_2 \in W[2]$, by
\begin{equation*}
[v_1,v_2]=-(v_1,v_2), \quad [w_1,w_2] = 0 \quad \mathrm{and} \quad [v,w]=0.
\end{equation*}
Then $L$ becomes a graded Lie algebra over $\mathbb{K}$. Let $J$ be the homogeneous ideal of the tensor algebra $T(L)$ that is generated by elements of the form
\begin{linenomath*} \begin{equation*}
v_1 \otimes v_2 + v_2\otimes v_1 + (v_1,v_2), \quad \quad v\otimes w - w\otimes v, \quad \quad w_1\otimes w_2 - w_2 \otimes w_1.
\end{equation*} \end{linenomath*}
By definition, the universal enveloping algebra of $L$ is $U(L)=T(L)/J$ and $U(L)$ carries a natural filtration: Let $S(W) = \bigoplus_{k\geq 0}S^k(W)$ be the
symmetric algebra of $W$ and define, for $p\geq 0$,
\begin{linenomath*} \begin{equation*}
U(L)_{-p} = \bigoplus_{k=0}^{[\frac{p}{2}]} \big(V^{\otimes(p-2k)}\otimes S^kW \big)/R_{p},
\end{equation*} \end{linenomath*}
where $R_p$ is the submodule generated by elements of the form
\begin{eqnarray*}
v_1\otimes \ldots \otimes v_i \otimes v_{i+1} \otimes \ldots \otimes v_{p-2k}\otimes w_1\ldots w_k  \nonumber \\
+\, v_1\otimes \ldots \otimes v_{i+1} \otimes v_{i} \otimes \ldots \otimes v_{p-2k}\otimes w_1\ldots w_k& \nonumber \\
+ \, v_1\otimes \ldots \otimes \hat{v}_i \otimes \hat{v}_{i+1} \otimes \ldots \otimes v_{p-2k}\otimes (v_i, v_{i+1})w_1\ldots w_k,
\end{eqnarray*}
with $k=0,\ldots,[\frac{p}{2}]$ and $i=1,\ldots, p-2k-1$. Since $U(L)$ is a (graded cocommutative) $\mathbb{K}$-coalgebra and $\mathcal{R}$ is a $\mathbb{K}$-algebra, the space $\mathcal{A}=\mathcal{A}(V,W; \mathcal{R})=\mathrm{Hom}_{\mathbb{K}}(U(L),\mathcal{R})$
is an associative algebra equipped with the convolution product and is called the \emph{convolution algebra of} $U(L)$. Since $U(L)$ is non-positively graded
one gets that $\mathcal{A}$ is non-negatively graded. Each element of $\mathcal{A}^p=\mathrm{Hom}_{\mathbb{K}}(U(L)_{-p}, \mathcal{R})$ is determined by $([\frac{p}{2}]+1)$-tuple
\begin{linenomath*} \begin{equation*}
\omega = (\omega_0, \omega_1, \ldots, \omega_{[\frac{p}{2}]})
\end{equation*} \end{linenomath*}
of homomorphisms
\begin{linenomath*} \begin{equation*}\label{map - omega}
\omega_k : V^{\otimes^{p-2k}}\otimes W^{\otimes^k} \to \mathcal{R}.
\end{equation*} \end{linenomath*}
By construction, any $\omega_k$ is symmetric in the $W$-arguments and satisfies
\begin{linenomath*} \begin{equation} \label{omega_k - omega_k+1}
\begin{array}{l}
\omega_k(\ldots, v_i, v_{i+1}, \ldots ; \ldots )\, + \,\omega_k(\ldots,v_{i+1}, v_i, \ldots; \ldots) \\
\\
 =  - \, \omega_{k+1}(\ldots, \hat{v}_i, \hat{v}_{i+1}, \ldots; (v_i, v_{i+1}), \ldots),
\end{array}
\end{equation} \end{linenomath*}
for $v_i, v_j \in V$ and $i=1,\ldots, p-2k$. Hence, each $\omega_k$ defines a map
\begin{linenomath*} \begin{equation*}\label{map-omega-hom}
\omega_k : V^{\otimes^{p-2k}} \to \mathrm{Hom}_{\mathbb{K}}(S^kW, \mathcal{R}).
\end{equation*} \end{linenomath*}
Similarly, since $S(W[2])$ is a coalgebra (concentrated in even non positive degrees), the space $\mathrm{Hom}_{\mathbb{K}}(S(W[2]),\mathcal{R})$ is an algebra with the shuffle product given, for $H \in \mathrm{Hom}_{\mathbb{K}}(S^p(W[2]),\mathcal{R})$, $K \in \mathrm{Hom}_{\mathbb{K}}(S^q(W[2]),\mathcal{R})$,  by
\begin{linenomath*} \begin{equation} \label{multiplication in A - f}
(H\cdot K)(w_1, \ldots, w_{p+q}) = \sum\limits_{ \sigma \in sh(p,\,q)}H(w_{ \sigma(1)},\ldots, w_{ \sigma(p)})K(w_{ \sigma(p+1)},\ldots,w_{ \sigma(p+q)}).
\end{equation} \end{linenomath*}
Here and henceforth, $sh(p,q)$ is the set of $(p,q)$-shuffle permutations of $1,\ldots, p+q$, i.e., of permutations $ \sigma$ such that $ \sigma(1)< \ldots < \sigma(p)$ and $ \sigma(p+1)<\ldots < \sigma(p+q)$. This leads to the following formula for the product in $\mathcal{A}$:
\begin{linenomath*} \begin{equation} \label{multiplication in A}
\begin{array}{l}
(\omega \cdot \eta)_k (v_1,\ldots,v_{p+q-2k}) = \\
\\
\sum\limits_{\begin{array}{c} i+j=k \\ i\leq [\frac{p}{2}] \\ j\leq [\frac{q}{2}]\end{array}}\sum\limits_{\sigma \in sh(p-2i,\, q-2j)}(-1)^{|\sigma|}\omega_i(v_{\sigma(1)},\ldots, v_{\sigma(p-2i)})\eta_j (v_{\sigma(p-2i + 1)},\ldots, v_{\sigma(p+q-2k)}),
\end{array}
\end{equation} \end{linenomath*}
where $(-1)^{|\sigma|}$ is the signature of $\sigma$ and the product in each summand takes place in $\mathrm{Hom}_{\mathbb{K}}(S(W[2]), \mathcal{R})$. In particular, for $k=0$,
\begin{linenomath*} \begin{equation*}
\begin{array}{l}
(\omega \cdot \eta)_0 (v_1,\ldots,v_{p+q}) = \\
\\
\sum\limits_{\sigma \in sh(p,\, q)}(-1)^{|\sigma|}\omega_0(v_{\sigma(1)},\ldots, v_{\sigma(p)})\eta_0 (v_{\sigma(p + 1)},\ldots, v_{\sigma(p+q)}),
\end{array}
\end{equation*} \end{linenomath*}
where the multiplication in each summand takes place in $\mathcal{R}$.

\subsubsection{The cohomology of the Courant-Dorfman algebra of a Courant algebroid}\label{cohomology CD}
We now recall the standard complex of a Courant-Dorfman algebra associated to a Courant algebroid $(E, \lcf \cdot, \cdot\rcf, \langle \cdot, \cdot \rangle, \rho)$.
Let $\mathcal{R} = \C$, $\mathcal{E} = \Gamma(E^\ast)\cong \Gamma(E)$, $\partial = d_E$, and $\Omega^1 = \Gamma(T^\ast M)$. We thus have a metric $\mathcal{R}$-module
$(\mathcal{E}, \langle \cdot,\cdot \rangle)$ and a symmetric bilinear form $(\cdot,\cdot) : \mathcal{E} \times \mathcal{E} \to \Omega^1$
determined by $(\cdot,\cdot) = d\langle \cdot, \cdot \rangle$. The graded $\mathcal{R}$-module $L = \mathcal{E}[1]\oplus \Omega^1[2]$ is a graded Lie algebra over $\mathcal{R}$ with the nontrivial brackets determined by $-(\cdot,\cdot)$. The universal enveloping algebra $U(L)$ and corresponding convolution algebra
$\mathcal{A} = \mathcal{A}(\mathcal{E},\Omega^1; \mathcal{R}) = \mathrm{Hom}_{\mathbb{K}}(U(L), \mathcal{R})$ are also defined as above.
In particular, we have
\vspace{1mm}
\noindent
\begin{center}
$\mathcal{A}^0 =\C$, $\mathcal{A}^1 = \Gamma(E^\ast)$, $\mathcal{A}^2 = \Gamma\big((\bigwedge^2 E^\ast) \oplus TM\big)$, $\mathcal{A}^3 = \Gamma((\bigwedge^3E^\ast) \oplus (E^\ast \otimes TM))$.
\end{center}
\vspace{1mm}
\noindent
Set $\mathcal{C}^0(\mathcal{E}; \mathcal{R}) = \mathcal{R}$ and, for each $p>0$, let $\mathcal{C}^p(\mathcal{E}; \mathcal{R}) \subset \mathcal{A}^p$ be the submodule consisting of elements $\bar{\omega} = (\bar{\omega}_0,\bar{\omega}_1, \ldots, \bar{\omega}_{[\frac{p}{2}]})$
such that each
\begin{linenomath*} \begin{equation*}
\bar{\omega}_k : \mathcal{E}^{\otimes^{p-2k}} \otimes {\Omega^1}^{\otimes^k} \to \mathcal{R}
\end{equation*} \end{linenomath*}
satisfies the following two additional conditions:
\begin{enumerate}
\item
$\bar{\omega}_k : \mathcal{E}^{\otimes^{p-2k}} \to \mathrm{Hom}_{\mathbb{K}}(S^k \Omega^1, \mathcal{R})$ takes values in $\mathrm{Hom}_{\mathcal{R}}(S_{\mathcal{R}}^k \Omega^1, \mathcal{R}) \subset \mathrm{Hom}_{\mathbb{K}}(S^k \Omega^1, \mathcal{R})$, where $S_{\mathcal{R}}^k \Omega^1$
is the $\mathcal{R}$-module of the $k$-symmetric power of the $\mathcal{R}$-module $\Omega^1$ and $\mathrm{Hom}_{\mathcal{R}}(S_{\mathcal{R}}^k \Omega^1, \mathcal{R})$ is
the space of $\mathcal{R}$-linear maps  $S_{\mathcal{R}}^k \Omega^1\to\mathcal{R}$.
\item
$\bar{\omega}_k : \mathcal{E}^{\otimes^{p-2k}} \to \mathrm{Hom}_{\mathcal{R}}(S_{\mathcal{R}}^k \Omega^1, \mathcal{R})$ is $\mathcal{R}$-linear in the $(p-2k)$-th argument of $\mathcal{E}^{\otimes^{p-2k}}$.
\end{enumerate}

\vspace{1mm}
\noindent
It is shown in \cite{ro2009}, by induction and using \eqref{omega_k - omega_k+1}, that, for all $1\leq i < p-2k$ and $f \in \C$,
\begin{eqnarray*}\label{f-order-lin omega}
\lefteqn{ \bar{\omega}_k(e_1,\ldots, fe_i,\ldots ; \ldots)  = f\bar{\omega}_k(e_1,\ldots, e_i,\ldots ; \ldots)}\nonumber \\
  & & +\, \sum_{j=1}^{p-2k-i} (-1)^j\langle e_i, e_{i+j}\rangle \bar{\omega}_{k+1}(e_1,\ldots, \hat{e}_i,\ldots, \hat{e}_{i+j},\ldots ; df,\ldots).
\end{eqnarray*}
From the above one concludes that each term of the sequence $\bar{\omega} = (\bar{\omega}_0,\bar{\omega}_1, \ldots, \bar{\omega}_{[\frac{p}{2}]})\in \mathcal{C}^p(\mathcal{E}; \mathcal{R})$ is a first-order differential
operator in the first $p-2k-1$ arguments (see Definition \ref{def-dif op - symbol}) and $\mathcal{R}$-linear in the $(p-2k)$-th argument.

\vspace{2mm}
\noindent
It is easy to see that the space $\mathrm{Hom}_{\mathcal{R}}(S_{\mathcal{R}}^k \Omega^1, \mathcal{R})$ is identified with the $\mathcal{R}$-module of symmetric $k$-derivations of $\mathcal{R}$ \cite{ro2009}. In other words this is the space of symmetric maps on $\mathcal{R}^{\otimes^k}$ with values in $\mathcal{R}$ which are derivations in each argument. Hence, the
image $\bar{\omega}_k (e_1,\ldots, e_{p-2k})$ of $(e_1,\ldots, e_{p-2k}) \in \mathcal{E}^{\otimes^{p-2k}}$ can be
viewed as either a symmetric $k$-derivation of $\mathcal{R}$ whose value on $f_1,\ldots,f_k \in \mathcal{R}$ will be denoted by
\begin{linenomath*} \begin{equation*}
\omega_k(e_1,\ldots,e_{p-2k}; f_1,\ldots,f_k),
\end{equation*} \end{linenomath*}
or as a symmetric $\mathcal{R}$-multilinear function on $S_{\mathcal{R}}^k \Omega^1$ whose value on a $k$-tuple $(\alpha_1,\ldots,\alpha_k)$ of elements of $\Omega^1$ will be denoted by
\begin{linenomath*} \begin{equation*}
\bar{\omega}_k (e_1,\ldots,e_{p-2k}; \alpha_1,\ldots,\alpha_k).
\end{equation*} \end{linenomath*}
It is then obvious that
\begin{linenomath*} \begin{equation}\label{convention df - f}
\bar{\omega}_k (e_1,\ldots,e_{p-2k}; df_1,\ldots,df_k) = \omega_k (e_1,\ldots,e_{p-2k}; f_1,\ldots,f_k)
\end{equation} \end{linenomath*}
and so in the following we will interchange between the two realizations of elements of $\mathcal{C}^p(\mathcal{E}; \mathcal{R})$ without other notice.

\begin{definition}
The graded subalgebra $\big(\mathcal{C}(\mathcal{E};\mathcal{R}), \cdot \big)=\big((\mathcal{C}^p(\mathcal{E}; \mathcal{R}))_{p\geq0}, \cdot\big)$ of $\big(\mathcal{A}(\mathcal{E},\Omega^1; \mathcal{R}), \cdot\big)$ is called \emph{the Courant-Dorfman algebra}
of the Courant algebroid $(E, \lcf \cdot, \cdot\rcf, \langle \cdot, \cdot \rangle, \rho)$.
\end{definition}

\noindent
Define the map
\begin{linenomath*} \begin{equation} \label{d}
d: \mathcal{C}^\bullet(\mathcal{E}; \mathcal{R}) \to\mathcal{C}^{\bullet+1}(\mathcal{E}; \mathcal{R})
\end{equation} \end{linenomath*}
by setting, for all $\omega = (\omega_0,\omega_1, \ldots, \omega_{[\frac{p}{2}]}) \in \mathcal{C}^p(\mathcal{E}; \mathcal{R})$, $p\geq 0$,
\begin{linenomath*} \begin{equation*}
d\omega =\big((d\omega)_0,(d\omega)_1,\ldots,(d\omega)_{[\frac{p+1}{2}]}\big)\in \mathcal{C}^{p+1}(\mathcal{E}; \mathcal{R}),
\end{equation*} \end{linenomath*}
where, for any $k= 0, \ldots, [\frac{p+1}{2}]$,
\begin{eqnarray}\label{formule - d}
\lefteqn{(d\omega)_k(e_1,\ldots,e_{p+1-2k};f_1,\ldots,f_k) \, =} \nonumber \\
& &\sum\limits_{\mu=1}^k \omega_{k-1}(d_E f_{\mu}, e_1,\ldots,e_{p+1-2k}; f_1,\ldots, \hat{f}_{\mu},\ldots, f_k) \nonumber \\
& & + \, \sum\limits_{i=1}^{p+1-2k}(-1)^{i-1}\langle e_i, d_E(\omega_k(e_1,\ldots, \hat{e}_i,\ldots,e_{p+1-2k}; f_1,\ldots,f_k)) \rangle \nonumber \\
& & + \,\sum\limits_{i<j} (-1)^i \omega_k(e_1,\ldots, \hat{e}_i,\ldots, \hat{e}_j, \lcf e_i,e_j\rcf,e_{j+1},\ldots,e_{p+1-2k}; f_1,\ldots,f_k).
\end{eqnarray}

\vspace{1mm}
\noindent
On the other hand, $\overline{d\omega} =\big((\overline{d\omega})_0,(\overline{d\omega})_1,\ldots,(\overline{d\omega})_{[\frac{p+1}{2}]}\big)\in \mathcal{C}^{p+1}(\mathcal{E}; \mathcal{R})$, $p\geq 0$, is described at \cite[Corollary 4.9]{ro2009}. Recalling that $g^{\flat} : E \to E^\ast$ is the vector bundle map defined by $\langle \cdot, \cdot \rangle$, then for any $\alpha_1, \ldots, \alpha_k \in \Omega^1$, $\overline{d\omega}$ is given by the formula
\begin{eqnarray*}
\lefteqn{(\overline{d\omega})_k(e_1,\ldots,e_{p+1-2k};\alpha_1,\ldots,\alpha_k) \, =} \nonumber \\
& &\sum\limits_{\mu=1}^k \bar{\omega}_{k-1}({g^{\flat}}^{-1}(\rho^\ast(\alpha_{\mu})), e_1,\ldots,e_{p+1-2k}; \alpha_1,\ldots, \hat{\alpha}_{\mu},\ldots, \alpha_k) \nonumber \\
& & + \, \sum\limits_{i=1}^{p+1-2k}(-1)^{i-1}\rho(e_i)(\bar{\omega}_k(e_1,\ldots, \hat{e}_i,\ldots,e_{p+1-2k}; \alpha_1,\ldots,\alpha_k))  \nonumber \\
& & + \, \sum\limits_{i=1}^{p+1-2k}\sum\limits_{\mu=1}^k(-1)^{i} \bar{\omega}_k(e_1,\ldots, \hat{e}_i,\ldots,e_{p+1-2k}; \alpha_1,\ldots, i_{\rho(e_i)}d\alpha_{\mu}\ldots,\alpha_k) \nonumber \\
& & + \,\sum\limits_{i<j} (-1)^i \bar{\omega}_k(e_1,\ldots, \hat{e}_i,\ldots, \hat{e}_j, \lcf e_i,e_j\rcf,e_{j+1},\ldots,e_{p+1-2k}; \alpha_1,\ldots,\alpha_k) .
\end{eqnarray*}

\vspace{1mm}
\noindent
The next lemma shows that $d$ is a graded derivation on $\big(\mathcal{C}(\mathcal{E}; \mathcal{R}), \cdot \big)$.

\begin{lemma}
Let $\omega = (\omega_0,\ldots,\omega_{[\frac{p}{2}]}) \in \mathcal{C}^p(\mathcal{E}; \mathcal{R})$ and $\eta = (\eta_0,\ldots,\eta_{[\frac{q}{2}]}) \in \mathcal{C}^q(\mathcal{E}; \mathcal{R})$. The map \eqref{d} satisfies the Leibniz identity
\begin{linenomath*} \begin{equation*}\label{Leibniz - general}
d(\omega\cdot \eta) = d\omega \cdot \eta + (-1)^p \omega \cdot d\eta.
\end{equation*} \end{linenomath*}
\end{lemma}
\begin{proof}
A straightforward computation shows that
\begin{linenomath*} \begin{equation*}
d(\omega\cdot \eta) = \big((d(\omega\cdot \eta))_0, \ldots, (d(\omega\cdot \eta))_{[\frac{p+q+1}{2}]}\big),
\end{equation*} \end{linenomath*}
where, for any $k=0, \ldots, [\frac{p+q+1}{2}]$,
\begin{linenomath*} \begin{equation*}\label{Leibniz - k}
(d(\omega\cdot \eta))_k = (d\omega \cdot \eta)_k + (-1)^p (\omega \cdot d\eta)_k.
\end{equation*} \end{linenomath*}
\end{proof}

\begin{proposition}[\cite{ro2009}] \label{CDdiff}
The operator $d$ is a derivation of degree $+1$ of $\mathcal{C}(\mathcal{E};\mathcal{R})$ and squares to zero.
\end{proposition}

\vspace{1mm}
\noindent
The complex $(\mathcal{C}(\mathcal{E};\mathcal{R}), d)$ is named \emph{standard complex of} $(\mathcal{E};\mathcal{R})$ and its $p$-th cohomology group is denoted by $H^p(\mathcal{E};\mathcal{R})$.

\vspace{2mm}
\noindent
Define next two inner products in $\mathcal{C}(\mathcal{E}; \mathcal{R})$. For $\alpha \in \Omega^1$, consider the operator
$i_{\alpha} : \mathcal{C}(\mathcal{E}; \mathcal{R}) \to \mathcal{C}(\mathcal{E}; \mathcal{R})$ defined, for $\bar{\omega} = (\bar{\omega}_0, \ldots, \bar{\omega}_{[\frac{p}{2}]}) \in \mathcal{C}^p(\mathcal{E}; \mathcal{R})$, by
\begin{linenomath*} \begin{equation*}
i_{\alpha}\bar{\omega} = ((i_{\alpha}\bar{\omega})_0, \ldots, (i_{\alpha}\bar{\omega})_{[\frac{p-2}{2}]}),
\end{equation*} \end{linenomath*}
where
\begin{linenomath*} \begin{equation} \label{operator i-alpha}
(i_{\alpha}\bar{\omega})_k (e_1,\ldots,e_{p-2-2k};\alpha_1,\ldots,\alpha_k) = \bar{\omega}_{k+1}(e_1,\ldots,e_{p-2(k+1)};\alpha,\alpha_1,\ldots,\alpha_k).
\end{equation} \end{linenomath*}
For $f\in \mathcal{R}$, define similarly the operator $i_f$ so that
\begin{linenomath*} \begin{equation} \label{operator i-f}
\overline{i_f\omega} = i_{df}\bar{\omega}.
\end{equation} \end{linenomath*}
Let also $i_e : \mathcal{C}(\mathcal{E}; \mathcal{R}) \to \mathcal{C}(\mathcal{E}; \mathcal{R})$ act on $\omega = (\omega_0, \ldots, \omega_{[\frac{p}{2}]})\in \mathcal{C}^p(\mathcal{E}; \mathcal{R})$ returning the sequence
\begin{linenomath*} \begin{equation} \label{operator i-e - 1}
i_e\omega = ((i_e\omega)_0, \ldots, (i_e\omega)_{[\frac{p-1}{2}]}),
\end{equation} \end{linenomath*}
where
\begin{linenomath*} \begin{equation} \label{operator i-e}
(i_e\omega)_k(e_1, \ldots, e_{p-1-2k}; f_1,\ldots,f_k) = \omega_k(e,e_1, \ldots, e_{p-1-2k}; f_1,\ldots,f_k).
\end{equation} \end{linenomath*}

\begin{proposition}[\cite{ro2009}]\label{roy to lod-pir}
The map $\omega \to \omega_0$ is a chain map from the complex $(\mathcal{C}(\mathcal{E}; \mathcal{R}), d)$ to the Loday-Pirashvili complex $(\mathcal{C}_{LP}(\mathcal{E}; \mathcal{R}), d_{LP})$ of the Loday algebra $\mathcal{E}$ with coefficients in the symmetric $\mathcal{E}$-module $\mathcal{R}$.
\end{proposition}

\begin{lemma}
Let $f\in \mathcal{R}$, $e\in \mathcal{E}$, $\omega = (\omega_0, \ldots, \omega_{[\frac{p}{2}]})\in \mathcal{C}^p(\mathcal{E}; \mathcal{R})$ and $\eta = (\eta_0, \ldots, \eta_{[\frac{q}{2}]})\in \mathcal{C}^q(\mathcal{E}; \mathcal{R})$. The operators $i_f, i_e$ are derivations of degree $-2$ and $-1$, respectively, satisfying the Leibniz rules
\begin{linenomath*} \begin{equation*} \label{Leibniz - i - f}
i_f(\omega \cdot \eta) = (i_f \omega)\cdot \eta + \omega \cdot (i_f\eta),
\end{equation*} \end{linenomath*}
\begin{linenomath*} \begin{equation} \label{Leibniz - i - e}
i_e(\omega \cdot \eta) = (i_e \omega)\cdot \eta + (-1)^p\omega \cdot (i_e\eta).
\end{equation} \end{linenomath*}
\end{lemma}
\begin{proof}
Use formulas (\ref{multiplication in A}), (\ref{operator i-f}), (\ref{operator i-alpha}) and (\ref{operator i-e}) to show that as an element of $\mathcal{C}(\mathcal{E}; \mathcal{R})$, it is
$i_f(\omega \cdot \eta) = \big((i_f(\omega \cdot \eta))_0, \ldots, (i_f(\omega \cdot \eta))_{[\frac{p+q-2}{2}]}\big)$. A straightforward computation shows that
\begin{linenomath*} \begin{equation*}
(i_f(\omega \cdot \eta))_k = ((i_f \omega)\cdot \eta)_k + (\omega \cdot (i_f\eta))_k,
\end{equation*} \end{linenomath*}
for $k=0,\ldots,[\frac{p+q-2}{2}]$. Similarly, it is $i_e(\omega \cdot \eta) = \big((i_e(\omega \cdot \eta))_0, \ldots, (i_e(\omega \cdot \eta))_{[\frac{p+q-1}{2}]}\big)$, where
\begin{linenomath*} \begin{equation*}
(i_e(\omega \cdot \eta))_k = ((i_e \omega)\cdot \eta)_k + (-1)^p(\omega \cdot (i_e\eta))_k.
\end{equation*} \end{linenomath*}
\end{proof}

\vspace{2mm}
\noindent
Recall that the $\mathbb{K}$-module $L'=\mathcal{E}[1]\oplus \mathcal{R}[2]$ is a graded Lie algebra over $\mathbb{K}$ with the nontrivial brackets given by $-\langle \cdot,\cdot \rangle$. Let $\{\cdot,\cdot\}$
be the graded commutator on the space of graded endomorphisms of $\mathcal{C}(\mathcal{E}; \mathcal{R})$. If $P$ and $Q$ are two graded endomorphisms of degree $p$ and $q$, respectively, then the graded endomorphism
\begin{linenomath*} \begin{equation} \label{commutator in C}
\{P,Q\} = P\circ Q -(-1)^{pq}Q\circ P
\end{equation} \end{linenomath*}
is of degree $p+q$. Thus, there is a graded Lie algebra representation $i : L' \to \mathrm{End}(\mathcal{C}(\mathcal{E}; \mathcal{R}))$ of $L'$ in the space $\mathrm{End}(\mathcal{C}(\mathcal{E}; \mathcal{R}))$
of endomorphisms of $\mathcal{C}(\mathcal{E}; \mathcal{R})$ defined by the assignments $\mathcal{E}\ni e \mapsto i_e$ and $\mathcal{R} \ni f\mapsto i_f$. By construction, the  commutation relations are
\begin{linenomath*} \begin{equation} \label{br-derivations}
\{i_{e_1},i_{e_2}\} = i_{-\langle e_1,e_2\rangle} = - i_{\langle e_1,e_2\rangle},
\end{equation} \end{linenomath*}
\begin{linenomath*} \begin{equation} \label{br-derivations-0}
\{i_e,i_f\} = - \{i_f, i_e\} = 0 \quad \mathrm{and} \quad \{i_f, i_g\} = - \{i_g,i_f\} = 0.
\end{equation} \end{linenomath*}
Commuting the inner products $i_e$ and $i_f$ with the derivation $d$, we define the corresponding Lie derivatives:
\begin{linenomath*} \begin{equation} \label{Lie derivatives}
\mathcal{L}_e = \{i_e, d\} = i_e \circ d + d\circ i_e \quad \quad \mathrm{and} \quad \quad \mathcal{L}_f = \{i_f, d \} = i_f\circ d - d\circ i_f.
\end{equation} \end{linenomath*}

\begin{lemma}\label{Cartan - C}
The following Cartan's commutation relations hold in the space of graded endomorphisms of $\mathcal{C}(\mathcal{E}; \mathcal{R})$:
\begin{enumerate}
\item
$\mathcal{L}_f  = i_{d_Ef}$
\item
$\{\mathcal{L}_f, i_e \}  =  \mathcal{L}_f \circ i_e + i_e\circ \mathcal{L}_f = i_{d_E f}\circ i_e + i_e \circ i_{d_E f} \stackrel{(\ref{br-derivations})}{=} i_{-\langle d_E f, \, e \rangle}$
\item
$\{\mathcal{L}_e, i_f \}  =  \mathcal{L}_e \circ i_f - i_f\circ \mathcal{L}_e \stackrel{(\ref{Lie derivatives})}{=} - \{\mathcal{L}_f, i_e\} - \{\{i_f,i_e\},d\} \stackrel{(\ref{br-derivations-0})}{=}i_{\langle d_E f, \,e \rangle}$
\item
$\{\mathcal{L}_{e_1}, i_{e_2} \}  =  i_{\lcf e_1, \, e_2 \rcf}$
\item
$\{\mathcal{L}_f, \mathcal{L}_g \}  =  0$
\item
$\{\mathcal{L}_e, \mathcal{L}_f \}  = \{\mathcal{L}_e, i_{d_E f}\} = i_{\lcf e, d_E f\rcf} = i_{d_E \langle e, d_E f \rangle} = \mathcal{L}_{\langle e, \,d_E f \rangle}$
\item
$\{\mathcal{L}_f, \mathcal{L}_e \}  =  - \{\mathcal{L}_e, \mathcal{L}_f \} = - i_{d_E \langle e, \,d_E f \rangle} = - \mathcal{L}_{\langle e, \,d_E f \rangle}$
\item
$\{\mathcal{L}_{e_1}, \mathcal{L}_{e_2} \}  =  \mathcal{L}_{\lcf e_1, \,e_2 \rcf}$
\end{enumerate}
\end{lemma}
\begin{proof}
Direct computation.
\end{proof}

\vspace{1mm}
\noindent

\begin{comment}\label{Theta}
\emph{There is an alternative description of the complex $(\mathcal{C}(\mathcal{E}; \mathcal{R}), d)$ of the Courant-Dorfman algebra $(\mathcal{C}(\mathcal{E}; \mathcal{R}), \cdot)$ of a Courant algebroid $(E, \lcf \cdot, \cdot\rcf, \langle \cdot, \cdot \rangle, \rho)$ that uses graded geometry. More precisely Roytenberg \cite{ro2002}, following ideas of Ševera \cite{sev1999}, proved that there is a $1-1$ correspondence between vector bundles equipped with pseudo-Euclidean forms and symplectic $N$-manifolds of degree $2$. In particular, Courant algebroids $E$ are in $1-1$ correspondence with symplectic $NQ$-manifolds $(\mathcal{M}, \Omega)$ of degree $2$; $\Omega$ is the symplectic form, and the data $(\lcf \cdot, \cdot\rcf, \langle \cdot, \cdot \rangle, \rho)$ of $E$ are encoded in a function $\Theta$ of degree $3$ satisfying the structure equation}
\begin{linenomath*} \begin{equation*}
\{\Theta, \Theta\}=0.
\end{equation*} \end{linenomath*}
\emph{In the latter, the bracket $\{\cdot,\cdot\}$ is the Poisson bracket of degree $-2$ defined by $\Omega$ on the graded algebra $\mathcal{O}$ of polynomial functions. The operator $D = \{\Theta, \cdot \}$ defines a differential on $\mathcal{O}$ and it can be shown that the standard complex of the Courant algebroid is isomorphic to the Poisson dg algebra $(\mathcal{O}, D)$ \cite{ro2009}.
On the other hand, there is a natural Poisson structure on $\mathcal{C}(\mathcal{E}; \mathcal{R})$ for which $(\mathcal{C}(\mathcal{E}; \mathcal{R}),d)\cong (\mathcal{O},D)$ as Poisson dg algebras. Details and proofs can be found in \cite{ro2002, ro2009}.}
\end{comment}

\vspace{2mm}
\noindent
In the context of the Courant-Dorfman algebra $(\mathcal{E}; \mathcal{R}, \lcf \cdot, \cdot \rcf, \langle \cdot, \cdot \rangle, d_E)$ of a Courant
algebroid $(E, \lcf \cdot, \cdot\rcf, \langle \cdot, \cdot \rangle, \rho)$, a Dirac structure $L\subset E$ defines an $\mathcal{R}$-submodule $\mathcal{L}\subset \mathcal{E}$, $\mathcal{L}=\Gamma(L)$, which
is isotropic with respect to $\langle \cdot, \cdot \rangle$ and closed under $\lcf \cdot, \cdot \rcf$. It is then called a \emph{Dirac submodule}.

\begin{proposition}[\cite{ro2009}]
If $\mathcal{L}$ is a Dirac submodule of a Courant-Dorfman algebra $(\mathcal{E}; \mathcal{R})$, then $(\mathcal{L}; \mathcal{R})$ is a Lie-Rinehart algebra under the restrictions
of the bracket $\lcf \cdot, \cdot \rcf$ (or, equivalently, $[\cdot, \cdot]$) and the anchor map $\rho$.
\end{proposition}

\vspace{1mm}
\noindent
The corresponding graded commutative $\mathcal{R}$-algebra $\mathcal{C}(\mathcal{L}; \mathcal{R}) = \big(\mathcal{C}^p(\mathcal{L}; \mathcal{R})\big)_{p\geq 0}$ is the space
of $\mathcal{R}$-multilinear alternating forms on $L$, i.e., $\mathcal{C}^p(\mathcal{L}; \mathcal{R}) = \Gamma(\bigwedge^p L^\ast)$. It is endowed with the differential $d$ which is the Chevalley-Eilenberg differential for Lie algebroids.
The cohomology of the resulting cochain complex $(\mathcal{C}(\mathcal{L}; \mathcal{R}), d)$ is the Lie algebroid cohomology of $L$ with values in $\mathcal{R}$ and is denoted by $H^\bullet(L; \mathcal{R})$. For more details, see \cite{rine, hbs}.

\subsubsection{The cohomology of the Courant-Dorfman algebra of differential operators of a Courant algebroid}\label{cohomology-CD-Diff}
For the rest of the paper, we focus on a special type of multidifferential operators whose algebraic structure is described as follows. Let $E_i \to M$, $i=1,\ldots,k$, and $B\to M$ be smooth vector bundles of constant rank over a smooth manifold $M$, and $\Gamma(E_i)$, $i=1,\ldots,k$, $\Gamma(B)$ the corresponding spaces of smooth sections viewed as $\mathcal{R}$-modules.
\begin{definition}\label{def-dif op - symbol}
We say that an operator $D\in \Gamma\big(\mathrm{Hom}_{\mathbb{K}}(E_1\otimes \ldots \otimes E_k, \, B)\big)$ is a \emph{differential operator of order $s$ in the $i$-argument}, $i=1, \ldots, k$, if its $i$-symbol\,\footnote{Term which is introduced in \cite{dub-viol-mass}.}
$\sigma_i(D)(f)\in \Gamma\big(\mathrm{Hom}_{\mathbb{K}}(E_1\otimes \ldots \otimes E_k, \, B)\big)$ given, for any $(e_1, \ldots, e_i,\ldots,e_k)\in \Gamma(E_1\otimes \ldots \otimes E_i \otimes \ldots \otimes E_k)$
and $f\in \mathcal{R}$, by
\begin{linenomath*} \begin{equation*}
\sigma_i(D)(f)(e_1, \ldots, e_i,\ldots,e_k) : = D(e_1, \ldots, fe_i,\ldots,e_k) - f D(e_1, \ldots, e_i,\ldots,e_k),
\end{equation*} \end{linenomath*}
is a $(s-1)$-order differential operator on the $i$-argument of $\Gamma(E_1\otimes \ldots \otimes E_k)$.
\end{definition}

\vspace{1mm}
\noindent
For each $p \in\N$ and $m\in \N^{\ast}$, consider the submodule $\mathfrak{D}^p_{m, m-1}(\mathcal{E}; \mathcal{R})$ of $\mathcal{A}^p$ consisting of elements $\bar{\omega} = (\bar{\omega}_0,\bar{\omega}_1, \ldots, \bar{\omega}_{[\frac{p}{2}]})$ such that each
\begin{linenomath*} \begin{equation*}
\bar{\omega}_k : \mathcal{E}^{p-2k}\otimes S^k \Omega^1 \to \mathcal{R}
\end{equation*} \end{linenomath*}
satisfies the following two additional conditions:
\begin{itemize}
\item[-]
\emph{Condition 1:}\hspace{2mm} $\bar{\omega}_k : \mathcal{E}^{\otimes^{p-2k}} \to \mathrm{Hom}_{\mathbb{K}}(S^k \Omega^1, \mathcal{R})$ takes values in $\mathfrak{Diff}_m(S_{\mathcal{R}}^k \Omega^1, \mathcal{R})\subset \mathrm{Hom}_{\mathbb{K}}(S^k \Omega^1, \mathcal{R})$, where $S_{\mathcal{R}}^k \Omega^1$ is the $\mathcal{R}$-module of the $k$-symmetric power of the $\mathcal{R}$-module $\Omega^1$ and
\begin{eqnarray*}
\lefteqn{\mathfrak{Diff}_m(S_{\mathcal{R}}^k \Omega^1, \mathcal{R}) = \big\{D \in \mathrm{Hom}_{\mathbb{K}}(S^k \Omega^1, \mathcal{R}) \, / \, D \,\, \mathrm{is \, \, a \,\, differential}} \nonumber \\
& & \hspace{20mm} {\mathrm{operator\,\,of \,\, order \,\, at\,\, most}\,\,m\,\, \mathrm{in \,\, each \,\, entry}\big\}}.
\end{eqnarray*}
\item[-]
\emph{Condition 2:}\hspace{2mm} $\bar{\omega}_k : \mathcal{E}^{\otimes^{p-2k}} \to \mathfrak{Diff}_m(S_{\mathcal{R}}^k \Omega^1, \mathcal{R})$ is a differential operator of order at most $m$ on the first $p-2k-1$ arguments and of order at most $m-1$ on the $(p-2k)$-th argument of $\mathcal{E}^{\otimes^{p-2k}}$.
\end{itemize}
Note that for $p=1$, it is $\bar{\omega}=(\bar{\omega}_0)$ and the unique argument of the map $\bar{\omega}_0 : \mathcal{E} \to \mathcal{R}$ is considered as \emph{first} argument.
In particular, the spaces $\mathfrak{D}^1_{m}(\mathcal{E}; \mathcal{R})$, $m\in \N^\ast$, of differential operators on $\mathcal{E}$ of order at most $m$ fit in the  series
\begin{linenomath*} \begin{equation*}
\mathcal{C}^1(\mathcal{E}; \mathcal{R}) =\mathfrak{D}^1_0(\mathcal{E}; \mathcal{R}) \subset \mathfrak{D}^1_1(\mathcal{E}; \mathcal{R}) \subset \ldots \subset \mathfrak{D}^1_{m}(\mathcal{E}; \mathcal{R}) \subset \ldots,
\end{equation*} \end{linenomath*}
and, for any $p>1$ and $m\in \N^\ast$,
\begin{linenomath*} \begin{equation*}
\mathcal{C}^p(\mathcal{E}; \mathcal{R}) = \mathfrak{D}^p_{1,0}(\mathcal{E}; \mathcal{R}) \subset \ldots \subset \mathfrak{D}^p_{m,m-1}(\mathcal{E}; \mathcal{R}) \subset \ldots.
\end{equation*} \end{linenomath*}

\vspace{1mm}
\noindent
For any $(e_1,\ldots, e_{p-2k})\in \mathcal{E}^{\otimes^{p-2k}}$, $\bar{\omega}_k(e_1,\ldots,e_{p-2k})$ can be viewed as a symmetric map on $\Omega^{1^{\otimes k}}$ acting as a differential operator of order at most $m$  on each entry; its value on a $k$-tuple $(\alpha_1,\ldots,\alpha_k)$ is denoted by
\begin{linenomath*} \begin{equation*}
\bar{\omega}_k(e_1,\ldots,e_{p-2k}; \alpha_1,\ldots,\alpha_k).
\end{equation*} \end{linenomath*}
We adopt the convention presented in the previous paragraph and for any $k$-tuple $(df_1, \ldots, df_k)$ with $f_1,\ldots,f_k \in \mathcal{R}$, we write
\begin{linenomath*} \begin{equation*}
\bar{\omega}_k(e_1,\ldots,e_{p-2k}; df_1,\ldots,df_k) = \omega_k(e_1,\ldots,e_{p-2k}; f_1,\ldots,f_k).
\end{equation*} \end{linenomath*}
In the following, we switch between the two realizations of elements of $\mathfrak{D}^p_{m,m-1}(\mathcal{E}; \mathcal{R})$ without other notice.

\vspace{2mm}
\noindent
Set $\mathfrak{D}(\mathcal{E}; \mathcal{R}) = \big(\mathfrak{D}^p_{m, m-1}(\mathcal{E}; \mathcal{R})\big)_{p\in \N, \, m\in \N^\ast}$ and for any $\bar{\omega} = (\bar{\omega}_0, \ldots, \bar{\omega}_{[\frac{p}{2}]})\in
\mathfrak{D}^p_{m, m-1}(\mathcal{E}; \mathcal{R})$, let
\begin{itemize}
\item[-]
$\sigma_i(\bar{\omega}_k)$ be the symbol of $\bar{\omega}_k : \mathcal{E}^{p-2k}\otimes S^k \Omega^1 \to \mathcal{R}$ with respect to the $i$-th position, $i=1,\ldots,p-2k$, of $\mathcal{E}$-arguments:
For any $(e_1,\ldots, e_{p-2k})\in \mathcal{E}^{\otimes^{p-2k}}$ and $f\in \mathcal{R}$,
\begin{eqnarray*}
\sigma_i(\bar{\omega}_k)(f)(e_1,\ldots, e_i, \ldots, e_{p-2k}; \ldots) :& = &\bar{\omega}_k (e_1,\ldots, fe_i, \ldots, e_{p-2k}; \ldots)  \nonumber \\
& & \, -\, f \bar{\omega}_k (e_1,\ldots, e_i, \ldots, e_{p-2k}; \ldots).
\end{eqnarray*}
\item[-]
$s_i(\bar{\omega}_k)$ be the symbol of $\bar{\omega}_k : \mathcal{E}^{p-2k}\otimes S^k \Omega^1 \to \mathcal{R}$ with respect to the $i$-th position, $i=1,\ldots,k$, of $\Omega^1$-arguments: For any $\alpha_1, \ldots, \alpha_k \in \Omega^1$ and $f \in \mathcal{R}$,
\begin{eqnarray*}
s_i(\bar{\omega}_k)(f)(\ldots ; \, \alpha_1, \ldots, \alpha_i, \ldots, \alpha_k ) : & = &  \bar{\omega}_k(\ldots; \, \alpha_1, \ldots, f\alpha_i, \ldots, \alpha_k ) \nonumber \\
& & - \,f \bar{\omega}_k(\ldots ;\, \alpha_1, \ldots, \alpha_i, \ldots, \alpha_k ).
\end{eqnarray*}
\end{itemize}

\begin{proposition}\label{proposition-symbols}
The space $(\mathfrak{D}(\mathcal{E}; \mathcal{R}), \, \cdot)$ is a graded subalgebra of $(\mathcal{A}, \, \cdot)$.
More precisely, for any $\bar{\omega} \in \mathfrak{D}^p_{m, m-1}(\mathcal{E}; \mathcal{R})$ and $\bar{\eta} \in \mathfrak{D}^q_{n, n-1}(\mathcal{E}; \mathcal{R})$,
the differential operator $\bar{\omega} \cdot \bar{\eta}$ is an element in $\mathfrak{D}^{p+q}_{\mathrm{max}\{m,n\},\, \mathrm{max}\{m,n\}-1}(\mathcal{E}; \mathcal{R})$.
\end{proposition}
\begin{proof}
We must show that $\mathfrak{D}(\mathcal{E}; \mathcal{R})$ is closed under the shuffle multiplication in $\mathcal{A}$ defined by \eqref{multiplication in A - f} and \eqref{multiplication in A}.
Let $\bar{\omega} \in \mathfrak{D}^p_{m, m-1}(\mathcal{E}; \mathcal{R})$ and $\bar{\eta} \in \mathfrak{D}^q_{n, n-1}(\mathcal{E}; \mathcal{R})$. Clearly, $\bar{\omega} \cdot \bar{\eta}$ is a $(p+q)$-form. We will prove that $\bar{\omega} \cdot \bar{\eta}$ verifies conditions 1
and 2 defining $\mathfrak{D}(\mathcal{E}; \mathcal{R})$. Combining the formulas \eqref{multiplication in A - f} and \eqref{multiplication in A}, we obtain that,
for any $(e_1, \ldots, e_{p+q-2k}) \in \mathcal{E}^{\otimes ^{p+q-2k}}$ and $(\alpha_1, \ldots, \alpha_{k}) \in S_{\mathcal{R}}^{k} \Omega^1$,
\begin{eqnarray*}
\lefteqn{(\bar{\omega}\cdot \bar{\eta})_k(e_1, \ldots, e_{p+q-2k}; \, \alpha_1, \ldots, \alpha_{k}) = }\nonumber   \\
\\
& & \sum\limits_{\begin{array}{c} r+t=k,\\ r\leq [\frac{p}{2}],\\ t\leq [\frac{q}{2}]\end{array}}\sum\limits_{\varrho \in sh(p-2r,\, q-2t)} \sum\limits_{\tau \in sh(r,\,t) }
(-1)^{|\varrho|}
\bar{\omega}_r(e_{\varrho(1)},\ldots, e_{\varrho(p-2r)}; \, \alpha_{\tau(1)}, \ldots, \alpha_{\tau(r)} )\cdot \nonumber \\
& & \hspace{3cm} \bar{\eta}_t (e_{\varrho(p-2r + 1)},\ldots, e_{\varrho(p+q-2k)};\, \alpha_{\tau(r+1)},\ldots,\alpha_{\tau(r+t)}).
\end{eqnarray*}

\vspace{1mm}
\noindent
\emph{Condition 1:} Fixing $(e_1, \ldots, e_{p+q-2k})\in \mathcal{E}^{\otimes^{p+q-2k}}$, $(\alpha_1,\ldots, \alpha_k) \in S^k\Omega^1$, $(r,t)$ with $r+t=k$, and $\varrho\in sh(p-2r,\, q-2t)$,  observe that $(\bar{\omega} \cdot \bar{\eta})_k(e_1, \ldots, e_{p+q-2k}; \alpha_1, \ldots, \alpha_{k})$ is a sum over shuffle permutations $\tau \in sh(r,\,t)$. Hence, the argument in the $i$-position of $\Omega^1$-arguments of $(\bar{\omega} \cdot \bar{\eta})_k$ occurs either as an argument in the $\Omega^1$-arguments of $\bar{\omega}_r$ or in the $\Omega^1$-arguments of $\bar{\eta}_t$. A simple calculation shows that the value $s_i\big((\bar{\omega} \cdot \bar{\eta})_k\big)(f)(\ldots;\, \alpha_1,\ldots, \alpha_i,\ldots,\alpha_k)$ of the symbol $s_i\big((\bar{\omega} \cdot \bar{\eta})_k\big)(f)$ is a sum of terms of the form
\begin{eqnarray*}
\lefteqn{s_l(\bar{\omega}_r)(f)(\ldots ;\,\alpha_{\tau(1)},\ldots,\alpha_{\tau(r)} )\bar{\eta}_t(\ldots; \alpha_{\tau(r+1)},\ldots,\alpha_{\tau(r+t)} )\, + }\nonumber \\
& &  \bar{\omega}_r(\ldots; \alpha_{\tau(1)},\ldots,\alpha_{\tau(r)} )s_{l'}(\bar{\eta}_t)(f)(\ldots;\alpha_{\tau(r+1)},\ldots,\alpha_{\tau(r+t)}),
\end{eqnarray*}
where $l=\tau(i)$, if $1\leq \tau(i)\leq r$, and $l'=\tau(i)$, if $r+1\leq \tau(i)\leq r+t$, for $\tau \in sh(r,\,t)$. The  symbol $s_i\big((\bar{\omega} \cdot \bar{\eta})_k\big)(f)$ is thus a differential operator of order at most $\mathrm{max}\{m-1, n-1\}$ in the $\Omega^1$-arguments and $(\bar{\omega} \cdot \bar{\eta})_k$ is a differential operator of order at most $\mathrm{max}\{m-1,n-1\} + 1=\mathrm{max}\{m,n\}$ in the $\Omega^1$-arguments.

\vspace{1mm}
\noindent
\emph{Condition 2:} Similarly, fixing $(e_1, \ldots, e_{p+q-2k})\in \mathcal{E}^{\otimes^{p+q-2k}}$, $(\alpha_1,\ldots, \alpha_k) \in S^k\Omega^1$, $(r,t)$ with $r+t=k$, and $\tau\in sh(r, t)$, we have that $(\bar{\omega} \cdot \bar{\eta})_k(e_1, \ldots, e_{p+q-2k}; \alpha_1,\ldots, \alpha_{k})$ is a sum over shuffle permutations $\varrho\in sh(p-2r,\, q-2t)$.
Hence, the $i$-argument in the $\mathcal{E}$-entries of $(\bar{\omega} \cdot \bar{\eta})_k$ occurs either as an argument in the $\mathcal{E}$-entries of $\bar{\omega}_r$ or in the $\mathcal{E}$-entries of $\bar{\eta}_t$. Moreover, the last term of $(\bar{\omega} \cdot \bar{\eta})_k$ arises either as the last term of $\bar{\omega}_r$ or the last term of $\bar{\eta}_t$. A simple calculation shows again that the value $\sigma_i\big((\bar{\omega} \cdot \bar{\eta})_k\big )(f)(e_1,\ldots, e_i,\ldots,e_{p+q-2k}; \ldots)$ of the symbol $\sigma_i\big((\bar{\omega} \cdot \bar{\eta})_k\big )(f)$ is a sum of terms of the form
\begin{eqnarray*}
\lefteqn{(-1)^{|\varrho|}\big(\sigma_l (\bar{\omega}_r)(f)(e_{\varrho(1)},\ldots, e_{\varrho(p-2r)})\bar{\eta}_t(e_{\varrho(p-2r+1)},\ldots,e_{\varrho(p +q-2k)}) \, +} \nonumber \\
& & \bar{\omega}_r(e_{\varrho(1)},\ldots, e_{\varrho(p-2r)}) \sigma_{l'}(\bar{\eta}_t)(f)(e_{\varrho(p-2r+1)},\ldots,e_{\varrho(p +q-2k)})\big),
\end{eqnarray*}
where $l=\varrho(i)$, if $1\leq \varrho(i)\leq p-2r$, and $l'=\varrho(i)$, if $p-2r+1\leq \varrho(i)\leq p+q-2k$, for $\varrho \in sh(p-2r,\,q-2t)$. Hence the symbol $\sigma_i\big((\bar{\omega} \cdot \bar{\eta})_k\big )$ is a differential operator of order at most $\mathrm{max}\{m-1, n-1\}$. Thus, $(\bar{\omega} \cdot \bar{\eta})_k$ is a differential operator of order at most $\mathrm{max}\{m-1,n-1\}+1 = \mathrm{max}\{m,n\}$
in the first $p+q-2k-1$ $\mathcal{E}$-entries. The order of the symbol $\sigma_{p+q-2k}(\bar{\omega} \cdot \bar{\eta})_k$ is $\mathrm{max}\{m-2,\, n-2\}$, so the order of the operator $(\bar{\omega} \cdot \bar{\eta})_k$ in the ($p+q-2k$)-entry is $\mathrm{max}\{m-2,n-2\}+1=\mathrm{max}\{m, n\}-1$.
\end{proof}

\begin{definition}\label{cd op}
The graded subalgebra $\big(\mathfrak{D}(\mathcal{E}; \mathcal{R}), \cdot \big) =\big((\mathfrak{D}^p_{m, m-1}(\mathcal{E}; \mathcal{R}))_{p > 0, \, m\in \N^\ast}, \cdot \big)$ of the algebra $\big(\mathcal{A}, \cdot\big)$ is called the \emph{Courant-Dorfman algebra of differential operators} of the Courant algebroid $(E, \lcf \cdot, \cdot\rcf, \langle \cdot, \cdot \rangle, \rho)$.
\end{definition}

\vspace{1mm}
\noindent
We extend the map \eqref{d} to a map, also denoted by $d$,
\begin{linenomath*} \begin{equation}\label{def - d on D}
d: \mathfrak{D}^p_{m, m-1}(\mathcal{E}; \mathcal{R}) \to \mathfrak{D}^{p+1}_{m+1, m}(\mathcal{E}; \mathcal{R}),
\end{equation} \end{linenomath*}
defined by \eqref{formule - d}. One can directly check that

\begin{proposition}\label{diff cdop}
The operator $d: \mathfrak{D}^p_{m, m-1}(\mathcal{E}; \mathcal{R}) \to \mathfrak{D}^{p+1}_{m+1, m}(\mathcal{E}; \mathcal{R})$ is a derivation of degree $+1$ of $\big(\mathfrak{D}(\mathcal{E}; \mathcal{R}), \cdot \big)$ that squares to zero.
\end{proposition}

\noindent
The complex $\big(\mathfrak{D}(\mathcal{E}; \mathcal{R}), d\big)$ will be called \emph{complex of differential operators of $(\mathcal{E}; \mathcal{R})$}. Its $p$-cohomology group will be denoted by $\mathfrak{H}^p(\mathcal{E}; \mathcal{R})$.

\vspace{3mm}
\noindent
Let $(\mathcal{D}(\mathcal{E}; \mathcal{R}), \cdot)= \big((\mathcal{D}^p(\mathcal{E}; \mathcal{R}))_{p\in \N}, \cdot \big)$ be the graded \emph{shuffle algebra of $E$} \cite{grab-k-pon}. For any $p\in \N$, $\mathcal{D}^p(\mathcal{E}; \mathcal{R})$ consists of all multidifferential operators $\omega : \mathcal{E}^p \to \mathcal{R}$. The space $\mathcal{D}(\mathcal{E}; \mathcal{R})$ is endowed with the \emph{Loday operator $\partial_L$}, which is a degree 1 graded derivation \cite[Eqs. 62 \& 63]{grab-k-pon}. Since $E$ is a Courant algebroid, a particular case of Loday algebroid, $\partial_L^2 =0$ and $(\mathcal{D}(\mathcal{E}; \mathcal{R}), \partial_L)$ is a cochain complex whose corresponding cohomology is called \emph{Loday algebroid cohomology}. Let $\mathcal{D}_{m, m-1}^p(\mathcal{E}; \mathcal{R})$ be the subspace of $\mathcal{D}^p(\mathcal{E}; \mathcal{R})$ whose elements are multidifferential operators of order $\leq m$ with respect the first $p-1$ arguments and of order $\leq m-1$ in the last argument. Equipping it with the shuffle product\footnote{It is easy to see that $(\mathcal{D}_{m, m-1}^p(\mathcal{E}; \mathcal{R}))_{p\in \N, m\in \N^\ast}$ is closed with respect to the shuffle product.}, the corresponding algebra $\big((\mathcal{D}_{m, m-1}^p(\mathcal{E}; \mathcal{R}))_{p\in \N, m\in \N^\ast}, \cdot\big)$ is simply the subalgebra of $(\mathfrak{D}(\mathcal{E}; \mathcal{R}), \cdot)$ composed of the elements $\bar{\omega}=(\bar{\omega}_0)$, and the restriction of $\partial_L$ on $(\mathcal{D}_{m, m-1}^p(\mathcal{E}; \mathcal{R}))_{p\in \N, m\in \N^\ast}$ coincides with $d$ of Proposition \ref{diff cdop} when $k=0$ (it is enough to compare \cite[Eq. 63]{grab-k-pon} and \eqref{formule - d} for $k=0$). The following generalizes Proposition \ref{roy to lod-pir}.

\begin{proposition}\label{prop bp to loday}
Let $(\mathfrak{D}(\mathcal{E}; \mathcal{R}), d)$ be as above. For any  $\bar{\omega} = (\bar{\omega}_0, \ldots, \bar{\omega}_{[\frac{p}{2}]})\in \mathfrak{D}^p_{m, m-1}(\mathcal{E}; \mathcal{R})$, $p\in \N$, $m\in \N^\ast$, the map $\bar{\omega} \to \bar{\omega}_0$ from the complex $(\mathfrak{D}(\mathcal{E}; \mathcal{R}), d)$ to the subcomplex $\big((\mathcal{D}_{m, m-1}^p(\mathcal{E}; \mathcal{R}))_{p\in \N, m\in \N^\ast}, \partial_L\big)$ of $\big(\mathcal{D}(\mathcal{E}; \mathcal{R}), \partial_L\big)$ is a chain map.
\end{proposition}

\begin{remark}\label{rem - L - diff}
\emph{Recall that a Dirac structure $L\subset E$ defines a Dirac submodule $\mathcal{L}$ of the Courant-Dorfman algebra $(\mathcal{E}; \mathcal{R})$. By definition,  $\mathcal{C}^p(\mathcal{L}; \mathcal{R})= \Gamma(\bigwedge^p L^\ast)$. In our notation, this is the space  $\mathfrak{D}^p_{0}(\mathcal{L}; \mathcal{R})$ of (multi)differential operators on $L^{\otimes p}$ which are of $0$-order (i.e. linear) in each argument. The derivation \eqref{def - d on D} maps elements of $\mathfrak{D}^p_{0}(\mathcal{L}; \mathcal{R})$ to elements of $\mathfrak{D}^{p+1}_{0}(\mathcal{L}; \mathcal{R})$ and coincides with the Chevalley-Eilenberg differential on $\Gamma(\bigwedge^\bullet L^\ast)$. As a result, the cohomology of $\big((\mathfrak{D}^p_{0}(\mathcal{L}; \mathcal{R}))_{p\in \N}, \, d\big)$ is the Lie algebroid cohomology $H^\bullet (L; \mathcal{R})$.}
\end{remark}

\noindent
For any $f\in \mathcal{R}$ and $e\in \mathcal{E}$, the operators $i_f: \mathcal{C}(\mathcal{E}; \mathcal{R}) \to \mathcal{C}(\mathcal{E}; \mathcal{R})$
and $i_e : \mathcal{C}(\mathcal{E}; \mathcal{R}) \to \mathcal{C}(\mathcal{E}; \mathcal{R})$ given respectively by
\eqref{operator i-f} and \eqref{operator i-e - 1}, extend in a similar way to operators $i_f : \mathfrak{D}(\mathcal{E}; \mathcal{R}) \to \mathfrak{D}[-2](\mathcal{E}; \mathcal{R})$
and $i_e: \mathfrak{D}(\mathcal{E}; \mathcal{R}) \to \mathfrak{D}[-1](\mathcal{E}; \mathcal{R})$, respectively. More precisely, they preserve the order of differential operators, that is
\begin{linenomath*} \begin{equation} \label{operator i-f-e-D}
i_f : \mathfrak{D}^p_{m, m-1}(\mathcal{E}; \mathcal{R}) \to \mathfrak{D}^{p-2}_{m, m-1}(\mathcal{E}; \mathcal{R}) \quad \mathrm{and} \quad i_e : \mathfrak{D}^p_{m, m-1}(\mathcal{E}; \mathcal{R}) \to
\mathfrak{D}^{p-1}_{m, m-1}(\mathcal{E}; \mathcal{R}).
\end{equation} \end{linenomath*}
The commutator $\{\cdot, \cdot\}$ on the space of graded endomorphisms of $\mathfrak{D}(\mathcal{E}; \mathcal{R})$ is similarly defined, and the relations \eqref{br-derivations} - \eqref{Lie derivatives} together with Cartan's formulas of Lemma \ref{Cartan - C} hold as well.

\vspace{1mm}
\noindent
\begin{example}
\emph{Let $(E, \lcf \cdot, \cdot\rcf, \langle \cdot, \cdot \rangle, \rho)$ be a Courant algebroid over a smooth manifold $M$ and $D : \Gamma(E)\times \Gamma(E) \to \Gamma(E)$ a linear $E$-connection on $E$ as defined in \cite{Andr-Scandalis, CM, {Gualt-branes}}. The curvature $R^D(e_1, e_2) = [D_{e_1}, D_{e_2}] - D_{\lcf e_1, e_2 \rcf}$ is a $1$-order differential operator on the first argument and $0$-order differential operator on the second argument with values in $\Gamma(\mathrm{End}(E))$. Consequently, the operator $A : \mathcal{E}^{\otimes 4} \to \mathcal{R}$ defined by
\[A(e_1,e_2,e_3,e_4) = \langle R^D(e_1, e_2)e_3 + c.p., \, e_4\rangle \]
is a $1$-order differential operator on the first three entries and $0$-order on the last entry, thus $A\in \mathfrak{D}^4_{1,0}(\mathcal{E}; \mathcal{R})$. Likewise, the naive torsion operator $T^D$ of $D$ defined by
\[T^D(e_1,e_2) = D_{e_1}e_2 - D_{e_2}e_1 - \lcf e_1, e_2\rcf, \quad \quad e_1, e_2 \in \Gamma(E),\]
is a $1$-order differential operator on the first argument and $0$-order on the second argument with values in $\Gamma(E)$. By coupling it with smooth sections of $E$ and taking the cyclic permutation, we produce the operator
\[\mathcal{T}^D(e_1, e_2, e_3) = \langle T^D(e_1, e_2), e_3 \rangle + c.p.,\]
which is an element of $\mathfrak{D}^3_1(\mathcal{E}; \mathcal{R})$.}

\noindent
\emph{The operators $R^D$ and $T^D$ are replaced or corrected to $\C$-linear substitutes in \cite{ABD, {Gualt-branes}}.}
\end{example}

\section{Dorfman connections}\label{section-Dorfman conn}

Dorfman connections were initially introduced for dull algebroids\footnote{A \emph{dull algebroid} is a vector bundle $Q$ over a smooth manifold $M$ endowed with an anchor
map $\rho_Q : Q \to TM$ and a bracket $[\cdot, \cdot]_Q$ on $\Gamma(Q)$ such that, for all $q_1, q_2 \in \Gamma(Q)$ and $f_1,f_2 \in \C$,
$\rho_Q[q_1,q_2]_Q = [\rho_Q(q_1), \rho_Q(q_2)]$ and satisfying the Leibniz identity in both terms:
\begin{linenomath*} \begin{equation*}
[f_1q_1, f_2q_2]_Q = f_1f_2[q_1,q_2]_Q +f_1\rho_Q(q_1)(f_2)q_2 - f_2\rho_Q(q_2)(f_1)q_1.
\end{equation*} \end{linenomath*} } in \cite{mjl} by Jotz Lean in order to study the standard Courant algebroid $TE \oplus T^\ast E$ over a vector bundle $E$. In particular, the author establishes a one-to-one correspondence between linear splittings of $TE\oplus T^\ast E$ and a special kind of Dorfman connections on $E\oplus T^\ast M$. This result generalizes in the context of Courant algebroids the well known result of Dieudonn\'e \cite{dieud} that a linear $TM$-connection on a smooth vector bundle $E$ defines a splitting of $TE$ in horizontal and vertical subbundles. In the following, we modify the above notion for Courant algebroids and develop its basic theory.

\subsection{Predual vector bundle}\label{predual}
We first discuss the notion of \emph{predual vector bundle $B$} for a Courant algebroid $E$  by adapting the corresponding notion for a dull algebroid given in \cite[Definition 3.1]{mjl} \footnote{The difference with the definition in \cite{mjl} is that we require the map $d_B : \C \to \Gamma(B)$ to be a derivation.}.
\begin{definition}\label{pre-dual}
Let $(E, \lcf \cdot, \cdot\rcf, \langle \cdot, \cdot \rangle, \rho)$ be a Courant algebroid over $M$, $B\rightarrow M$ a smooth vector bundle of constant rank,
$\lan \cdot, \cdot \ran : \Gamma(E) \times \Gamma(B) \to \C$ a fiberwise bilinear pairing, and $d_B : \C \to \Gamma(B)$, a derivation of $\C$ into the $\C$-bimodule $\Gamma(B)$, such that, for any $e\in \Gamma(E)$ and $f\in \C$,
\begin{linenomath*} \begin{equation} \label{def-dB}
\lan e, d_Bf \ran = \rho(e)(f).
\end{equation} \end{linenomath*}
The triple $(B, d_B, \lan \cdot, \cdot \ran)$ is called \emph{predual} of $E$, and $E$ and $B$ are said to be \emph{paired by} $\lan \cdot, \cdot \ran$.
\end{definition}

\begin{remark}\label{d_Bf}
\emph{Let $A$ be a commutative associative algebra over $\R$ and $\mathfrak{M}$ a certain category of (bi)modules over $A$.
According to the theory of derivations (over $\R$) of $A$  into $\mathfrak{M}$ \cite{cartan}, there is a pair $(P, \delta)$, where $P$ is an object of $\mathfrak{M}$ and $\delta : A \to P$
is a derivation from $A$ to $P$, which has the \emph{universal property}: For any (bi)module $P'\in \mathfrak{M}$ and any derivation $\delta' : A \to P'$,
there exists a unique homomorphism $\alpha : P \to P'$ such that $\delta' = \alpha \circ \delta$. Considering the commutative associative algebra $A=\C$ over $\R$
and the space $\Omega^1 = \Gamma(T^{\ast} M)$, we have that the pair $(\Omega^1, d)$, where $d : \C \to \Omega^1$ is the usual derivation of smooth functions on $M$,
is the universal derivation in the category of geometric $\C$-modules \cite{nestr}. Moreover, the following two propositions are true \cite{nestr}:
\begin{enumerate}
\item
A $\C$-module $P$ is geometric, i.e. $\bigcap_{x\in M}\mu_x P = \{0\}$, if and only if $P$ is isomorphic to $\Gamma(P)$. The latter is the set of all sections of the (pseudo)bundle $|P| = \bigcup_{x\in M}P_x$ over $M$, where $P_x = P/\mu_xP$ and $\mu_x =\{f\in \C \, / \, f(x)=0\}$ is the ideal of functions in $\C$ vanishing at $x$.
\item
For any smooth vector bundle $B \to M$ the sequence
\[0 \to \mu_x \Gamma(B) \to \Gamma(B) \to B_x \to 0,\]
where the first arrow is the inclusion while the second assigns to every section its value at $x\in M$, is exact. Hence, $\Gamma(B)/\mu_x\Gamma(B) \cong B_x$.
\end{enumerate}
Applying this to predual vector bundles $B$ we have that
\begin{enumerate}
\item
$\Gamma(B)$ is a geometric $\C$-module, since
\[|\Gamma(B)| =  \bigcup_{x\in M}\Gamma(B)_x = \bigcup_{x\in M} \Gamma(B)/\mu_x\Gamma(B)\cong \bigcup_{x\in M} B_x = B.\]
\item
By the universal property of $(\Omega^1, d)$ in the category of geometric $\C$-modules \cite[Theorem 11.43]{nestr} and for the derivation $d_B: \C \to \Gamma(B)$ of Definition \ref{pre-dual}, there is a unique homomorphism $\alpha \in \mathrm{Hom}_{\C}(\Omega^1, \Gamma(B))$ such that
\begin{linenomath*} \begin{equation*}\label{definition d_Bf}
d_B = \alpha \circ d.
\end{equation*} \end{linenomath*}
\end{enumerate}
Moreover, since $\alpha : \Omega^1 \to \Gamma(B)$ is $\C$-linear, we have that there is a smooth bundle map $\alpha : T^\ast M \to B$ over $M$
such that $\alpha (\eta) = \alpha \circ \eta$, for all $\eta \in \Omega^1$ \cite[Lemma 10.29]{lee}.}

\vspace{1mm}
\noindent
\emph{When $B = E$, it is $d_B=d_E$ and $\alpha = {g^{\flat}}^{-1}\circ \rho^\ast $ (see Definition \ref{def Courant algebroid}).}
\end{remark}

\begin{example}\label{predual-T*M-Ex}
\emph{Let $(E, \lcf \cdot, \cdot\rcf, \langle \cdot, \cdot \rangle, \rho)$ be a Courant algebroid over $M$. Consider the pairing  $\lan \cdot, \cdot \ran : \Gamma(E) \times \Gamma(T^\ast M) \to \C$ given, for any $e \in \Gamma(E)$ and $\eta \in \Gamma(T^\ast M)$, by
\begin{equation}\label{predual-T*M}
\lan e, \eta \ran : = \langle \eta, \rho(e)\rangle,
\end{equation}
where the pairing at the right-hand side is the usual one between $1$-forms and vector fields.} The triple $(T^\ast M, d, \lan \cdot, \cdot \ran)$ is a predual of $E$ and $\alpha=\mathrm{Id}_{T^\ast M}$.
\end{example}

\begin{remark}\label{rank a}
\emph{Note that the rank of $\alpha : \Gamma(T^\ast M) \to \Gamma(B)$ might not be maximal or even constant on $M$. However, under mild assumptions and using results  in \cite{AZ} and \cite{Andr-Scandalis}\footnote{Although these results concern projective singular foliations they also apply to projective finitely generated $\C$-submodules of the $\C$-module of smooth sections of any real smooth vector bundle.}, one can show that it is constant almost everywhere as we now explain.}

\vspace{1mm}
\noindent
\emph{Consider the $\C$-submodule $\mathcal{B}=\alpha(\Omega^1_c)$ of $\Gamma(B)$, where $\Omega^1_c$ denotes the compactly supported sections of $T^\ast M$,
and suppose that $\mathcal{B}$ is locally finitely generated and projective\footnote{This setting contains the regular case and quite many singular situations, namely the ones now described as \emph{almost regular}. The general case is more subtle but one can reproduce the proof along the same lines presented here.}.}

\vspace{1mm}
\noindent
\emph{Let $x\in M$.}
\begin{itemize}
\item[-]
\emph{The quotient $\tilde{\mathcal{B}}_x : = \mathcal{B}/\mu_x\mathcal{B}$ is a finite dimensional vector space and $\dim \mathcal{B}_x \leq \dim \tilde{\mathcal{B}}_x = m \leq \dim T_x^\ast M =n$, where $\mathcal{B}_x=\alpha_x(T^\ast_x M)$.}
\item[-]
\emph{The evaluation map $ev_x: \tilde{\mathcal{B}}_x \to \mathcal{B}_x$ given, for any $[\sigma] \in \tilde{\mathcal{B}}_x$, by
\[ev_x([\sigma])=\sigma(x), \quad \quad \sigma \in \mathrm{Im}\alpha,\]
is a surjective homomorphism.}
\item[-]
\emph{The map $\alpha : \Omega^1_c \to \mathcal{B}$ induces a surjective homomorphism $\tilde{\alpha}_x : T^\ast_x M \to \tilde{\mathcal{B}}_x$ and one has the following commutative diagram:}
\begin{linenomath*} \begin{equation*}
\begin{CD}
T^\ast_x M @>{\alpha_x} >> \mathcal{B}_x \\
@V\tilde{\alpha}_x VV          @AA ev_x A \\
\tilde{\mathcal{B}}_x @ = \tilde{\mathcal{B}}_x,
\end{CD}
\end{equation*} \end{linenomath*}
\emph{i.e.}
\begin{linenomath*} \begin{equation}\label{eval}
ev_x \circ \tilde{\alpha}_x = \alpha_x.
\end{equation} \end{linenomath*}
\emph{The set $M_0$ of points of $M$ where $ev$ is bijective, i.e. the set of continuity of $x \mapsto \dim \mathcal{B}_x$, is open and dense in $M$. Moreover, for any $x\in M_0$, $\tilde{\mathcal{B}}_x = \mathcal{B}_x$.}
\item[-]
\emph{The function $x\mapsto\dim \tilde{\mathcal{B}}_x $ is constant since $\mathcal{B}$ is supposed to be projective \cite[Lemma 1.6]{AZ}.}
\item[-]
\emph{A basis of $\tilde{\mathcal{B}}_x$ is lifted to a set of generators of $\mathcal{B}\vert_U$, where $U$ is a neighborhood of $x$. For any $y\in U$, the values of these generators at $y$ give us a family of generators of $\mathcal{B}_y$ whose projections on $\tilde{\mathcal{B}}_y$ produces a basis of $\tilde{\mathcal{B}}_y$ since $\dim \tilde{\mathcal{B}}_y = \dim \tilde{\mathcal{B}}_x$.}
\item[-]
\emph{$\tilde{\mathcal{B}}\vert_U = \bigcup_{x \in U}\tilde{\mathcal{B}}_x$ is a trivial vector bundle over $U$ of rank $m$. By considering an open cover of $M$ and by gluing the trivial vector bundles constructed as above (this is possible, see \cite{AZ, Andr-Scandalis}) on the open sets of the cover, we take a vector bundle $\tilde{\mathcal{B}}= \bigcup_{x \in M}\tilde{\mathcal{B}}_x$ over $M$ of constant rank $m$. The map $ev : \tilde{\mathcal{B}} \to B$ is a vector bundle map inducing an isomorphism $\Gamma_c(\tilde{\mathcal{B}}) \cong \mathcal{B}$ at the level of smooth sections.}
\end{itemize}
\emph{Thus, $\mathrm{rank}\tilde{\alpha}_x = \dim \tilde{\mathcal{B}}_x = m$ is a constant function of $x$ on $M_0$ and $\mathrm{rank}\tilde{\alpha}_x = \mathrm{rank}\alpha_x$ on $M_0$. Consequently, $\alpha$ is of constant rank $m$ on the open and dense subset $M_0$ of $M$.}
\end{remark}

\begin{lemma}\label{anchored B}
To every predual $(B, d_B, \lan \cdot, \cdot \ran)$ of $E$ corresponds a $\C$-linear map $\varrho : \Gamma(B) \to \Gamma(TM)$.
\end{lemma}
\begin{proof}
Let $(B, d_B, \lan \cdot, \cdot \ran)$ be a predual of $E$. For any $b\in \Gamma(B)$, the $\R$-linear endomorphism $\lan d_E \cdot, b \ran :\C \to \C$ is, evidently, a derivation of $\C$, therefore corresponds to a vector field on the base manifold $M$, noted as $\varrho(b)$. Hence, we have a $\C$-linear map $\varrho : \Gamma(B) \to \Gamma(TM)$ such that, for any $f\in \C$,
\begin{equation}\label{anchor-B}
\varrho(b)(f) = \lan d_E f, b \ran.
\end{equation}
\end{proof}

\noindent
The induced vector bundle map $\varrho : B \to TM$ will be called \emph{the anchor map of $B$}.

\begin{remark}\label{remark vbm}
\emph{The vector bundle maps involved in the definition of a predual bundle of a Courant algebroid are related as follows. Let $g^\flat : E \to E^\ast$ be the isomorphism defined by $\langle \cdot, \cdot \rangle$ and $p : B \to E^\ast$ the morphism defined by $\lan \cdot, \cdot \ran$. We have that: }
\begin{itemize}
\item[-]
\emph{For all $e^1$, $e^2$ in the same fiber of $E^\ast$},
\begin{equation*}
\langle {g^\flat}^{-1}(e^1), ({g^\flat}^{-1})^\ast(e^2) \rangle = e^1 (({g^\flat}^{-1})^\ast(e^2)) = e^1 ({g^\flat}^{-1}(e^2)) = \langle {g^\flat}^{-1}(e^1), {g^\flat}^{-1}(e^2) \rangle,
\end{equation*}
\emph{hence we get $({g^\flat}^{-1})^\ast ={g^\flat}^{-1}$.}
\item[-]
\emph{For any $f\in \C$ and $e\in E$,}
\begin{equation*}
p(d_Bf) = (p\circ \alpha)(df) \quad \mathrm{and} \quad p(d_Bf)(e) = \lan e, d_Bf\ran = \rho(e)(f) = \rho^\ast(df)(e),
\end{equation*}
\emph{thus $p\circ \alpha = \rho^\ast$.}
\end{itemize}
\emph{Also, for any $f\in \C$ and $b\in B$, we have}
\begin{equation*}
\varrho(b)(f) = \lan d_Ef, b \ran = p(b)(({g^\flat}^{-1}\circ \rho^\ast)(df)) = (\rho \circ ({g^\flat}^{-1})^\ast \circ p )(b)(f) = (\rho \circ {g^\flat}^{-1}\circ p )(b)(f),
\end{equation*}
\emph{so $\varrho = \rho \circ {g^\flat}^{-1}\circ p$. In particular,}
\begin{equation*}
\varrho \circ \alpha = (\rho \circ {g^\flat}^{-1})\circ (p \circ \alpha) = \rho \circ {g^\flat}^{-1}\circ \rho^\ast.
\end{equation*}
\emph{Hence, $(\varrho \circ \alpha)(df) = \rho(d_Ef)=0$, for any $f\in \C$, and we get $\varrho \circ \alpha = 0$. The  relations above make the following diagram commutative:}
\begin{equation*}
\begin{tikzcd}
T^\ast M \arrow[r, "\alpha"]\arrow[rr, bend left, "p\circ \alpha"]\arrow["0", dr] & B \arrow[r,"p"]\arrow[d, "\varrho"] & E^\ast \arrow[d, "{g^\flat}^{-1}"]\\
& TM & E\arrow[l,"\rho"].
\end{tikzcd}
\end{equation*}
\end{remark}

\subsection{Dorfman connections}
In this subsection we present and study a modified definition of the notion of \emph{Dorfman connection} from the one introduced in \cite{mjl}. The difference is that we consider a Courant algebroid $(E, \lcf \cdot, \cdot\rcf, \langle \cdot, \cdot \rangle, \rho)$ in the place of the dull algebroid $(Q, [\cdot, \cdot]_Q, \rho_Q)$  \cite[Definition 3.3]{mjl} and that the predual structure is given by Definition \ref{pre-dual}. 

\begin{definition}\label{dorfman connection}
Let $(E, \lcf \cdot, \cdot\rcf, \langle \cdot, \cdot \rangle, \rho)$ be a Courant algebroid over a smooth manifold $M$ and $(B, d_B, \lan \cdot, \cdot \ran)$ a predual of $E$. An \emph{$E$-Dorfman connection on $B$} is an $\R$-bilinear map
\begin{linenomath*} \begin{equation*}
\nabla : \Gamma(E) \times \Gamma(B) \to \Gamma(B)
\end{equation*} \end{linenomath*}
such that, for all $e \in \Gamma(E)$, $b\in \Gamma(B)$, and $f\in \C$, the following three properties hold:
\begin{enumerate}
\item
$\nabla\!_{fe}b = f\nabla_eb + \lan e, b\ran d_Bf$,
\item
$\nabla\!_{e}(fb) = f\nabla_eb + \rho(e)(f)b$,
\item
$\nabla\!_{e}(d_Bf) = d_B (\mathcal{L}_{\rho(e)}f)$.
\end{enumerate}
\end{definition}

\begin{example}\label{dorfman on T*M}
\emph{Consider a Courant algebroid $(E, \lcf \cdot, \cdot\rcf, \langle \cdot, \cdot \rangle, \rho)$ over $M$ and its predual
$(T^\ast M, d, \lan \cdot, \cdot \ran)$ as in Example \ref{predual-T*M-Ex}. One can easily check that the map
\begin{eqnarray}\label{dorfman on T*M - eq}
\triangledown : \Gamma(E) \times \Gamma(T^\ast M) & \to & \Gamma(T^\ast M) \nonumber \\
(e, \eta) & \mapsto & \triangledown_e \eta : = \mathcal{L}_{\rho(e)}\eta
\end{eqnarray}
defines an $E$-Dorfman connection on $T^\ast M$.}
\end{example}

\begin{remarks}\label{rem on dorf con}

\noindent
\begin{enumerate}
\item
\emph{There are two equivalent ways to read the conditions of Definition \ref{dorfman connection}.}
\begin{enumerate}
\item[(i)]
\emph{The spaces of sections $\Gamma(E\otimes B)\cong \Gamma(E)\otimes \Gamma(B)$ and $\Gamma(B)$ being $\C$-bimodules,
the first two conditions of Definition \ref{dorfman connection} say that a Dorfman connection $\nabla$ is a first-order differential operator on the
bimodules $\Gamma(E\otimes B)$ and $\Gamma(B)$ in the sense of \cite{dub-viol-mass}. The third condition says that the subspace $\mathrm{Im}d_B \subset \Gamma(B)$ is invariant under the map $\nabla_e$, for any $e\in \Gamma(E)$.}
\item[(ii)] \label{napla as dif op}
\emph{Considering the Atiyah algebroid $\mathbb{A}(B)$ over $M$\footnote{Recall that the smooth sections of the Atiyah algebroid $\mathbb{A}(B)$ of a vector bundle $B\to M$,  are the derivative endomorphisms of $\Gamma(B)$, i.e. the maps $D: \Gamma(B) \to \Gamma(B)$ for which there exists a vector field $\sigma_D$ on $M$ such that, for any $f\in \C,\;b\in \Gamma(B)$, it is $D(fb)=fD(b) + \sigma_D(f)b$. The bracket on $\mathbb{A}(B)$ is the commutator bracket of endomorphisms and the anchor map is $D \mapsto \sigma_D$ \cite{yks-mck, mck}.}, the conditions of Definition \ref{dorfman connection} imply that an $E$-Dorfman connection on $B$ can be also viewed as an additive operator $\nabla : \Gamma(E)\to \Gamma(\mathbb{A}(B))$, $e\mapsto \nabla(e):=\nabla_e$, which is a differential operator of order 1 in the sense of \cite[Definition 2.2]{grab-k-pon} and whose values $\nabla_e$ leave the space $\mathrm{Im}d_B$ invariant. In fact, for any $f\in \C$, let $m^{E}_f$ and $m^{\mathbb{A}(B)}_f$ be the linear operators, multiplication by $f$, provided by the $\C$-module structure of $\Gamma(E)$ and $\Gamma(\mathbb{A}(B))$, respectively. The symbol $\sigma_{\nabla}(f) : \Gamma(E) \to \Gamma(\mathbb{A}(B))$ of $\nabla$,}
\[ \sigma_{\nabla}(f) : = \nabla \circ m^E_f - m^{\mathbb{A}(B)}_f \circ \nabla,\]
\emph{is a 0- order differential operator of $\Gamma(E)$ into $\Gamma(\mathbb{A}(B))$. We have}
\[\sigma_{\nabla}(f)(e) = \nabla(fe)-f\nabla(e) = \nabla_{fe}-f\nabla_e = \lan e, \cdot \ran d_Bf,\]
\emph{and $\sigma_{\nabla}(f) \circ m^E_g - m^{\mathbb{A}(B)}_g \circ \sigma_{\nabla}(f)$ vanishes for all $g\in \C$.}
\end{enumerate}
\item
\emph{Another way to read the third condition of Definition \ref{dorfman connection}, is the following. Consider the $E$-Dorfman connection $\triangledown$ on $T^\ast M$ defined in Example \ref{dorfman on T*M}.
Writting $d_B = \alpha \circ d$ (see Remark \ref{d_Bf}), we have}
\begin{equation*}
\nabla\!_{e}(d_Bf) = d_B (\mathcal{L}_{\rho(e)}f) \Leftrightarrow \nabla\!_{e}\alpha (df) = \alpha (\mathcal{L}_{\rho(e)}df) \Leftrightarrow \nabla\!_{e}\alpha (df) = \alpha (\triangledown_e df) ,
\end{equation*}
\emph{meaning that the following diagram is commutative:}
\begin{equation*}
\begin{CD}
\Gamma(T^\ast M) \supset \mathrm{Im}d @>{\alpha} >> \Gamma(B) \\
@V{\triangledown_e = \mathcal{L}_{\rho(e)}}VV          @VV{\nabla_e}V \\
\Gamma(T^\ast M)\supset \mathrm{Im}d @ >>{\alpha}> \Gamma(B).
\end{CD}
\end{equation*}
\item
\emph{The failure of an $E$-Dorfman connection to be a linear $E$-connection is controlled by the derivation $d_B=\alpha \circ d$. If $\alpha=0$, $\nabla$ is a linear $E$-connection on $B$ in the sense of \cite{Al-Xu, CM}.}
\end{enumerate}
\end{remarks}

\begin{proposition}\label{existence Dorfman connection}
Let $(E, B)$ be as in Definition \ref{dorfman connection}. The set $\mathfrak{C}(E,B)$ of $E$-Dorfman connections on $B$ is not empty.
\end{proposition}
\begin{proof} The proof of this proposition is organized in several steps.

\vspace{1mm}
\noindent
\textit{1st Step: Construction of $E$-Dorfman connections locally - Part A.}

\vspace{1mm}
\noindent
Let $(U, x^1, \ldots, x^n)$ be a local coordinate system of $M$. Over $U$, the vector bundles $E$, $E^\ast$, $B$ and $B^\ast$ can be trivilized, and we can choose
$(e_1, \ldots, e_r)$ a local frame of smooth sections of $E\vert_U$ with $(e^1, \ldots, e^r)$ its dual frame of local smooth sections of $E^\ast\vert_U$, $r = \mathrm{rank}E$, and $(b_1,\ldots,b_s)$ a local frame of smooth sections of $B\vert_U$ with $(b^1, \ldots, b^s)$ its dual frame of local sections of $B^\ast\vert_U$, $s = \mathrm{rank}B$.
With respect to these choices, the local expressions of the homomorphisms $\alpha : \Gamma(T^\ast M\vert_U)\to \Gamma(B\vert_U)$, $\lan \cdot, \cdot \ran : \Gamma(E\vert_U) \times \Gamma(B\vert_U) \to C^{\infty}(U, \R)$ and $\rho : \Gamma(E\vert_U) \to \Gamma(TM\vert_U)$ are,
respectively\footnote{We adopt this writing for the maps $\alpha$, $\lan \cdot, \cdot \ran$ and $\rho$
so that the matrices of their coefficients correspond to those used for the calculation of their images.},
\begin{equation*}
\alpha = \alpha^{ij}b_i\otimes \frac{\partial}{\partial x^j}, \quad \quad \lan \cdot, \cdot \ran = p_{ij}b^i\otimes e^j,
\quad \quad \rho = \rho^i_j \frac{\partial}{\partial x^i}\otimes e^j.
\end{equation*}
Their coefficients are smooth functions on $U$. With the same notation for their associated matrices, $A = (\alpha^{ij})$, $P= (p_{ij})$, and $\rho = (\rho^i_j)$,
we have that \eqref{def-dB} is equivalent to
\begin{linenomath*} \begin{equation}\label{A-P-R}
A^TP = \rho \Leftrightarrow P^TA = \rho^T
\end{equation} \end{linenomath*}
(see, also, Remark \ref{remark vbm}). Let $\pi: E \to M$ and $\tau : B \to M$ be the projections of $E$ and $B$ on $M$, respectively. Consider the $\R$-bilinear map $\nabla^0 : \Gamma(E\vert_U) \times \Gamma(B\vert_U) \to \Gamma(B\vert_U)$ defined, for any $e\in \pi^{-1}(U)\subset \Gamma(E\vert_U)$ and $b\in \tau^{-1}(U)\subset \Gamma(B\vert_U)$, $b = f^1b_1 + \ldots + f^sb_s$ with $f^i\in C^\infty(U,\R)$, by
\begin{linenomath*} \begin{equation} \label{dorfman trivial}
\nabla^0\!_e b = \nabla^0\!_e (f^1b_1 +  \ldots + f^sb_s) : = (\mathcal{L}_{\rho(e)}f^1)b_1 + \ldots + (\mathcal{L}_{\rho(e)}f^s)b_s + \sum_{i=1}^s f^i d_B(P_i(e)),
\end{equation} \end{linenomath*}
where $P_i(e)$ is the $i$-component function of the local section $\lan e, \cdot \ran$ of $B^{\ast}\vert_U$. One can easily check that \eqref{dorfman trivial} satisfies the first and second
condition of Definition \ref{dorfman connection}, but not the third. The latter is satisfied under the following strong condition on $\alpha$: For any $k = 1, \ldots, s$,
\begin{linenomath*} \begin{equation*}\label{diff equation for a}
\rho^l_t \frac{\partial \alpha^{kj}}{\partial x^l} + \alpha^{ij}\alpha^{km}\frac{\partial p_{it}}{\partial x^m } = \alpha^{kj}\frac{\partial \rho^l_t }{\partial x^j}\,
\stackrel{\eqref{A-P-R}}{\Leftrightarrow}\,
\alpha^{it}\frac{\partial \alpha^{kl}}{\partial x^t} = \alpha^{kt}\frac{\partial \alpha^{il}}{\partial x^t}.
\end{equation*} \end{linenomath*}
For this reason, we search for an homomorphism $C : \Gamma(E\vert_U\otimes B\vert_U)\cong \Gamma(E\vert_U)\otimes \Gamma(B\vert_U) \to \Gamma(B\vert_U)$  such that the map $\nabla : \Gamma(E\vert_U)\times \Gamma(B\vert_U) \to \Gamma(B\vert_U)$, with
\begin{linenomath*} \begin{equation} \label{def-napla}
\nabla = \nabla^0 +C,
\end{equation} \end{linenomath*}
is a Dorfman connection.

\vspace{1mm}
\noindent
\textit{2nd Step: Existence of $C$.}

\vspace{1mm}
\noindent
We now show that it is always possible to find locally such a homomorphism $C$. In the local coordinates of $M$ and the local frames of $E$ and $B$ considered above, $C$ is written as $C = C_{ij}^k b_k\otimes b^i \otimes e^j$, with $C_{ij}^k \in C^\infty(U, \R)$.
Obviously, \eqref{def-napla} verifies the first and second condition of Definition \ref{dorfman connection}, while the third is satisfied if and only if, for any $t=1,\ldots,r$ and $f\in C^\infty(U,\R)$,
\begin{linenomath*}\begin{eqnarray}\label{diff equation for a, c, f}
\nabla_{e_t}d_Bf = d_B(\mathcal{L}_{\rho(e_t)}f) & \Leftrightarrow & \big[\alpha^{ij}C_{it}^k + \rho^l_t \frac{\partial \alpha^{kj}}{\partial x^l} + \alpha^{ij}\alpha^{km}\frac{\partial p_{it}}{\partial x^m } - \alpha^{km}\frac{\partial \rho^j_t }{\partial x^m}\big]\frac{\partial f}{\partial x^j}b_k =0 \nonumber \\
& \Leftrightarrow & \alpha^{ij}C_{it}^k + \rho^l_t \frac{\partial \alpha^{kj}}{\partial x^l} + \alpha^{ij}\alpha^{km}\frac{\partial p_{it}}{\partial x^m } - \alpha^{km}\frac{\partial \rho^j_t }{\partial x^m} = 0 \nonumber \\
& \stackrel{\eqref{A-P-R}}{\Leftrightarrow} & \alpha^{ij}C_{it}^k + p_{it} \big(\alpha^{il}\frac{\partial \alpha^{kj}}{\partial x^l}- \alpha^{km}\frac{\partial \alpha^{ij}}{\partial x^m}\big)=0,
\end{eqnarray}\end{linenomath*}
for any $k=1,\ldots,s$ and $j=1,\ldots,n$. Setting $C_t = (C_{ti}^k)$ to be the matrix corresponding to the homomorphism $C(e_t, \cdot): \Gamma(B\vert_U) \to \Gamma(B\vert_U)$, and $N_t = (N_t^{kj})$, with
\begin{linenomath*} \begin{equation*}
N_t^{kj} = - p_{it} \big(\alpha^{il}\frac{\partial \alpha^{kj}}{\partial x^l}- \alpha^{km}\frac{\partial \alpha^{ij}}{\partial x^m}\big),
\end{equation*} \end{linenomath*}
equation \eqref{diff equation for a, c, f} is written, in matrix form, as
\begin{linenomath*} \begin{equation} \label{matrix A C_k-n}
 C_tA = N_t.
\end{equation} \end{linenomath*}
The left hand side in the last equation corresponds to the homomorphism $C(e_t, \cdot) \circ \alpha : \Gamma(T^\ast M\vert_U) \to \Gamma(B\vert_U)$. Hence, $\nabla = \nabla^0 + C$ verifies the third condition of Definition \ref{dorfman connection} if and only if, for any $t=1,\ldots, r$, there exists a matrix $C_t$ on $U$ satisfying \eqref{matrix A C_k-n}.

\vspace{1mm}
\noindent
In order to solve the last equation, we remark that $\alpha$ is locally of constant rank. In fact, let $A_x$ be the matrix of the map $\alpha_x : T^\ast_x M \to B_x$, $x\in U$. Since $\mathrm{rank}A_x = \dim \mathrm{Im}\alpha_x$, it is $\mathrm{rank}A_x \leq n$, so let $m = \max \{\mathrm{rank}A_x \, / \, x\in U\}$. Let $x_0 \in U$ be such that $\mathrm{rank}A_{x_0}=m$ and $A'_{x_0}$ be an $m\times m$ submatrix of $A_{x_0}$  such that $\det A'_{x_0} \neq 0$. Due to continuity of the function $\det$, it is $\det A'_y \neq 0$ for any $y$ in a neighborhood $U_0\subset U$ of $x_0$. Thus, $A$ is of constant rank $m$ on $U_0$. Hence, by restricting $U$ if necessary, we can suppose that the matrix $A$, and so the vector bundle map $\alpha$, is of constant rank $m$ on $U$. By reordering the coordinates and the frame of smooth sections of $B$, we may assume that the first $m$ rows of $A$ are linearly independent on $U$. In this case, since $m\leq \min \{s, n\}$, the matrix $A$ can be
written in block form as $A = \left(\begin{array}{c}
 A^1_{m\times n} \\
   A^2_{(s-m)\times n} \\
\end{array} \right)$ with $\mathrm{rank}A^1_{m\times n} = m$ on $U$. Similarly, the matrix $C_t$ is written in block form as $C_t = \left(
                                                                                       \begin{array}{cc}
                                                                                         C^1_{t, s\times m} & C^2_{t, s \times (s-m)} \\
                                                                                       \end{array}
                                                                                     \right)
$. Then, equation \eqref{matrix A C_k-n} is equivalent to
\begin{equation}\label{new-C-A-new}
C^1_t A^1 + C_t^2A^2 = N_t \Leftrightarrow C^1_t A^1 = N_t - C_t^2A^2.
\end{equation}
Since $\mathrm{rank}A^1_{m\times n} = m$ on $U$, it has a right inverse matrix $A^{1^R} = A^{1^T} (A^1 A^{1^T})^{-1}$. So, \eqref{new-C-A-new} gives us
\begin{equation*}
C^1_t = (N_t - C_t^2A^2)A^{1^R}.
\end{equation*}
Consequently, for any choice of the block $C^2_t$ of $C_t$, the last equation determines the block $C_t^1$ and thus
the matrix $C_t$, for $t=1, \ldots, r$, and, equivalently, the homomorphism $C$.

\vspace{1mm}
\noindent
\textit{3rd Step: Construction of $E$-Dorfman connections globally.}

\vspace{1mm}
\noindent
Choose an open cover $(U_i)_{i \in I}$ of $M$ such that, for any $i\in I$, $E \vert_{U_i}$ and $B \vert_{U_i}$ are
trivial bundles, and pick a smooth partition of unity $(\psi_j)_{j\in J}$ subordinate to this cover, i.e. $\mathrm{supp}\psi_j \subset U_i$ for some $i=i(j)$. By the previous discussion, each bundle $B \vert_{U_i}$ admits at least one $E \vert_{U_i}$-Dorfman connection $\nabla^i$
of type \eqref{def-napla}. Denote such a connection by $\nabla^i$ and, for any $e\in \Gamma(E)$ and $b\in \Gamma(B)$, set
\begin{linenomath*} \begin{equation*}\label{conn-affin}
\nabla_eb = \sum_{j\in J}\psi_j(\nabla^i_{e\vert_{U_i}}b\vert_{U_i}).
\end{equation*} \end{linenomath*}
Since, the set  of supports $\{\mathrm{supp}\psi_j\;/ j\in J\}$  is locally finite, the sum $\sum_{j\in J}\psi_j(x) =1$, for all $x\in M$, has only finitely many nonzero terms in a neighborhood of each point. Thus, it is easy to check that $\nabla$ is an $E$-Dorfman connection on $B$. So, $\mathfrak{C}(E,B)$ is non empty.
\end{proof}

\begin{remark}
\emph{Under the assumptions in Remark \ref{rank a}, the second step of the last proof, the solution of \eqref{matrix A C_k-n}, can be made more precise. Suppose that the submodule $\mathcal{B} = \alpha(\Omega_c^1)$ of $\Gamma(B)$ is finitely generated and projective.}

\vspace{1mm}
\noindent
\emph{Let $x\in U$, $A_x$ be the matrix of the map $\alpha_x : T^\ast_x M \to B_x$, $\tilde{A}_x$ the matrix of $\tilde{\alpha}_x : T^\ast_xM \to \tilde{\mathcal{B}}_x$ and $\mathcal{EV}_x$ the matrix of $ev_x$ in appropriate bases of the corresponding spaces. Equation \eqref{eval} is then written in matrix form as $\mathcal{EV}_x \tilde{A}_x = A_x$. Thus, $C_{tx}\mathcal{EV}_x \tilde{A}_x = C_{tx}A_x$, and \eqref{matrix A C_k-n} is equivalent to
\begin{linenomath*} \begin{equation}\label{new-C-A}
C_{t x}\mathcal{EV}_x \tilde{A}_x = N_{t x}.
\end{equation} \end{linenomath*}
Since $\tilde{A}_x$ is an $m \times n$ matrix with $\mathrm{rank}\tilde{A}_x=m$, which is a constant function of $x$ on $U$, it has a right inverse $\tilde{A}_x^R = \tilde{A}_x^T(\tilde{A}_x \tilde{A}_x^T)^{-1}$. So, \eqref{new-C-A} is equivalent to
\begin{linenomath*} \begin{equation}\label{new-C-A-solution}
C_{t x}\mathcal{EV}_x = N_{t x} \tilde{A}_x^R.
\end{equation}\end{linenomath*}
Clearly, $\mathcal{EV}_x$, $N_{t x}$ and $\tilde{A}_x^R$ depend smoothly on $x\in U$. At each point $x\in U_0=M_0 \cap U$, equation \eqref{new-C-A-solution} has a unique solution $C_{t x}=(C_{ti}^k(x))$,
since $\mathcal{EV}_x$ is invertible on $U_0$, and depends smoothly on $x\in U_0$. Each function $C_{ti}^k : U_0 \to \R$ has a
smooth extension $\tilde{C}_{ti}^k$ on $U$ such that $\tilde{C}_{ti}^k\vert_{U_0}= C_{ti}^k$\footnote{The set $U_0$ as open subset of $M$ is an open submanifold of $M$ which is considered to be an embedded submanifold of codimension zero. Then, by \emph{Extension Lemma for functions on submanifolds} \cite[Lemma 5.34]{lee}, every smooth function $f$ on $U_0$ has a smooth extension $\tilde{f}$ on $U$ such that $\tilde{f}\vert_{U_0} = f$.}. Hence, the matrix $\tilde{C}_{tx}=(\tilde{C}_{ti}^k(x))$ is a solution of \eqref{new-C-A-solution} at $x\in U$. Indeed,
\begin{itemize}
\item[-]
if $x\in U_0$, the above claim is true.
\item[-]
if $x\in U\setminus U_0$, since $U_0$ is an open and dense subset of $U$, every point $x \in U\setminus U_0$ is a limit point of $U_0$, i.e.,
there exist a sequence $(x_n)_{n\in \N^\ast}$ of points of $U_0$ such that $\lim_{n\to \infty}x_n = x$. By continuity of the functions $\tilde{C}_{ti}^k$ on $U$,
we have that $\lim_{n\to \infty} \tilde{C}_{ti}^k(x_n) = \tilde{C}_{ti}^k(x) \Leftrightarrow \lim_{n\to \infty} C_{ti}^k(x_n) = \tilde{C}_{ti}^k(x)$.
Also, at each point $x_n$, $n\in \N^\ast$, \eqref{new-C-A-solution} is true, i.e. $C_{t x_n}\mathcal{EV}_{x_n}= N_{t x_n} \tilde{A}_{x_n}^R$. Taking limits in the last equation, we have the required result.
\end{itemize}
Therefore, it is always possible to find, locally (on $U$), a homomorphism $C : \Gamma(E\vert_U)\times \Gamma(B\vert_U) \to \Gamma(B\vert_U)$ such that $\nabla = \nabla^0 + C$ defines an $E$-Dorfman connection on $B$.}
\end{remark}

\begin{proposition}\label{affine space}
The set $\mathfrak{C}(E,B)$ carries a natural affine structure with corresponding linear space
\begin{equation}\label{cal-S}
\mathcal{S} = \big\{S \in \Gamma\big(\mathrm{Hom}_{\C}(E\otimes B, B)\big ) \; /\; \mathrm{Im}d_B \subseteq \bigcap_{e\in \Gamma(E)}\ker S(e, \cdot)\big\}.
\end{equation}
\end{proposition}
\begin{proof}
Let $\nabla^0$, $\nabla^1$ be two $E$-Dorfman connections on $B$. Then, for any $g \in \C$, the affine combination
\begin{linenomath*} \begin{equation*}
\nabla : = (1-g)\nabla^0 + g \nabla^1
\end{equation*} \end{linenomath*}
is an $E$-Dorfman connection on $B$. Thus, $\mathfrak{C}(E,B)$ is an affine space.

\vspace{1mm}
\noindent
Let $\nabla$ be an $E$-Dorfman connection on $B$ and $S$ an element of $\Gamma\big(\mathrm{Hom}_{\C}(E\otimes B, B)\big ) \cong \Gamma(E^\ast) \otimes \Gamma(B^\ast) \otimes \Gamma(B) \cong \mathcal{C}^1(\mathcal{E}; \mathcal{R}) \otimes \Gamma(\mathrm{End}(B))$. It is easy to check that $\nabla + S$ is an $E$-Dorfman connection on $B$ if and only if, for any $e\in \Gamma (E)$, $\mathrm{Im}d_B \subseteq \ker S(e, \cdot)$, i.e.
\[\mathrm{Im}d_B \subseteq \bigcap_{e\in \Gamma(E)}\ker S(e, \cdot).\]
Consider the space $\mathcal{S}$ defined in \eqref{cal-S}. Clearly, $\mathcal{S}$ is a non empty (because the zero homomorphism $0 : \Gamma(E)\otimes\Gamma(B) \to \Gamma(B)$, $0(e, \cdot) =0$, belongs to $\mathcal{S}$) linear space. On the other hand, if $\nabla$ and $\nabla'$ are two $E$-Dorfman connections on $B$, by the first condition of Definition \ref{dorfman connection} we get that, for any $f\in \C$, $e\in \Gamma(E)$ and $b\in \Gamma(B)$,
\begin{linenomath*} \begin{equation*}
\nabla\!_{fe}b - \nabla'\!_{fe}b = f(\nabla\!_{e} - \nabla'\!_{e})b.
\end{equation*} \end{linenomath*}
Thus $\nabla b - \nabla' b$ is a $\C$-linear homomorphism from $\Gamma(E)$ to $\Gamma(B)$. Using the Leibniz rule and the third condition of
Definition \ref{dorfman connection} for $\nabla$ and $\nabla'$ we obtain, respectively, that
\begin{linenomath*} \begin{equation*}
\nabla\!_{e}(fb) - \nabla'\!_{e}(fb) = f(\nabla\!_{e} - \nabla'\!_{e})b \quad \quad \mathrm{and} \quad \quad \nabla\!_{e}d_Bf - \nabla'\!_{e}d_Bf = 0.
\end{equation*} \end{linenomath*}
This shows that $\nabla\!_{e} - \nabla'\!_{e}$ is a $\C$-linear endomorphism of $\Gamma(B)$ vanishing identically on $\mathrm{Im}d_B$. As a result, $\nabla - \nabla'$ is an element of $\mathcal{S}$, which means that $\mathfrak{C}(E,B)$ is an affine space modeled on $\mathcal{S}$.
\end{proof}

\vspace{2mm}
\noindent
In the above framework and taking into account that $B$ is an anchored smooth vector bundle (Lemma \ref{anchored B}) with anchor map \eqref{anchor-B}, we have:

\begin{proposition}\label{prop - D - connection}
Every $E$-Dorfman connection $\nabla: \Gamma(E)\times\Gamma(B)\to \Gamma(B)$ on $B$ defines a linear $B$-connection $D:\Gamma(B)\times\Gamma(E)\to\Gamma(E)$ on $E$, in the sense of \cite{Cantrijn-Lan}\,\footnote{Let $A\to M$ be a vector bundle over a smooth manifold $M$ endowed with an anchor map $a : A \to TM$. According to \cite{Cantrijn-Lan}, a linear $A$-connection on a vector bundle $E\to M$ is a $\R$-bilinear operator $D : \Gamma(A)\times \Gamma(E) \to \Gamma(E)$ such that, for any $s\in \Gamma(A)$, $e\in \Gamma(E)$ and $f\in \C$,
\begin{equation*}
D_{fs}e = fD_se \quad \quad \mathrm{and} \quad \quad D_s(fe) = fD_se + a(s)(f)e.
\end{equation*}}, through the equation
\begin{equation}\label{b-connection}
\langle D_be,e'\rangle = \lan \lcf e,e'\rcf, b\ran + \lan e',\nabla_eb\ran-\rho(e)\lan e',b\ran,
\end{equation}
that has the property
\begin{equation}\label{property B-connection}
D_{d_Bf} = 0, \quad \quad \forall \, f\in \C.
\end{equation}
Conversely, every linear $B$-connection $D:\Gamma(B)\times\Gamma(E)\to\Gamma(E)$ on $E$ with the property \eqref{property B-connection} defines an $E$-Dorfman connection $\nabla : \Gamma(E) \times \Gamma(B) \to \Gamma(B)$ on $B$ through the relation
\begin{equation}\label{b-connection2}
\lan e',\nabla_eb\ran =  \rho(e)\lan e',b\ran  - \lan \lcf e,e'\rcf,b\ran + \langle D_be,e'\rangle,
\end{equation}
modulo some element $S \in \mathcal{S}\cap \mathcal{S}'$, where $\mathcal{S}$ is the set \eqref{cal-S} and
\[\mathcal{S}' = \big\{S \in \Gamma\big(\mathrm{Hom}_{\C}(E\otimes B, B)\big ) \; /\; \mathrm{Im}S \subseteq \bigcap_{e'\in \Gamma(E)}\ker \lan e', \cdot \ran\big\}.\]
\end{proposition}
\begin{proof}
Compute directly that
\begin{eqnarray*}
\langle D_{fb}e,e'\rangle & = &\lan \lcf e,e'\rcf,fb\ran+\lan e',\nabla_efb\ran-\rho(e)\lan e',fb\ran \nonumber \\
& = & f\lan \lcf e,e'\rcf,b\ran+\lan e',f\nabla_eb+\rho(e)(f)b\ran-\rho(e)\big(f\lan e',b\ran\big) \nonumber\\
& = & f\lan \lcf e,e'\rcf,b\ran+f\lan e',\nabla_eb\ran-f\rho(e)\lan e',b\ran \\
& = & f\langle D_be,e'\rangle\nonumber\\
& = & \langle fD_be,e'\rangle.
\end{eqnarray*}
Since $\langle\cdot,\cdot\rangle$ is nondegenerate, we then get $D_{fb}e=fD_be$. On the other hand,
\begin{eqnarray*}
\langle D_b(fe),e'\rangle & = &\lan \lcf fe,e'\rcf,b\ran+\lan e',\nabla_{fe}b\ran-\rho(fe)\lan e',b\ran \nonumber \\
& = & \lan f\lcf e,e'\rcf-\rho(e')(f)e+\langle e,e'\rangle d_Ef,b\ran \nonumber\\
& & +\, \lan e',f\nabla_eb+\lan e, b\ran d_Bf\ran -f\rho(e)\lan e',b\ran \nonumber\\
& = & f\lan \lcf e,e'\rcf,b\ran -\rho(e')(f)\lan e,b\ran +\langle e,e'\rangle\lan d_Ef,b\ran \nonumber \\
& & +\, f\lan e',\nabla_eb\ran +\lan e,b\ran\lan e',d_Bf\ran - f\rho(e)\lan e',b\ran\nonumber \\
& = & f\langle D_be,e'\rangle+\lan d_Ef,b\ran\langle e,e'\rangle\nonumber\\
& = & \langle  fD_be+\lan d_Ef,b\ran e,e'\rangle \nonumber \\
&\stackrel{\eqref{anchor-B}}{=} & \langle  fD_be + \varrho(b)(f)e,e'\rangle.
\end{eqnarray*}
By the nondegeneracy of $\langle\cdot,\cdot\rangle$, we get $D_b(fe)=fD_be + \varrho(b)(f)e$. In addition, for any $f\in \C$ and $e, e'\in \Gamma(E)$, we have
\begin{eqnarray*}
\langle D_{d_Bf}e,e'\rangle & = &\lan \lcf e,e'\rcf,d_Bf\ran+\lan e',\nabla_ed_Bf\ran-\rho(e)\lan e',d_Bf\ran \nonumber \\
& = & \lan \lcf e,e'\rcf,d_Bf\ran+\lan e', d_B(\mathcal{L}_{\rho(e)}f)\ran-\rho(e)\lan e',d_Bf\ran \nonumber \\
& \stackrel{\eqref{def-dB}}{=} & \rho(\lcf e,e'\rcf)(f)+\rho(e')\rho(e)(f)-\rho(e)\rho(e')(f)\nonumber\\
& = & \rho(\lcf e,e'\rcf)(f) - [\rho(e), \rho(e')](f)=0.
\end{eqnarray*}
Again, due to the nondegeneracy of $\langle\cdot,\cdot\rangle$, we conclude that $D$ verifies \eqref{property B-connection}.
As a consequence, \eqref{b-connection} defines a linear $B$-connection $D$ on $E$ that verifies \eqref{property B-connection}.
Hence, the first claim of the proposition is established.

\vspace{1mm}
\noindent
For the inverse direction we will show that a linear $B$-connection $D:\Gamma(B)\times\Gamma(E)\to\Gamma(E)$ on $E$ satisfying \eqref{property B-connection}, defines a Dorfman connection of $E$ on $B$ modulo an element $S\in \mathcal{S}\cap \mathcal{S}'$.
The connection $\nabla$ is defined through equation \eqref{b-connection2}. Indeed, it can be checked directly that this satisfies the properties in the Definition \ref{dorfman connection} of a Dorfman connection, as
\begin{eqnarray*}
\lan e',\nabla_{fe}b\ran & = & \rho(fe)\lan e',b\ran - \lan\lcf fe,e'\rcf,b\ran+\langle D_b(fe),e'\rangle \nonumber \\
& = & f\rho(e)\lan e',b\ran - \lan f\lcf e,e'\rcf,b\ran+ \rho(e')(f)\lan e,b\ran\nonumber \\
 &  & -\langle e,e'\rangle\lan d_Ef,b\ran + \langle fD_be,e'\rangle +\langle e,e'\rangle\lan d_Ef,b\ran \nonumber \\
  & = & \lan e',f\nabla_eb+\lan e,b\ran d_Bf\ran,
\end{eqnarray*}
and for the other two conditions it is
\begin{eqnarray*}
\lan e',\nabla_e(fb)\ran & = & \rho(e)\lan e',fb\ran - \lan\lcf e,e'\rcf,fb\ran+\langle D_{fb}e,e'\rangle \nonumber \\
& = & \rho(e)(f)\lan e',b\ran + f\rho(e)\lan e',b\ran - f\lan \lcf e, e'\rcf,b\ran + \langle fD_be,e'\rangle\nonumber \\
  & = & \lan e',f\nabla_eb+\rho(e)(f)b\ran,
\end{eqnarray*}
and
\begin{eqnarray*}
\lan e',\nabla_ed_Bf\ran & = & \rho(e)\lan e',d_Bf\ran - \lan\lcf e,e'\rcf,d_Bf\ran+\langle D_{d_Bf}e,e'\rangle \nonumber \\
& \stackrel{\eqref{def-dB}}{=} & \rho(e)\rho(e')(f) -\rho(\lcf e,e'\rcf)(f)\nonumber \\
& = & \rho(e)\rho(e')(f)-[\rho(e), \rho(e')](f)\nonumber \\
& = & \rho(e')\rho(e)(f)\nonumber \\
& \stackrel{\eqref{def-dB}}{=}  & \lan e',d_B(\mathcal{L}_{\rho(e)}f)\ran.\nonumber\\
\end{eqnarray*}
Thus $D$ defines an $E$-Dorfman connection $\nabla$ on $B$ modulo an element $S\in \mathcal{S}\cap \mathcal{S}'$.
\end{proof}

\begin{remark}
\vspace{2mm}
\emph{When $B=E$ and $\lan \cdot, \cdot \ran$ is the symmetric bilinear form  $\langle \cdot, \cdot \rangle$ of Definition \ref{def Courant algebroid}, the linear $E$-connection $D$ on $E$ defined by an $E$-Dorfman connection $\nabla$ on $E$ via \eqref{b-connection} is \begin{equation*}
D_{e_1}e_2 = \nabla_{e_2}e_1 - \lcf e_2, e_1\rcf.
\end{equation*}
Note also that $D_{d_Ef}e = \nabla_ed_Ef - \lcf e, d_Ef\rcf = d_E(\mathcal{L}_{\rho(e)}f) - d_E(\mathcal{L}_{\rho(e)}f)=0$, as required.}

\vspace{1mm}
\noindent
\emph{It is natural to ask when $\nabla$ is compatible with $\langle \cdot, \cdot \rangle$, i.e. when
\[\rho(e)\langle e_1, e_2\rangle = \langle \nabla_ee_1, e_2\rangle + \langle e_1, \nabla_e e_2\rangle.\]
A direct check says that this is the case if and only if}
\[\langle D_{e_1}e, e_2\rangle + \langle e_1, D_{e_2}e\rangle =0.\]

\vspace{1mm}
\noindent
\emph{On the other hand, $D$ is compatible with $\langle \cdot, \cdot\rangle$, i.e.
\[\rho(e)\langle e_1, e_2\rangle = \langle D_ee_1, e_2\rangle + \langle e_1, D_e e_2\rangle,\]
 if and only if
\[\langle \nabla_{e_1}e, e_2\rangle + \langle e_1, \nabla_{e_2}e\rangle =\rho(e_1)\langle e, e_2\rangle + \rho(e_2)\langle e_1,e\rangle.\]
The two compatibilities hold simultaneously when $E$ is a bundle of quadratic Lie algebras ($d_E=0)$.}
\end{remark}

\begin{remark}
\vspace{2mm}
\emph{As in the classical context, $D$ defines a linear $B$-connection $D^\ast$ on the dual vector bundle $E^\ast$ of $E$ via the relation
\begin{equation*}
\varrho(b) \langle e^\ast, e \rangle = \langle D^\ast_b e^\ast, e\rangle + \langle e^\ast, D_b e\rangle, \quad \quad e^\ast \in \Gamma(E^\ast), \; e\in \Gamma(E).
\end{equation*}
By a straightforward calculation, taking into account that $E^\ast \cong E$, $\lcf e^\ast,e\rcf + \lcf e,e^\ast\rcf = d_E \langle e^\ast, e\rangle$ and \eqref{b-connection}, we have that $D^\ast$ is given through the formula
\begin{equation}\label{B-ast-connection}
\langle D^\ast_b e^\ast, e\rangle = \lan \lcf e^\ast, e\rcf, b\ran - \lan e^\ast, \nabla_eb \ran + \rho(e)\lan e^\ast,b\ran.
\end{equation}
Furthermore, it is easy to see that, for any $f\in \C$,
\begin{equation}\label{D*-d_Bf}
D^\ast_{d_Bf}=0.
\end{equation}}
\end{remark}

\subsection{Curvature of Dorfman connections}\label{section-curvature}
For each $p\in \N$ and $m \in \N^{\ast}$, let $\mathfrak{D}_{m, m-1}^p(\mathcal{E};\Gamma(B))$ be the space of sequences of differential
operators $H= (H_0,H_1,\ldots,H_{[\frac{p}{2}]})$ determined by a $([\frac{p}{2}]+1)$-tuple of homomorphisms
\begin{linenomath*} \begin{equation*}
H_k : \mathcal{E}^{\otimes^{p-2k}}\otimes S^k\Omega^1 \to \Gamma(B)
\end{equation*} \end{linenomath*}
such that
\begin{enumerate}
\item
$H_k : \mathcal{E}^{\otimes^{p-2k}} \to \mathrm{Hom}_{\mathbb{K}}(S^k \Omega^1, \Gamma(B))$ takes values in the space $\mathfrak{Diff}_m(S_{\mathcal{R}}^k \Omega^1, \Gamma(B))\subset \mathrm{Hom}_{\mathbb{K}}(S^k \Omega^1, \Gamma(B))$, where $S_{\mathcal{R}}^k \Omega^1$ is the $\mathcal{R}$-module of the $k$-symmetric power of the $\mathcal{R}$-module $\Omega^1$ and
\begin{eqnarray*}
\lefteqn{\mathfrak{Diff}_m(S_{\mathcal{R}}^k \Omega^1, \Gamma(B)) = \big{\{}D \in \mathrm{Hom}_{\mathbb{K}}(S^k \Omega^1, \Gamma(B)) \, / \,D \,\, \mathrm{is \, \, a \,\, differential}} \nonumber \\
& & \,\mathrm{operator \,\, on} \,\,S^k \Omega^1\,\, \mathrm{of \,\, order \,\, at\,\, most} \,\,m\,\, \mathrm{in \,\, each \,\, entry\,\, with \,\, values \,\, in \,\, \Gamma(B)}\big{\}}.
\end{eqnarray*}
\item
$H_k : \mathcal{E}^{\otimes^{p-2k}} \to \mathfrak{Diff}_m(S_{\mathcal{R}}^k \Omega^1, \Gamma(B))$ is a differential operator of order at most $m$ on the first $p-2k-1$ arguments and of order at most $m-1$ on the $(p-2k)$-th argument of $\mathcal{E}^{\otimes^{p-2k}}$.
\end{enumerate}
Note that in the case $p=1$, the elements of $\mathfrak{D}^1_m(\mathcal{E}; \Gamma(B))$ are of type $H=(H_0)$ and the map $H_0 :\mathcal{E} \to \Gamma(B)$
is viewed as a differential operator of order at most $m$ on $\mathcal{E}$ with values in $\Gamma(B)$. Set $\mathfrak{D}(\mathcal{E}; \Gamma(B))=\big(\mathfrak{D}_{m, m-1}^p(\mathcal{E}; \Gamma(B))\big)_{p \in \N, \, m\in \N^\ast}$.

\begin{definition}\label{covariant Dorfman derivation}
An $E$-Dorfman connection $\nabla$ on $B$ defines a map
\begin{linenomath*} \begin{eqnarray*}
d^{\nabla} : \Gamma(B) & \to & \mathfrak{D}_1^1(\mathcal{E};\Gamma(B))\\
b &\mapsto & d^{\nabla}b : = \nabla_{\cdot}b,
\end{eqnarray*} \end{linenomath*}
called \emph{Dorfman covariant derivation}, satisfying, for any $f\in \C, b\in \Gamma(B)$, the Leibniz rule
\begin{linenomath*} \begin{equation*}
d^{\nabla}(fb) = d_Ef \otimes b + fd^{\nabla}b.
\end{equation*} \end{linenomath*}
\end{definition}

\vspace{2mm}
\noindent
We extend the multiplication $\cdot$ in $\mathfrak{D}(\mathcal{E}; \mathcal{R})$ (as subalgebra of $\mathcal{A}$, see formula (\ref{multiplication in A})) to a multiplication, also denoted by $\cdot$,
between elements of $\mathfrak{D}(\mathcal{E}; \mathcal{R})$ and of $\mathfrak{D}(\mathcal{E}; \Gamma(B))$, setting, for any $\omega \in \mathfrak{D}_{m, m-1}^p(\mathcal{E}; \mathcal{R})$
and $\eta \otimes b = (\eta_0,\eta_1, \ldots, \eta_{[\frac{q}{2}]})\otimes b \in \mathfrak{D}_{n, n-1}^q(\mathcal{E}; \Gamma(B))$,
\begin{linenomath*} \begin{equation*}
\omega \cdot (\eta \otimes b) : = (\omega \cdot \eta)\otimes b.
\end{equation*} \end{linenomath*}

\vspace{2mm}
\noindent
As in the theory of linear connections, one can prove  that the Dorfman covariant derivation $d^{\nabla}$ extends uniquely to an operator of degree $+1$, denoted also by $d^{\nabla}$,
\begin{linenomath*} \begin{equation} \label{ext - d - napla}
d^{\nabla} : \mathfrak{D}_{m,m-1}^p(\mathcal{E}; \Gamma(B)) \to \mathfrak{D}_{m+1, m}^{p+1}(\mathcal{E};\Gamma(B)),
\end{equation} \end{linenomath*}
satisfying the Leibniz rule
\begin{linenomath*} \begin{equation} \label{Leibniz - d - napla}
d^{\nabla}(\omega \otimes b) = d\omega \otimes b + (-1)^p\omega \cdot d^{\nabla}b,
\end{equation} \end{linenomath*}
for all $\omega \otimes b\in \mathfrak{D}_{m,m-1}^p(\mathcal{E};\Gamma(B))$. Taking into account \eqref{formule - d} and the second axiom of Definition \ref{dorfman connection}, we have that
for any $H= (H_0,H_1,\ldots,H_{[\frac{p}{2}]}) \in \mathfrak{D}_{m,m-1}^p(\mathcal{E}; \Gamma(B))$, its covariant derivative $d^{\nabla}H \in \mathfrak{D}_{m+1, m}^{p+1}(\mathcal{E}; \Gamma(B))$ is given by the $([\frac{p+1}{2}]+1)$-tuple $((d^{\nabla}H)_0, \big(d^{\nabla}H)_1, \ldots, (d^{\nabla}H)_{[\frac{p+1}{2}]})$ with
\begin{eqnarray*}\label{formula - d-napla}
\lefteqn{(d^{\nabla}H)_k(e_1,\ldots,e_{p+1-2k};f_1,\ldots,f_k) \, =} \nonumber \\
& & \sum\limits_{\mu=1}^k H_{k-1}(d_E f_{\mu}, e_1,\ldots,e_{p+1-2k}; f_1,\ldots, \hat{f}_{\mu},\ldots, f_k) \nonumber \\
& & + \, \sum\limits_{i=1}^{p+1-2k}(-1)^{i-1} \nabla_{e_i}(H_k(e_1,\ldots, \hat{e}_i,\ldots,e_{p+1-2k}; f_1,\ldots,f_k)) \nonumber \\
& & + \,\sum\limits_{i<j} (-1)^i H_k(e_1,\ldots, \hat{e}_i,\ldots, \hat{e}_j, \lcf e_i,e_j\rcf,e_{j+1},\ldots,e_{p+1-2k}; f_1,\ldots,f_k).
\end{eqnarray*}
For each $e_i \in \mathcal{E}$, $\nabla_{e_i}$ acts as a derivation on $H_k(e_1,\ldots, \hat{e}_i,\ldots,e_{p+1-2k}; f_1,\ldots,f_k)$ and it is a first order differential operator with respect to $e_i$. Furthermore, $\lcf e_i,e_j\rcf$ is a first order differential operator on the first item.
Hence, the order of the differential operators of the tuple $((d^{\nabla}H)_0, \big(d^{\nabla}H)_1, \ldots, (d^{\nabla}H)_{[\frac{p+1}{2}]})$ is increased by $1$ in each
argument with respect to the order of the operators of $H= (H_0,H_1,\ldots,H_{[\frac{p}{2}]})$ on the same argument.

\vspace{2mm}
\noindent
Elements of the space $(d^{\nabla})^m(\mathcal{C}^s(\mathcal{E}; \Gamma(B)))$, $m+s = p$, are elements of $\mathfrak{D}_{m, m-1}^p(\mathcal{E}; \Gamma(B))$. One simply applies successively ($m$ times)
the Dorfman covariant derivative $d^{\nabla}$ on $\mathcal{C}^s(\mathcal{E}; \Gamma(B))$. For $\omega \otimes b \in \mathcal{C}^s(\mathcal{E}; \Gamma(B))$, it is $(d^{\nabla})^m(\omega \otimes b) = \big(((d^{\nabla})^m(\omega \otimes b))_0,\,\ldots,\,((d^{\nabla})^m(\omega \otimes b))_{[\frac{p}{2}]}\big)$ with
\begin{linenomath*} \begin{equation*}
((d^{\nabla})^m(\omega \otimes b))_k : \mathcal{E}^{\otimes^{p-2k}}\otimes S^k\Omega^1 \to \Gamma(B).
\end{equation*} \end{linenomath*}
By induction on $m\in \N^\ast$, one can show that for any $k=0, \ldots, [\frac{p}{2}]$,
\begin{linenomath*} \begin{equation*}
((d^{\nabla})^m(\omega \otimes b))_k : \mathcal{E}^{\otimes^{p-2k}} \to \mathfrak{Diff}_{m-2k}(S_{\mathcal{R}}^k \Omega^1; \Gamma(B))
\end{equation*} \end{linenomath*}
is a differential operator of order $m-2k$ on the first $p-2k-1$ arguments and of $(m-2k-1)$-order on the last argument.

\vspace{1mm}
\noindent
Moreover, we may also prove by induction on $m\in \N$ that
\begin{itemize}
\item
if $m=2r$, $r\in \N$, $(d^{\nabla})^{2r}(\omega \otimes b) = \omega \cdot (d^{\nabla})^{2r}b$;
\item
if $m = 2r +1$, $r \in \N$, $(d^{\nabla})^{2r+1}(\omega \otimes b) = d\omega \cdot (d^{\nabla})^{2r}b + (-1)^p \omega \cdot (d^{\nabla})^{2r+1}(b)$.
\end{itemize}
As a result,  the elements of the spaces $(d^{\nabla})^m(\mathcal{C}^s(\mathcal{E};\Gamma(B)))$, $m+s =p$, are either of type $\omega\cdot (d^\nabla)^{m}b$
or of type $d\omega \cdot (d^\nabla)^{m-1}b$, with $\omega \in \mathcal{C}^s(\mathcal{E}; \mathcal{R})$.

\vspace{2mm}
\noindent
Extend the operators \eqref{operator i-f-e-D} to operators on $\mathfrak{D}_{\star, \, \star - 1}^\bullet(\mathcal{E};\Gamma(B))$ by
\begin{eqnarray}\label{generalization - i - f}
i_f : \mathfrak{D}_{m,\, m-1}^p(\mathcal{E}; \Gamma(B)) & \to  & \mathfrak{D}_{m,\, m-1}^{p-2}(\mathcal{E}; \Gamma(B))  \nonumber \\
                 \omega \otimes b & \mapsto & (i_f \omega)\otimes b
\end{eqnarray}
and
\begin{eqnarray}\label{generalization - i - e}
i_e : \mathfrak{D}_{m,\, m-1}^p(\mathcal{E};\Gamma(B)) & \to  & \mathfrak{D}_{m,\, m-1}^{p-1}(\mathcal{E}; \Gamma(B))  \nonumber \\
                 \omega \otimes b & \mapsto & (i_e \omega)\otimes b.
\end{eqnarray}
More explicitly, for an element of type $\omega \cdot (d^{\nabla})^{m}b \in \mathfrak{D}_{m, \,m-1}^p(\mathcal{E}; \Gamma(B))$, with $\omega \in \mathcal{C}^s(\mathcal{E}; \mathcal{R})$ and $m+s = p$,
\begin{linenomath*} \begin{equation} \label{generalization - i - f - D}
i_f (\omega \cdot (d^{\nabla})^{m}b) = (i_f\omega)\cdot (d^{\nabla})^{m}b + \omega \cdot (i_f(d^{\nabla})^{m}b)
\end{equation} \end{linenomath*}
and
\begin{linenomath*} \begin{equation} \label{generalization - i - e - D}
i_e (\omega \cdot (d^{\nabla})^{m}b) = (i_e \omega) \cdot (d^{\nabla})^{m}b + (-1)^p \omega \cdot (i_e (d^{\nabla})^{m}b).
\end{equation} \end{linenomath*}

\vspace{2mm}
\noindent
The commutator (\ref{commutator in C}) naturally extends to the space of graded endomorphisms of $\mathfrak{D}(\mathcal{E}; \Gamma(B))$. This way we obtain the operators
\begin{linenomath*} \begin{equation} \label{ext-Lie}
\nabla_e= \{i_e, d^\nabla\} = i_e\circ d^\nabla + d^\nabla \circ i_e,
\end{equation} \end{linenomath*}
\begin{linenomath*} \begin{equation} \label{ext-Lie-f}
\mathcal{L}_f^{\nabla} = \{i_f, d^\nabla \} = i_f \circ d^\nabla - d^\nabla \circ i_f,
\end{equation} \end{linenomath*}
satisfying the identities
\begin{linenomath*} \begin{equation*}
\{i_{e_1},i_{e_2}\} = i_{e_1}\circ i_{e_2} + i_{e_2}\circ i_{e_1} = i_{-\langle e_1,e_2 \rangle},
\end{equation*} \end{linenomath*}
\begin{linenomath*} \begin{equation} \label{ext-br}
\{\nabla_{e_1}, i_{e_2}\} = \nabla_{e_1} \circ i_{e_2} - i_{e_2} \circ \nabla_{e_1} = i_{\lcf e_1, e_2\rcf}.
\end{equation} \end{linenomath*}

\vspace{2mm}
\noindent
\begin{proposition}
The following identities hold for elements of $\mathfrak{D}_{m, \, m-1}^p(\mathcal{E}; \Gamma(B))$:
\begin{linenomath*} \begin{equation} \label{Lie derivatives connection}
\nabla_e (\omega \otimes b) = (\mathcal{L}_e\omega)\otimes b + \omega \cdot \nabla_e b \quad \quad \mathrm{and} \quad \quad \mathcal{L}_f^{\nabla}(\omega \otimes b) = (i_{d_E f}\omega)\otimes b.
\end{equation} \end{linenomath*}
\end{proposition}
\begin{proof}
Let $\omega \otimes b \in \mathfrak{D}_{m, \, m-1}^p(\mathcal{E};\Gamma(B))$. Then
\begin{eqnarray*}
\nabla_e (\omega \otimes b) & \stackrel{(\ref{ext-Lie})}{=} & (i_e\circ d^\nabla)(\omega \otimes b)+(d^\nabla \circ i_e)(\omega \otimes b) \nonumber \\
 & \stackrel{(\ref{Leibniz - d - napla})}{=} & i_e(d\omega\otimes b + (-1)^p\omega\cdot d^{\nabla}b) + d(i_e\omega) \otimes b + (-1)^{p-1}(i_e\omega)\cdot d^{\nabla}b \nonumber \\
 & \stackrel{(\ref{Leibniz - i - e}),(\ref{generalization - i - e}), (\ref{generalization - i - e - D})}{=} & (i_ed\omega + d i_e\omega)\otimes b + (-1)^p(i_e\omega)\cdot d^{\nabla}b + (-1)^{2p}\omega \otimes (i_e d^{\nabla}b) \nonumber \\
 & & + \,(-1)^{p-1}(i_e\omega)\cdot d^{\nabla}b \nonumber \\
 & \stackrel{(\ref{Lie derivatives})}{=} & \mathcal{L}_e\omega \otimes b + \omega \cdot \nabla_eb.
\end{eqnarray*}
Also,
\begin{eqnarray*}
\mathcal{L}_f^{\nabla}(\omega \otimes b) & \stackrel{(\ref{ext-Lie-f})}{=} & (i_f \circ d^\nabla - d^\nabla \circ i_f)(\omega \otimes b) \nonumber \\
& \stackrel{(\ref{Leibniz - d - napla})}{=} & i_f(d\omega \otimes b + (-1)^p \omega\cdot d^{\nabla}b) - d^\nabla((i_f \omega)\otimes b) \nonumber \\
& \stackrel{(\ref{Leibniz - d - napla}), (\ref{generalization - i - f}), (\ref{generalization - i - f - D})}{=} & (i_fd\omega)\otimes b + (-1)^p(i_f \omega)\cdot d^{\nabla}b - d(i_f\omega)\otimes b - (-1)^{p-2}(i_f\omega)\cdot d^{\nabla}b \nonumber \\
& \stackrel{(\ref{Lie derivatives})}{=} & \mathcal{L}_f\omega \otimes b \nonumber \\
& = & (i_{d_E f}\omega)\otimes b.
\end{eqnarray*}
\end{proof}

\begin{lemma}
The map
\begin{linenomath*} \begin{equation*}
(d^\nabla)^2 : \mathfrak{D}_{m, \, m-1}^p(\mathcal{E};\Gamma(B)) \to \mathfrak{D}_{m+2, \, m+1}^{p+2}(\mathcal{E}; \Gamma(B))
\end{equation*} \end{linenomath*}
is $\C$-linear on the sections of $B$.
\end{lemma}
\begin{proof}
For $\omega \otimes b \in \mathfrak{D}_{m, \, m-1}^p (\mathcal{E};\Gamma(B))$ and $f\in \mathcal{R}$,
\begin{eqnarray*}
(d^{\nabla})^2(\omega \otimes (fb)) & = & \omega\cdot (d^{\nabla})^2(fb) = \omega\cdot d^{\nabla}(d_Ef\otimes b + f d^{\nabla}b) \nonumber \\
& = &\omega \cdot \big(d(d_Ef)\otimes b - d_Ef\cdot d^{\nabla}b + d_Ef\cdot d^{\nabla}b + f (d^{\nabla})^2b\big) \nonumber \\
& = & \omega\cdot \big(f (d^{\nabla})^2b\big).
\end{eqnarray*}
\end{proof}

\vspace{3mm}
\noindent
It follows that the map $(d^{\nabla})^2 : \Gamma(B) \to (d^{\nabla})^2(\Gamma(B)) \subset \mathfrak{D}_{2,\, 1}^2(\mathcal{E}; \Gamma(B))$ is an $\Gamma(\mathrm{End}(B))$-valued element
of $\mathfrak{D}_{2, \, 1}^2(\mathcal{E}; \mathcal{R})$. We identify it with an element $R^{\nabla}=(R^{\nabla}_0, R^{\nabla}_1)$ of $\mathfrak{D}_{2,\,1}^2(\mathcal{E}; \Gamma(\mathrm{End}(B)))$.

\begin{definition}\label{definition curvature}
For any $E$-Dorfman connection $\nabla$ on a predual $B\to M$ of $E$,
the element $(d^\nabla)^2 \in \mathfrak{D}_{2,\,1}^2(\mathcal{E}; \Gamma(\mathrm{End}(B)))$ is called the \emph{curvature of} $\nabla$.
An $E$-Dorfman connection $\nabla$ on $B$ whose curvature $(d^\nabla)^2$ is identically zero is called \emph{flat} or \emph{$E$-Dorfman action on $B$}. In this case, $B$ is also called \emph{$E$-Dorfman module}.

\vspace{1mm}
\noindent
If $\nabla$ is flat, $d^\nabla$ is a differential on $\mathfrak{D}(\mathcal{E}; \Gamma(B))$ and we denote by $\mathfrak{H}^p(\mathcal{E}; \Gamma(B))$ the $p$-cohomology group of the cochain complex
$(\mathfrak{D}(\mathcal{E}; \Gamma(B)),\,d^\nabla)$.
\end{definition}

\begin{proposition}\label{prop-curvature}
The curvature $(d^\nabla)^2 : \Gamma(B) \to \mathfrak{D}_{2,\,1}^2(\mathcal{E}; \Gamma(B))$ of an $E$-Dorfman connection $\nabla$ on a predual vector bundle $B$ of $E$ satisfies the following identities:
\begin{linenomath*} \begin{equation*}\label{curvature formula e}
i_{e_2}\circ i_{e_1}((d^\nabla)^2b) = \nabla_{e_1}\nabla_{e_2}b - \nabla_{e_2}\nabla_{e_1}b - \nabla_{\lcf e_1, e_2\rcf}b,
\end{equation*} \end{linenomath*}
\begin{linenomath*} \begin{equation*}
i_f((d^\nabla)^2 b) = \nabla_{d_Ef}b.
\end{equation*} \end{linenomath*}
Furthermore, the restriction of $(d^\nabla)^2$ on $\mathrm{Im}d_B$ vanishes.
\end{proposition}
\begin{proof}
Let $b \in \Gamma(B)$, $e_1, e_2 \in \Gamma(E)$, and $f \in \mathcal{R} \cong \C$. Then one has
\begin{eqnarray*}
i_{e_2}\circ i_{e_1}((d^\nabla)^2 b) & = & i_{e_2}\circ (i_{e_1}\circ d^\nabla)(d^\nabla b) \nonumber \\
& \stackrel{(\ref{ext-Lie})}{=} & i_{e_2}\circ (\nabla_{e_1} - d^\nabla \circ i_{e_1})(d^\nabla b) \nonumber \\
& = & (i_{e_2}\circ \nabla_{e_1})(d^\nabla b) - (i_{e_2} \circ d^\nabla)\nabla_{e_1}b \nonumber \\
& \stackrel{(\ref{ext-br})(\ref{ext-Lie})}{=} & (\nabla_{e_1} \circ i_{e_2} - i_{\lcf e_1, e_2\rcf})(d^\nabla b) - (\nabla_{e_2} - d^\nabla \circ i_{e_2})\nabla_{e_1}b \nonumber \\
& = & (\nabla_{e_1}\nabla_{e_2} - \nabla_{e_2}\nabla_{e_1} - \nabla_{\lcf e_1, e_2\rcf})(b),
\end{eqnarray*}
which is the well known formula of curvature. We also write $i_{e_2}\circ i_{e_1}((d^\nabla)^2 b) = R_0^{\nabla}(e_1,e_2)b$, and so
\begin{linenomath*} \begin{equation} \label{curvature formula classic}
R_0^{\nabla}(e_1,e_2)b = \nabla_{e_1}\nabla_{e_2}b - \nabla_{e_2}\nabla_{e_1}b - \nabla_{\lcf e_1, e_2\rcf}b.
\end{equation} \end{linenomath*}
Since $i_f$ is of degree $-2$, it is
\begin{eqnarray*}
i_f((d^\nabla)^2 b) & = & (i_f \circ d^\nabla)(d^\nabla b) \stackrel{(\ref{ext-Lie-f})}{=} (\mathcal{L}_f^{\nabla}+ d^\nabla \circ i_f)(d^\nabla b) \nonumber \\
& = & \mathcal{L}_f^{\nabla} (d^\nabla b) \stackrel{\eqref{Lie derivatives connection}}{=} i_{d_E f}(d^\nabla b) = \nabla_{d_Ef}b,
\end{eqnarray*}
and so
\begin{linenomath*} \begin{equation} \label{curvature formula f}
R_1^{\nabla}(f) b = \nabla_{d_Ef}b.
\end{equation} \end{linenomath*}
Clearly
\begin{linenomath*} \begin{equation} \label{R_0-R_1}
R_0(e_1,e_2) + R_0(e_2,e_1) = -R_1(\langle e_1, e_2\rangle),
\end{equation} \end{linenomath*}
$R_0^{\nabla}(e_1,e_2)d_Bg =0$ and $R_1^{\nabla}(f)d_Bg =0$, for any $g \in \C$.

\vspace{2mm}
\noindent
One can justify the claim that $(d^\nabla)^2$ is an element of $\mathfrak{D}_{2, \,1}^2(\mathcal{E}; \Gamma(\mathrm{End}(B)))$,
by computing the symbols of $R^{\nabla}_0$ and $R^{\nabla}_1$ (Definition \ref{def-dif op - symbol}). After a straightforward calculation, we get that, for any $e_1, e_2 \in \Gamma(E)$, $b\in \Gamma(B)$, $\alpha \in \Omega^1$ and $f\in \C$:
\begin{enumerate}
\item
The symbol $\sigma_1(R_0^\nabla)(f)$ of $R_0^\nabla$ is a $1$-order differential operator in the first argument since
\begin{eqnarray*}
\sigma_1(R_0^\nabla)(f)(e_1,e_2)b & = & R_0^\nabla(fe_1,e_2)b - fR_0^\nabla(e_1,e_2)b \nonumber \\
& = &  \big(\lan \lcf e_2, e_1 \rcf, b \ran + \lan e_1, \nabla_{e_2}b\ran - \rho(e_2)(\lan e_1,b \ran) \big)d_Bf \nonumber \\
& &  -\, \lan d_E\langle e_1, e_2\rangle,\, b \ran d_Bf - \nabla_{\langle e_1, e_2\rangle d_Ef} b \nonumber \\
&\stackrel{\eqref{B-ast-connection}}{=} & -\langle D^\ast_b e_1, e_2\rangle d_Bf - \nabla_{\langle e_1, e_2\rangle d_Ef}b\nonumber \\
& = & -\langle D^\ast_b e_1, e_2\rangle d_Bf - \langle e_1, e_2\rangle\nabla_{ d_Ef}b - \varrho(b)(f)d_B\langle e_1, e_2\rangle.
\end{eqnarray*}
\item
The symbol $\sigma_2(R_0^\nabla)(f)$ of $R_0^\nabla$ is a $0$-order (i.e. $\C$-linear) differential operator in the second argument since
\begin{eqnarray*}
\sigma_2(R_0^\nabla)(f)(e_1,e_2)b & = & R_0^\nabla(e_1,fe_2)b -fR_0^\nabla(e_1,e_2)b \nonumber \\
 & =  & \big(\rho(e_1)\lan e_2,b \ran - \lan \lcf e_1, e_2 \rcf, b \ran - \lan e_2, \nabla_{e_1}b\ran  \big)d_Bf \nonumber \\
 & \stackrel{\eqref{b-connection}}{=}&- \langle D_be_1, e_2\rangle d_Bf.
\end{eqnarray*}
\item
The symbol $s_1(R_1^\nabla)(f)$ of $R_1^\nabla$ is a first order differential operator since
\begin{eqnarray*}
s_1(R_1^\nabla)(f)(\alpha)b & = & R_1^\nabla(f\alpha)b - f R_1^\nabla(\alpha)b \nonumber \\
& = &  \nabla_{f{g^\flat}^{-1}(\rho^\ast \alpha)}b -f \nabla_{{g^\flat}^{-1}(\rho^\ast \alpha)}b \nonumber \\
& = &  \lan {g^\flat}^{-1}(\rho^\ast \alpha), b\ran d_Bf.
\end{eqnarray*}
\end{enumerate}
\end{proof}

\begin{example}
\emph{Applying formulas \eqref{curvature formula classic} and \eqref{curvature formula f} to the $E$-Dorfman connection $\triangledown$ defined in Example \ref{dorfman on T*M}, we find that it is flat. Thus, $T^\ast M$ is an $E$-module.}
\end{example}

\begin{remark}\label{rem-napla as loday morph}
\emph{Using the operator $\nabla : \Gamma(E) \to \Gamma(\mathbb{A}(B))$ defined in item \ref{napla as dif op} of Remarks \ref{rem on dorf con}, the flatness condition $R^{\nabla}=(R_0^{\nabla},R_1^{\nabla}) = (0,0)$ is written as}
\[\nabla_{\lcf e_1, e_2\rcf} = [\nabla_{e_1}, \nabla_{e_2}] \quad \quad \mathrm{and} \quad \quad \nabla_{d_Ef}=0.\]
\emph{Hence flatness implies that $\nabla$ is a Loday algebra homomorphism from $(\Gamma(E), \lcf \cdot, \cdot \rcf)$ to $(\Gamma(\mathbb{A}(B)), [ \cdot, \cdot ])$ vanishing on $\mathrm{Im}d_E$.}
\end{remark}

\subsection{Induced connections and Bianchi identity}
\subsubsection{Connection on the dual bundle}
An $E$-Dorfman connection $\nabla$ of $E$  on a predual vector bundle $B$  induces a map
\begin{linenomath*} \begin{equation}\label{dual Dorfman connection}
\nabla^{\ast} : \Gamma(E) \times \Gamma(B^{\ast}) \to \Gamma(B^{\ast}),
\end{equation} \end{linenomath*}
$B^{\ast}$ being the dual vector bundle of $B$,
completely characterized by the relation
\begin{linenomath*} \begin{equation*}
\rho(e)\langle b^{\ast}, b\rangle = \langle \nabla^{\ast}_eb^{\ast}, b\rangle + \langle b^{\ast}, \nabla_e b\rangle,
\end{equation*} \end{linenomath*}
where $\langle \cdot , \cdot \rangle$ denotes the duality pairing between $\Gamma(B^\ast)$ and $\Gamma(B)$, $e\in \Gamma(E)$, $b\in \Gamma(B)$ and $b^{\ast}\in \Gamma(B^{\ast})$. By the properties of $\nabla$ (see Definition \ref{dorfman connection}) one has that, for any $f\in \C$,
\begin{enumerate}
\item
$\nabla^{\ast}_{fe}b^{\ast} = f\nabla^{\ast}_eb^\ast - \langle b^\ast, d_Bf \rangle \lan e, \cdot \ran$, and
\item
$\nabla^{\ast}_{e}(fb^\ast) = f\nabla^{\ast}_eb^\ast + \rho(e)(f)b^\ast$.
\end{enumerate}

\vspace{2mm}
\noindent
Note that $\nabla^\ast$ is not a Dorfman connection as it does not satisfy the conditions of Definition  \ref{dorfman connection}. It is a nonlinear connection constructed by a Dorfman connection,
which we will call the \textit{dual connection of the Dorfman connection $\nabla$}. When the Dorfman connection is clear from the context or notation, we will simply
refer to it as the \textit{dual connection}. The curvature $R^{\nabla^\ast} = (R_0^{\nabla^\ast}, R_1^{\nabla^\ast})$ of $\nabla^\ast$ is then defined by the relations
\begin{linenomath*} \begin{equation} \label{dual Dorfman curvature}
\langle R_0^{\nabla^\ast}(e_1,e_2)b^\ast, b\rangle + \langle b^\ast, R_0^{\nabla}(e_1,e_2)b\rangle = 0 \quad \mathrm{and} \quad \langle R_1^{\nabla^\ast}(f)b^\ast, b\rangle + \langle b^\ast, R_1^{\nabla}(f)b\rangle = 0.
\end{equation} \end{linenomath*}

\noindent
Clearly, if $\nabla$ is a flat $E$-Dorfman connection on $B$, then $\nabla^{\ast}$ is also a flat connection on $B^\ast$.

\subsubsection{Connection on the endomorphism bundle}
As in the classical case, $\nabla$ induces a (nonlinear) connection on any tensor bundle constructed from $B$.
In particular, the pair $(\nabla^\ast, \nabla)$ induces a tensor product nonlinear connection $\widetilde{\nabla} = \nabla^\ast \otimes \mathrm{id}_B + \mathrm{id}_{B^\ast}\otimes \nabla$ of $\Gamma(E)$ on $\Gamma(B^\ast \otimes B) \cong \Gamma(\mathrm{End}(B))$ by
\begin{linenomath*} \begin{equation*}
\widetilde{\nabla}_e (b^\ast \otimes b) = \nabla^\ast_e b^\ast \otimes b + b^\ast \otimes \nabla_eb.
\end{equation*} \end{linenomath*}
It is then a simple calculation to check that the connection $\widetilde{\nabla}$ has the following properties:
\begin{enumerate}
\item
$\widetilde{\nabla}_{fe}(b^\ast \otimes b) = f \widetilde{\nabla}_{e}(b^\ast \otimes b) - \langle b^\ast, d_Bf \rangle \lan e, \cdot \ran \otimes b + \lan e, b \ran b^\ast \otimes d_Bf$,
\item
$\widetilde{\nabla}_{e}f(b^\ast \otimes b) = f\widetilde{\nabla}_{e}(b^\ast \otimes b) + \rho(e)(f)(b^\ast \otimes b)$,
\item
$\widetilde{\nabla}_{e}(b^\ast \otimes d_Bf) = \nabla^\ast_e b^\ast \otimes d_Bf + b^\ast \otimes d_B(\mathcal{L}_{\rho(e)}f)$.
\end{enumerate}

\noindent
A simple computation based on these properties provides the following, that can also serve as an alternative definition of $\widetilde{\nabla}$.

\begin{proposition}
Let $e\in \Gamma(E)$ and $\tau \in \Gamma(\mathrm{End}(B))$. Then $\widetilde{\nabla}$ satisfies
\begin{linenomath*} \begin{equation} \label{definition widetilde nabla}
\widetilde{\nabla}_{e}\tau = [\nabla_e , \tau] = \nabla_e \circ \tau - \tau \circ \nabla_e.
\end{equation} \end{linenomath*}
\end{proposition}

\vspace{2mm}
\noindent
Set $\mathfrak{D}(\mathcal{E}; \Gamma(\mathrm{End}(B))) = \big(\mathfrak{D}_{m,\, m-1}^p(\mathcal{E}; \Gamma(\mathrm{End}(B)))\big)_{p\in \N,\, m\in \N^\ast}$. The space $ \mathfrak{D}_{m,\, m-1}^p(\mathcal{E}; \Gamma(\mathrm{End}(B)))$ is the space of
$([\frac{p}{2}]+1)$-tuples $\Phi= (\Phi_0,\Phi_1,\ldots,\Phi_{[\frac{p}{2}]})$ of homomorphisms
\begin{linenomath*} \begin{equation*}
\Phi_k : \mathcal{E}^{\otimes^{p-2k}}\otimes S^k\Omega^1 \to \Gamma(\mathrm{End}(B))
\end{equation*} \end{linenomath*}
characterized by similar conditions to those of operators $\bar{\omega}_k \in \mathfrak{D}_{m, \, m-1}^p(\mathcal{E}; \mathcal{R})$ and $H_k \in \mathfrak{D}_{m, \, m-1}^p(\mathcal{E}; \Gamma(B))$ (see, respectively, sections \ref{cohomology-CD-Diff} and \ref{section-curvature}).

\vspace{1mm}
\noindent
Note that for $p=1$, it is $\Phi = (\Phi_0)$ and the unique argument of the map $\Phi_0 : \mathcal{E} \to \Gamma(\mathrm{End}(B))$ is considered as \emph{first} argument. Hence, we denote by $\mathfrak{D}^1_{m}(\mathcal{E};\Gamma(\mathrm{End}(B)))$
the space of differential operators on $\mathcal{E}$ of order at most $m$, $m\in \N^\ast$, with values in $\Gamma(\mathrm{End}(B))$.

\vspace{2mm}
\noindent
The connection $\widetilde{\nabla}$ on $\mathrm{End}(B)$ defines a covariant derivation operator
\begin{linenomath*} \begin{equation*}
d^{\widetilde{\nabla}} : \Gamma(\mathrm{End}(B)) \to \mathfrak{D}_1^1(\mathcal{E}; \Gamma(\mathrm{End}(B)))
\end{equation*} \end{linenomath*}
such that, for any $f\in \C$ and $\tau \in \Gamma(\mathrm{End}(B))$, it is
\begin{linenomath*} \begin{equation*}
d^{\widetilde{\nabla}}(f\tau) = d_Ef\otimes \tau + fd^{\widetilde{\nabla}}\tau.
\end{equation*} \end{linenomath*}
This extends uniquely to an operator of degree $+1$, denoted also by $d^{\widetilde{\nabla}}$, on the space $\mathfrak{D}(\mathcal{E}; \Gamma(\mathrm{End}(B)))$, namely
\begin{linenomath*} \begin{equation*}
d^{\widetilde{\nabla}} : \mathfrak{D}_{m,\,m-1}^p(\mathcal{E}; \Gamma(\mathrm{End}(B))) \to \mathfrak{D}_{m+1, \, m}^{p+1}(\mathcal{E};\Gamma(\mathrm{End}(B))).
\end{equation*} \end{linenomath*}
More precisely, the image of an element $\Phi = (\Phi_0,\Phi_1,\ldots,\Phi_{[\frac{p}{2}]}) \in \mathfrak{D}_{m,\, m-1}^p(\mathcal{E}; \Gamma(\mathrm{End}(B)))$,
is $d^{\widetilde{\nabla}}\Phi = ((d^{\widetilde{\nabla}}\Phi)_0, (d^{\widetilde{\nabla}}\Phi)_1, \ldots, (d^{\widetilde{\nabla}}\Phi)_{[\frac{p+1}{2}]}) \in \mathfrak{D}_{m+1,\, m}^{p+1}(\mathcal{E}; \Gamma(\mathrm{End}(B))$, where
\begin{eqnarray}\label{formula - d-widetildenapla}
\lefteqn{(d^{\widetilde{\nabla}}\Phi)_k(e_1,\ldots,e_{p+1-2k};f_1,\ldots,f_k) \, =} \nonumber \\
& & \sum\limits_{\mu=1}^k \Phi_{k-1}(d_E f_{\mu}, e_1,\ldots,e_{p+1-2k}; f_1,\ldots, \hat{f}_{\mu},\ldots, f_k) \nonumber \\
& & + \, \sum\limits_{i=1}^{p+1-2k}(-1)^{i-1} \widetilde{\nabla}_{e_i}(\Phi_k(e_1,\ldots, \hat{e}_i,\ldots,e_{p+1-2k}; f_1,\ldots,f_k)) \nonumber \\
& & + \,\sum\limits_{i<j} (-1)^i \Phi_k(e_1,\ldots, \hat{e}_i,\ldots, \hat{e}_j, \lcf e_i,e_j\rcf,e_{j+1},\ldots,e_{p+1-2k}; f_1,\ldots,f_k).
\end{eqnarray}

\vspace{1mm}
\noindent
The curvature $R^{\widetilde{\nabla}} = (R_0^{\widetilde{\nabla}}, R_1^{\widetilde{\nabla}})$ of $\widetilde{\nabla}$ is given, for any $e_1,e_2 \in \Gamma(E)$, $f\in \C$, and $\tau \in \Gamma(\mathrm{End}(B))$, by
\begin{linenomath*} \begin{equation} \label{curvature napla tilde}
R_0^{\widetilde{\nabla}}(e_1,e_2)\tau = R_0^{\nabla}(e_1,e_2)\circ \tau - \tau \circ R_0^{\nabla}(e_1,e_2) \quad \mathrm{and} \quad R_1^{\widetilde{\nabla}}(f)\tau = R_1^{\nabla}(f)\circ \tau - \tau \circ R_1^{\nabla}(f).
\end{equation} \end{linenomath*}
Hence, if $\nabla$ is a flat $E$-Dorfman connection on $B$, then $\widetilde{\nabla}$ is a flat connection on $B^\ast \otimes B \cong \mathrm{End}(B)$. In this case, $d^{\widetilde{\nabla}}$ is a differential and $\big(\mathfrak{D}(\mathcal{E}; \Gamma(\mathrm{End}(B))),d^{\widetilde{\nabla}}\big)$ is a cochain complex.

\begin{proposition}[Bianchi identity]\label{Bianchi}
Let $(E, B, \nabla, R^{\nabla})$ be as above. Then
\begin{linenomath*} \begin{equation*}
d^{\widetilde{\nabla}} (R^{\nabla}) = 0.
\end{equation*} \end{linenomath*}
\end{proposition}
\begin{proof}
Consider the curvature $R^{\nabla} = (R^{\nabla}_0, R^{\nabla}_1) \in \mathfrak{D}_{2,\,1}^2(\mathcal{E};\Gamma(\mathrm{End}(B)))$. Its image through $d^{\widetilde{\nabla}}$ is
$d^{\widetilde{\nabla}} (R^{\nabla}) = ((d^{\widetilde{\nabla}} (R^{\nabla}))_0, (d^{\widetilde{\nabla}} (R^{\nabla}))_1) \in \mathfrak{D}_{3,\,2}^3(\mathcal{E};\Gamma(\mathrm{End}(B)))$ and for any $e_1, e_2, e_3 \in \Gamma(E)$ it is
\begin{eqnarray}\label{bianchi - e}
(d^{\widetilde{\nabla}} (R^{\nabla}))_0(e_1,e_2,e_3) & \stackrel{\eqref{formula - d-widetildenapla}}{=} & \widetilde{\nabla}_{e_1}(R_0^{\nabla}(e_2,e_3)) - \widetilde{\nabla}_{e_2}(R_0^{\nabla}(e_1,e_3)) + \widetilde{\nabla}_{e_3}(R_0^{\nabla}(e_1,e_2)) \nonumber \\
& & - R_0^{\nabla}(\lcf e_1,e_2 \rcf, e_3) - R_0^{\nabla}(e_2,\lcf e_1,e_3 \rcf) + R_0^{\nabla}(e_1, \lcf e_2,e_3 \rcf) \nonumber \\
&\stackrel{\eqref{definition widetilde nabla}}{=}& \nabla_{e_1}\circ R_0^{\nabla}(e_2,e_3) - R_0^{\nabla}(e_2,e_3)\circ\nabla_{e_1} - \nabla_{e_2}\circ R_0^{\nabla}(e_1,e_3) \nonumber \\
& & +\, R_0^{\nabla}(e_1,e_3)\circ\nabla_{e_2} + \nabla_{e_3}\circ R_0^{\nabla}(e_1,e_2) - R_0^{\nabla}(e_1,e_2)\circ\nabla_{e_3} \nonumber \\
& & -\, R_0^{\nabla}(\lcf e_1,e_2 \rcf, e_3) - R_0^{\nabla}(e_2,\lcf e_1,e_3 \rcf) + R_0^{\nabla}(e_1, \lcf e_2,e_3 \rcf) \nonumber \\
& = & 0.
\end{eqnarray}
For the last equation use the curvature expression (\ref{curvature formula classic}) and the fact that the bracket $\lcf \cdot, \cdot \rcf$ verifies the Jacobi identity (\ref{Jacobi-Courant}). Similarly, for $(d^{\widetilde{\nabla}} (R^{\nabla}))_1$ one gets
\begin{eqnarray}\label{bianchi - f}
(d^{\widetilde{\nabla}} (R^{\nabla}))_1(e;f) & \stackrel{\eqref{formula - d-widetildenapla}}{=} & R_0^{\nabla}(d_Ef, e) + \widetilde{\nabla}_{e}(R_1^{\nabla}(f)) \nonumber \\
& \stackrel{\eqref{definition widetilde nabla}}{=}& R_0^{\nabla}(d_Ef, e) + \nabla_e \circ R_1^{\nabla}(f) - R_1^{\nabla}(f) \circ \nabla_e \nonumber \\
& \stackrel{\eqref{curvature formula classic}, \eqref{curvature formula f}}{=} & \nabla_{d_Ef}\nabla_e - \nabla_e\nabla_{d_Ef} - \nabla_{\lcf d_Ef, e \rcf} + \nabla_e\nabla_{d_Ef} - \nabla_{d_Ef}\nabla_e \nonumber \\
& = & 0.
\end{eqnarray}
\end{proof}

\subsubsection{Connections on modules of differential operators}
Given predual vector bundles $B$ and $F$ of $E$ with Dorfman connections $\nabla^B$ and $\nabla^F$, respectively, we can assign covariant derivation laws to $\C$-modules obtained by $B$ and $F$ in a canonical way. For later use, we will discuss how to assign a covariant derivation law on the $\C$-module of differential operators acting on $\Gamma(B)$ and with values in $\Gamma(\mathrm{End}(F))$.

\vspace{1mm}
\noindent
Let $\mathfrak{D}^1_k(\Gamma(B); \Gamma(\mathrm{End}(F)))$ be the set of all differential operators $\psi : \Gamma(B)\to \Gamma(\mathrm{End}(F))$ of order at most $k\in \N$. It is a left $\C$-module with respect to the natural module structure of $\Gamma(\mathrm{End}(F))$; the operator $f\cdot \psi : \Gamma(B)\to \Gamma(\mathrm{End}(F))$, noted below by simply $f\psi$,  is defined by $(f \psi)(b)=f(\psi(b))$. Define a nonlinear $E$-connection $\widehat{\nabla}$ on $\mathfrak{D}^1_k(\Gamma(B); \Gamma(\mathrm{End}(F)))$ by setting, for every $e\in \Gamma(E)$, $\psi \in \mathfrak{D}^1_k(\Gamma(B); \Gamma(\mathrm{End}(F)))$ and $b\in \Gamma(B)$,
\begin{equation}\label{connection hat}
\widehat{\nabla}_e\psi (b) = \widetilde{\nabla}^F_e(\psi(b)) - \psi(\nabla_e^Bb).
\end{equation}
Nonlinearity is calculated directly by \eqref{connection hat}; for any $f\in \C$, we have

\begin{eqnarray}\label{formule-hat-symbol}
\widehat{\nabla}_{fe}\psi(b) & = & f \widehat{\nabla}_e\psi(b) + \lan e, \psi(b)(\cdot)\ran d_Ff - \lan e, \cdot \ran \psi(b)(d_Ff) \nonumber \\
& & - \, \sigma(\psi)(f)(\nabla_e^Bb)-\lan e, b \ran \psi (d_B f) - \sigma(\psi)(\lan e, b\ran)(d_Bf),
\end{eqnarray}
where $\sigma (\psi)$ denotes the symbol of the differential operator $\psi$, and
\[\widehat{\nabla}_ef\psi = f\widehat{\nabla}_e\psi + \rho(e)(f)\psi.\]
The connection $\widehat{\nabla}$ induces a covariant derivation law of degree $+1$:
\begin{eqnarray*}
d^{\widehat{\nabla}} : \mathfrak{D}_k^1(\Gamma(B); \Gamma(\mathrm{End}(F)))& \to &\mathfrak{D}^1_{{k+1}}(\mathcal{E} ; \mathfrak{D}_k^1(\Gamma(B); \Gamma(\mathrm{End}(F)))) \nonumber \\
\psi & \mapsto & d^{\widehat{\nabla}}\psi (e) : = \widehat{\nabla}_e\psi.
\end{eqnarray*}
For any $\psi \in \mathfrak{D}_k^1(\Gamma(B); \Gamma(\mathrm{End}(F)))$, the image $d^{\widehat{\nabla}}\psi$ is a differential operator of at most $(k+1)$-order in the $\mathcal{E}$-argument as can by justified by \eqref{formule-hat-symbol}:  $\nabla^B_\cdot b$ is a $1$-order differential operator in the $\mathcal{E}$-argument, and since $\psi$ is a differential operator of at most $k$-order, its symbol $\sigma(\psi)$ is a differential operator of at most $(k-1)$-order. Thus their composition $\sigma(\psi)(f)\circ \nabla_\cdot^Bb$ is a differential operator of at most $k$-order in the $\mathcal{E}$-argument. The other terms in \eqref{formule-hat-symbol} are differential operators of order less than $k$ in the $\mathcal{E}$-argument. Consequently our assertion is true. Moreover, for any $e\in \Gamma(E)$, $\widehat{\nabla}_e\psi : \Gamma(B) \to \Gamma(\mathrm{End}(F))$ is a differential operator of the same order as $\psi$, at most $k$. In fact, for any $f\in \C$ and $b\in \Gamma(B)$,
\begin{equation*}
\widehat{\nabla}_e\psi (fb)  =  f \widehat{\nabla}_e\psi (b) + \widehat{\nabla}_e \sigma(\psi)(f)(b) - \sigma(\psi)(\rho(e)(f))(b).
\end{equation*}
By induction on the order of $\psi$ on $B$-argument, we establish our claim.

\vspace{1mm}
\noindent
As in previous cases, $d^{\widehat{\nabla}}$ extends uniquely to an operator of degree $+1$, denoted also by $d^{\widehat{\nabla}}$, on the space of $\mathfrak{D}^1_k(\Gamma(B); \Gamma(\mathrm{End}(F)))$-valued differential operators, $k\in \N$. In particular, because of the Proposition \ref{proposition-symbols},
\begin{eqnarray*}
\lefteqn{ d^{\widehat{\nabla}} : \mathfrak{D}^p_{m, m-1}(\mathcal{E}; \mathfrak{D}_k^1(\Gamma(B); \Gamma(\mathrm{End}(F)))) \quad \rightarrow } \\
& & \quad \mathfrak{D}^{p+1}_{\mathrm{max}\{m+1, k+1\},\, \mathrm{max}\{m, k+1\}}(\mathcal{E}; \mathfrak{D}_{k}^1(\Gamma(B); \Gamma(\mathrm{End}(F)))).
\end{eqnarray*}
Developing the usual calculation, we find that the curvature $R^{\widehat{\nabla}} = (R^{\widehat{\nabla}}_0, R^{\widehat{\nabla}}_1)$ of $\widehat{\nabla}$ is given, for any $e_1,e_2 \in \Gamma(E), b\in\Gamma(B)$, $f\in \C$, and $\psi \in \mathfrak{D}_k^1(\Gamma(B); \Gamma(\mathrm{End}(F)))$, by
\begin{eqnarray}\label{courbure-modules}
\big(R^{\widehat{\nabla}}_0(e_1,e_2)\psi\big)b & = & R_0^{\widetilde{\nabla}^F}(e_1,e_2)\psi (b) - \psi(R_0^{\nabla^B}(e_1,e_2)b), \nonumber \\
\big(R^{\widehat{\nabla}}_1(f)\psi\big)b & = & R_1^{\widetilde{\nabla}^F}(f)\psi (b) - \psi(R_1^{\nabla^B}(f)b).
\end{eqnarray}
From the above expression of $R^{\widehat{\nabla}}$ and taking into account \eqref{curvature napla tilde}, we get that, if $\nabla^B$ and $\nabla^F$ are flat, then $\widehat{\nabla}$ is also flat. In this case, $\big(\mathfrak{D}(\mathcal{E}; \mathfrak{D}_k^1(\Gamma(B); \Gamma(\mathrm{End}(F)))), \; d^{\widehat{\nabla}}\big)$ is a cochain complex and its $p$-cohomology group is denoted by $\mathfrak{H}^p(\mathcal{E}; \mathfrak{D}^1_k(\Gamma(B); \Gamma(\mathrm{End}(F))))$.

\subsection{Examples of Dorfman connections}\label{section examples}
\begin{example}\label{Dorfman E on E stand}
\emph{Let $(E, \lcf \cdot, \cdot\rcf, \langle \cdot, \cdot \rangle, \rho)$ be a Courant algebroid over a smooth manifold $M$ and $D : \Gamma(E)\times \Gamma(E) \to \Gamma(E)$ a linear $E$-connection on $E$, as defined in \cite{Al-Xu, CM, {Gualt-branes}}, such that $D_{d_Ef} = 0$, for any $f\in \C$. Consider the map $\nabla : \Gamma(E) \times \Gamma(E) \to \Gamma(E)$, given, for any pair $(e,e')$ of sections of $E$, by}
\begin{equation*}
\nabla_e e' = \lcf e, e' \rcf + D_{e'}e.
\end{equation*}
\emph{We can easily verify that $\nabla$ defines an $E$-Dorfman connection on $E$. The corresponding, via the Proposition \ref{prop - D - connection}, linear $E$-connection $D$ on $E$ is the given one.}
\end{example}

\begin{example}\label{example - mjl}
{\rm{This example is inspired by \cite[Example 4.2]{mjl} concerning a Dorfman connection of a dull algebroid on a vector bundle. Consider the standard Courant algebroid $E=TM\oplus T^\ast M$ (Example \ref{Standard Courant algebroid}) and a linear $TM$-connection $\triangle$ on $TM$\footnote{For a vector bundle $A\to M$, a linear $TM$-connection $\triangle$ on $A$ is an $\R$-bilinear map $\triangle : \Gamma(TM) \times \Gamma(A) \to \Gamma(A)$ such that:
(i) $\triangle_{fX}a = f\triangle_Xa$, (ii) $\triangle_Xfa = f\triangle_Xa + X(f)a$. It is also called \emph{Koszul connection} and it always exists \cite[p. 185]{mck}.}. Let $\triangle^\ast$ be its dual connection on $T^\ast M$. The map}}
\begin{linenomath*} \begin{equation*}
\nabla : \Gamma(E) \times \Gamma(E)  \to  \Gamma(E),
\end{equation*} \end{linenomath*}
{\rm{defined, for any $X+\zeta, Y+\eta \in \Gamma(E)$, by}}
\begin{linenomath*} \begin{equation*}
\nabla_{X+\zeta}(Y+\eta) = \triangle_X Y + (\mathcal{L}_X\eta + \langle \triangle^\ast\!_{\cdot}\zeta, Y\rangle),
\end{equation*} \end{linenomath*}
{\rm{defines an $E$-Dorfman connection on $(E, d_E, \langle \cdot, \cdot \rangle)$. Its dual connection $\nabla^\ast$ on $E^\ast = T^\ast M \oplus TM$ is given, for any $X+\zeta\in \Gamma(E)$ and $\eta + Y\in \Gamma(E^\ast)$, by}}
\begin{linenomath*} \begin{equation*}
\nabla^\ast_{X+\zeta}(\eta + Y) = (\triangle^\ast_X\eta - \triangle^\ast_Y\zeta) + (\triangle_X Y + \langle \triangle^\ast _X\cdot - \mathcal{L}_X\cdot , Y\rangle ).
\end{equation*} \end{linenomath*}
\end{example}

\begin{examples}[Regular Courant algebroids]{\rm{\cite{csx1}}}
{\rm{Let $(E, \lcf \cdot, \cdot\rcf, \langle \cdot, \cdot \rangle, \rho)$ be a \emph{regular Cou\-rant algebroid}, i.e. $F:= \rho(E)\subseteq TM$ is an integrable
distribution of constant rank on the base manifold $M$ and so defines a regular foliation of $M$. Then, $\ker \rho$ and its orthogonal $(\ker \rho)^\bot$,
with respect to the metric $\langle \cdot, \cdot \rangle$, are constant rank smooth subbundles of $E$.
It can be checked that $\mathcal{G}=\ker \rho/(\ker \rho)^\perp$ is a bundle of quadratic Lie algebras over $M$ and, as it was proved in \cite{csx1}, $E$ is isomorphic to $F^\ast \oplus \mathcal{G}\oplus F$. Precisely, for a given splitting $\lambda : F\to E$ of the short exact sequence
\[0 \to \ker \rho \to E \to F \to 0,\]
whose image $\lambda(F)$ is isotropic in $E$, there exists a unique splitting $\sigma_{\lambda} : \mathcal{G} \to \ker \rho$ of the short exact sequence
\[0 \to (\ker \rho)^{\perp} \to \ker \rho \to \mathcal{G} \to 0\]
with image $\sigma_{\lambda}(\mathcal{G})$ orthogonal to $\lambda(F)$ in $E$, and for the pair of splittings $(\lambda, \,\sigma_{\lambda})$, the map $\Psi_{\lambda} : F^\ast \oplus \mathcal{G}\oplus F \to E$ defined, for any $\xi \in \Gamma(F^{\ast})$, $r\in \Gamma(\mathcal{G})$, $x\in \Gamma(F)$, by
\begin{equation*}
\Psi_{\lambda}(\xi + r + x) = \frac{1}{2}({g^{\flat}}^{-1}\circ \rho^\ast)(\xi) + \sigma_{\lambda}(r) + \lambda(x),
\end{equation*}
is an isomorphism. In this case, $d_Ef = d_{F^\ast}f+0+0$, where $d_{F^\ast} : \C \to \Gamma(F^\ast)$ denotes the leafwise de Rham differential. The Courant algebroid structure on $F^{\ast}\oplus \mathcal{G}\oplus F$ is
completely determined by a linear $F$-connection $\lhd : \Gamma(F) \times \Gamma(\mathcal{G}) \to \Gamma(\mathcal{G})$ on $\mathcal{G}$, a bundle map $R : \bigwedge^2 F \to \mathcal{G}$, and a $3$-form $\mathcal{H} \in \Gamma(\bigwedge^3 F^{\ast})$ satisfying some compatibility conditions. In the following, we construct two examples of Dorfman connections in the framework of regular Courant algebroids.}}

\vspace{2mm}
\noindent
Example 1. {\rm{Consider the vector bundle of constant rank $B=F\oplus F^\ast$ endowed with the natural predual structure $(\lan \cdot, \cdot \ran, d_B)$ of $E\cong F^\ast \oplus \mathcal{G}\oplus F$. More precisely, $\lan \cdot, \cdot \ran : (F^\ast \oplus \mathcal{G}\oplus F) \times_M B \to \R$ is given, for any $\zeta + r + X \in \Gamma(F^\ast \oplus \mathcal{G}\oplus F)$ and $Y+\eta \in \Gamma(B)$, by
\begin{linenomath*} \begin{equation*}
\lan \zeta + r + X, Y+\eta \ran = \langle \eta, X\rangle + \langle \zeta, Y\rangle
\end{equation*} \end{linenomath*}
and $d_Bf = 0 + d_{F^\ast}f$. Let $pr_{F^\ast} : F\oplus F^\ast \to F^\ast$ be the projection onto the second summand and $\mathcal{Q} : \Gamma(F)\times \Gamma(\mathcal{G}) \to \Gamma(F^\ast)$ be the $\C$-bilinear map defined, for any $x, y\in \Gamma(F)$ and $r\in \Gamma(\mathcal{G})$, by
\begin{equation*}
\langle Q(x,r), \, y\rangle = \langle r,\, R(x,y)\rangle_{\mathcal{G}},
\end{equation*}
where $\langle \cdot, \cdot\rangle_{\mathcal{G}}$ denotes the nondegenerate, ad-invariant, pseudo-metric on the bundle of quadratic Lie algebras $\mathcal{G}$, \cite[Lemma 2.1]{csx1}. Choose a classical $F$-connection $\triangle$ on $F$ (there always exists one) and denote by $\triangle^\ast$ its dual connection on $F^\ast$. One can then check directly that the map
\[\nabla : \Gamma(F^\ast \oplus \mathcal{G}\oplus F) \times \Gamma(F\oplus F^\ast) \to \Gamma(F\oplus F^\ast)\]
defined, for $\zeta + r + X \in \Gamma(F^\ast \oplus \mathcal{G}\oplus F)$ and $Y+\eta \in \Gamma(F\oplus F^\ast)$, by
\begin{linenomath*} \begin{equation*}
\nabla_{\zeta + r + X}(Y+\eta) = ([X,Y]+\triangle_YX) + pr_{F^\ast}(\mathcal{L}_X\eta - i(Y)d\zeta) + \triangle^\ast_Y\zeta + \mathcal{Q}(Y,r)
\end{equation*} \end{linenomath*}
 is a $F^\ast \oplus \mathcal{G}\oplus F$ - Dorfman connection on $F\oplus F^\ast$.}}

\vspace{3mm}
\noindent
Example 2. {\rm{According to Proposition 4.12 in \cite{grtz-stn}, the bundle of a quadratic Lie algebras $\mathcal{G}$, endowed with the induced Courant algebroid structure
from the one of $F^\ast \oplus \mathcal{G}\oplus F$, is a Courant algebroid. More precisely, if $\iota : \mathcal{G} \to F^\ast \oplus \mathcal{G}\oplus F$ is the injection of $\mathcal{G}$ into $F^\ast \oplus \mathcal{G}\oplus F$ and $pr_{\mathcal{G}}: F^\ast \oplus \mathcal{G}\oplus F \to \mathcal{G}$ the projection of $F^\ast \oplus \mathcal{G}\oplus F$ on the second summand, the Dorfman bracket on $\Gamma(\mathcal{G})$ is given, for any $r_1, r_2 \in \Gamma(\mathcal{G})$, by $[r_1,r_2]_{\mathcal{G}} = pr_{\mathcal{G}}(\lcf \iota(r_1), \iota(r_2)\rcf)$, the anchor map by $\rho_{\mathcal{G}}=\rho \circ \iota = 0$, and the inner product by $\langle r_1, r_2\rangle_{\mathcal{G}} = \langle \iota(r_1),\iota(r_2) \rangle$. Consider the map}}
\begin{linenomath*} \begin{equation*}
\nabla : \Gamma(\mathcal{G}) \times \Gamma(F^\ast \oplus \mathcal{G}\oplus F) \to \Gamma(F^\ast \oplus \mathcal{G}\oplus F)
\end{equation*} \end{linenomath*}
{\rm{defined, for any $r\in \Gamma(\mathcal{G})$ and $\xi + s + X \in \Gamma(F^\ast \oplus \mathcal{G}\oplus F)$, by
\begin{linenomath*}
 \begin{equation} \label{connection G - F*+G+F}
\nabla_r(\xi + s + X) = \lcf \iota(r), \xi + s + X \rcf + 2\mathcal{Q}(X,r) - \lhd_Xr,
\end{equation} \end{linenomath*}
where $\mathcal{Q}$ is the map mentioned in the previous example and $\lhd : \Gamma(F) \times \Gamma(\mathcal{G}) \to \Gamma(\mathcal{G})$ is the linear $F$-connection on $\mathcal{G}$
provided by the Courant algebroid structure on $F^\ast \oplus \mathcal{G}\oplus F$, \cite{csx1}. Then, \eqref{connection G - F*+G+F} yields a $\mathcal{G}$-Dorfman connection on $F^\ast \oplus \mathcal{G}\oplus F$.}}
\end{examples}

\begin{proposition}[\cite{mjl}]\label{bott}
Let $(E,L)$ be a Manin pair and $L^0$ the annihilator of $L$ in $E^\ast$.
\begin{enumerate}
\item
The quotient $E/L$ is an $L$-Dorfman module with respect to the Dorfman action $\nabla^L : \Gamma(L) \times \Gamma(E/L) \to \Gamma(E/L)$ defined, for any $l\in \Gamma(L)$ and $\bar{e}\in \Gamma(E/L)$, by
\begin{linenomath*} \begin{equation*}
\nabla^L_l\bar{e} = \overline{\lcf l, e\rcf}.
\end{equation*} \end{linenomath*}
\item
The space $(E/L)^\ast \cong L^0$ is also an $L$-module relative to the dual connection $(\nabla^L)^\ast$ of $\nabla^L$.
\end{enumerate}
\end{proposition}
\begin{proof}
For the first item, note that, since $L$ is a Dirac subbundle of $E$, $(L, \lcf \cdot, \cdot \rcf)$ is a Lie algebroid and can be considered as a special case of a Courant algebroid. Also, since it is Lagrangian, $L = L^\bot$, the symmetric bilinear form of the Courant algebroid structure on $E$ induces a nondegenerate pairing
$\langle \cdot, \cdot \rangle : L \times_M E/L \to \R$. Further equip the vector bundle $E/L \to M$ with the map $d_{E/L} : \C \to \Gamma(E/L)$ given by $d_{E/L}(f) = \overline{d_Ef}$. Then $E/L$ is a predual bundle of $L$ in the sense of Definition \ref{pre-dual}. It is easy to check that $\nabla^L$ satisfies the axioms of a Dorfman connection. As an element of $\Gamma(\bigwedge^2 L^\ast)\otimes \Gamma(\mathrm{End}(E/L))$, the curvature $R^{\nabla^L}$
vanishes identically on the sections of $L$. For $l_1,l_2 \in \Gamma(L), \bar{e}\in \Gamma(E/L)$, we have
\begin{eqnarray*}
R^{\nabla^L}(l_1,l_2)\bar{e} & = & \nabla^L_{l_1}\nabla^L_{l_2}\bar{e} - \nabla^L_{l_2}\nabla^L_{l_1}\bar{e} - \nabla^L_{\lcf l_1, l_2\rcf}\bar{e} \nonumber \\
& = & \overline{\lcf l_1, \lcf l_2, e\rcf\rcf} - \overline{\lcf l_2, \lcf l_1, e\rcf\rcf} - \overline{\lcf \lcf l_1, l_2\rcf, e\rcf} \nonumber \\
& = & \overline{\lcf l_1, \lcf l_2, e\rcf\rcf -\lcf l_2, \lcf l_1, e\rcf\rcf -\lcf \lcf l_1, l_2\rcf, e\rcf} \nonumber \\
& \stackrel{(\ref{Jacobi-Courant})}{=}& 0.
\end{eqnarray*}
Also, for any $f\in \C$ such that $d_Ef \in \Gamma(L)$, i.e. for $f$ constant along the leaves of $\mathrm{Im}\rho(L)$, $\nabla^L_{d_Ef}\bar{e} = \overline{\lcf d_Ef, e\rcf}=0$.

\vspace{1mm}
\noindent
For the second item, let $(\nabla^L)^\ast : \Gamma(L) \times \Gamma((E/L)^\ast) \to \Gamma((E/L)^\ast)$ be the dual
connection of $\nabla^L$ and $R^{(\nabla^L)^\ast}$ be its curvature. By \eqref{dual Dorfman curvature} it follows that $(E/L)^\ast$ is an $L$-module.
\end{proof}

\vspace{2mm}
\noindent
In \cite{mjl} it was noted that the Dorfman connection $\nabla^L$ of the last Proposition is analogous to the Bott connection defined by an involutive subbundle of $TM$.
For this reason, it is named \emph{Bott-Dorfman connection associated to $L$}.

\begin{example}[Courant algebroid related to port-Hamiltonian systems]
{\rm{Port--Ha\-mil\-tonian systems are a generalization of Hamiltonian systems that aim to describe the dynamics of a Hamiltonian system in interaction with control units,
energy dissipa\-ting or energy storing units (ports) \cite{vdsch}. The state space of such a system is modeled by a manifold $M$ endowed with a Dirac structure $L$ in a
Courant algebroid \cite{merker}. More specifically, start with the standard Courant algebroid $TM\oplus T^\ast M$, a vector bundle $V$ over $M$ endowed with a
flat linear $TM$-connection $\triangle$, and its dual bundle $V^\ast$ equipped with the dual connection $\triangle^\ast$. The sections $\lambda_{out}+\lambda_{in}$ of $V\oplus V^\ast$
model the output and input of the port. The vector bundle $E = TM\oplus T^\ast M \oplus V \oplus V^\ast$ equipped with the projection $\rho : E\to TM$ as anchor map, the symmetric nondegenerate bilinear form}}
\begin{linenomath*} \begin{equation*}
\langle X+ \zeta + \lambda_{out} + \lambda_{in}, Y+\eta + \mu_{out}+ \mu_{in}\rangle = \langle \eta, X\rangle + \langle \zeta, Y\rangle + \langle \mu_{in}, \lambda_{out}\rangle + \langle \lambda_{in},\mu_{out}\rangle,
\end{equation*} \end{linenomath*}
{\rm{and the bracket}}
\begin{eqnarray*}
\lefteqn{\lcf X+ \zeta+ \lambda_{out}+\lambda_{in}, Y+ \eta + \mu_{out}+ \mu_{in}\rcf  =  }\nonumber \\
& & [X,Y] + (\mathcal{L}_X\eta - i(Y)d\zeta + \langle \triangle^\ast_{\cdot}\lambda_{in}, \mu_{out}\rangle + \langle \mu_{in}, \triangle_{\cdot}\lambda_{out}\rangle) \nonumber \\
& & + (\triangle_X\mu_{out} - \triangle_Y\lambda_{out}) + (\triangle^\ast_X\mu_{in}-\triangle^\ast_Y\lambda_{in})
\end{eqnarray*}
{\rm{is a Courant algebroid over $M$. The dynamics of the system are determined by a Hamiltonian $H$ via the Hamiltonian condition $\dot{x}+ dH + \lambda_{out}+\lambda_{in}\in \Gamma(L)$.}}

\vspace{1mm}
\noindent
{\rm{Consider the vector bundle $B = T^\ast M \oplus V$. Clearly, the natural coupling of $T^\ast M$ with $TM$ and of $V$ with $V^\ast$
gives a coupling $\lan \cdot, \cdot \ran$ between $E$ and $B$. Then the pair $(\lan \cdot, \cdot \ran, \, d_B)$, where $d_B : \C \to \Gamma(B)$ is the derivation $d_Bf = (df, 0)$,
defines a predual structure of $E$ on $B$ and the map $\nabla : \Gamma(E) \times \Gamma(B) \to \Gamma(B)$ defined,
for any $e = X+\zeta + \lambda_{out}+ \lambda_{in}\in \Gamma(E)$ and $b = \eta + \mu_{out} \in \Gamma(B)$, by}}
\begin{linenomath*} \begin{equation*}
\nabla_eb = \mathcal{L}_X\eta  + \langle \triangle^\ast_{\cdot}\lambda_{in}, \mu_{out}\rangle + \triangle_X\mu_{out},
\end{equation*} \end{linenomath*}
{\rm{establishes an $E$-Dorfman connection on $B$. It is the restriction of $\lcf \cdot, \cdot \rcf$ on $\Gamma(E) \times \Gamma(B)$ taking values in $\Gamma(B)$.}}

\vspace{1mm}
\noindent
{\rm{We arrive at a similar result, if we consider $B'= T^\ast M\oplus V^\ast$, which evidently is a predual  of $E$, and take $\nabla' : \Gamma(E) \times \Gamma(B') \to \Gamma(B')$ to be
the restriction of $\lcf \cdot, \cdot \rcf$ on $\Gamma(E) \times \Gamma(B')$ taking values in $\Gamma(B')$. More precisely, for any $e = X + \zeta + \lambda_{out}+ \lambda_{in}\in \Gamma(E)$
and $b' = \eta+ \mu_{in} \in \Gamma(B')$,}}
\begin{linenomath*} \begin{equation*}
\nabla'_eb' = (\mathcal{L}_X\eta + \langle \mu_{in}, \triangle_{\cdot}\lambda_{out}\rangle) + \triangle^\ast_X\mu_{in}.
\end{equation*} \end{linenomath*}
\end{example}

\vspace{5mm}
\noindent
\emph{E-mails}
\\
\\
\noindent
Panagiotis Batakidis: \emph{batakidis@math.auth.gr}
\\
\\
\noindent
Fani Petalidou: \emph{petalido@math.auth.gr}

\end{document}